\documentclass[11pt]{article}
\usepackage{amsfonts}
\usepackage{amsgen,amsmath,amstext,amsbsy,amsopn,amssymb,subfigure}
\usepackage[dvips]{graphicx}
\usepackage{comment}
\usepackage{enumitem}
\usepackage{booktabs} 
\usepackage{blkarray}
\usepackage{algorithm}
\usepackage{algpseudocode}
\usepackage{url}
\usepackage{longtable}

\usepackage[usenames]{color}

 \allowdisplaybreaks

\newtheorem{Definition}{Definition}
\newtheorem{Theorem}{Theorem}
\newtheorem{Assumption}{Assumption}

\newtheorem{Lemma}{Lemma}
\newtheorem{Remark}{Remark}

\newcommand{\be}{\begin{equation}}
\newcommand{\ee}{\end{equation}}
\newcommand{\bea}{\begin{eqnarray}}
\newcommand{\eea}{\end{eqnarray}}
\newcommand{\beas}{\begin{eqnarray*}}
\newcommand{\eeas}{\end{eqnarray*}}

\newcommand{\bbR}{\mathbb{R}}

\newcommand{\Var}{{\rm Var}}

\newcommand{\sgn}{{\rm sgn}}

\newcommand{\bbP}{\mathbb{P}}

\newcommand{\argmin}{\mathop{\rm arg\min}}
\newcommand{\argmax}{\mathop{\rm arg\max}}

\newcommand{\supp}{{\rm supp}}

\allowdisplaybreaks

\begin{document}

\title{Sparse Group Lasso: Optimal Sample Complexity, Convergence Rate, and Statistical Inference}
\author{T. Tony Cai$^{1}$, ~ Anru R. Zhang$^{2, 4}$, ~ and ~ Yuchen Zhou$^{3, 4}$}
\date{}
\maketitle

\begin{abstract}
	We study sparse group Lasso for high-dimensional double sparse linear regression, where the parameter of interest is simultaneously element-wise and group-wise sparse. This problem is an important instance of the simultaneously structured model -- an actively studied topic in statistics and machine learning.  In the noiseless case, matching upper and lower bounds on sample complexity are established for the exact recovery of sparse vectors and for stable estimation of approximately sparse vectors, respectively. In the noisy case,  upper and matching minimax lower bounds for estimation error are obtained. We also consider the debiased sparse group Lasso and investigate its asymptotic property for the purpose of statistical inference. Finally, numerical studies are provided to support the theoretical results.
\end{abstract}


\footnotetext[1]{Department of Statistics \& Data Science, The Wharton School, University of Pennsylvania, Philadelphia, PA 19104. The research of Tony Cai was supported in part by NSF grants DMS-1712735 and DMS-2015259  and NIH grants R01-GM129781 and R01-GM123056.}  
\footnotetext[2]{Departments of Biostatistics \& Bioinformatics, Computer Science, Mathematics, and Statistical Science, Duke University, Durham, NC 27710.}
\footnotetext[3]{Department of Statistics \& Data Science, The Wharton School, University of Pennsylvania, Philadelphia, PA 19104.}
\footnotetext[4]{The research of Anru R. Zhang and Yuchen Zhou was supported in part by NSF Grants  CAREER-2203741 and DMS-1811868 and NIH grant R01-GM131399-01. This work was done while Anru R. Zhang and Yuchen Zhou were at the Department of Statistics, University of Wisconsin-Madison, Madison, WI 53706.}

\newpage
%%%%%%%%%%%%%%%%%%%%%%%%%%%%%
\section{Introduction}\label{sec:intro}
%%%%%%%%%%%%%%%%%%%%%%%%%%%%%

Consider the \emph{high-dimensional double sparse regression} with simultaneously group-wise and element-wise sparsity structures
\begin{equation}\label{eq:model}
y = X \beta^* + \varepsilon, \quad \text{or equivalently}\quad y_i = X_i^\top \beta^* + \varepsilon_i, \quad i=1,\ldots, n.
\end{equation}
Here, the covariates $X \in \mathbb{R}^{n\times p}$ and parameter $\beta^*$ are divided into $d$ known groups, where the $j$th group contains $b_j$ variables,
\begin{equation}\label{eq:X-group}
X = [X_{(1)} ~ \cdots ~ X_{(d)}],\quad \beta^*=\left((\beta^*_{(1)})^\top, \cdots (\beta^*_{(d)})^\top\right)^\top, \quad X_{(j)} \in \mathbb{R}^{n\times b_j}, \beta_{(j)}^*\in \mathbb{R}^{b_j};
\end{equation}
$\beta^*$ is a \emph{$(s, s_g)$-sparse vector} in the sense that
\begin{equation}\label{eq:sparsity}
\|\beta^*\|_{0, 2} := \sum_{j=1}^{d} 1_{\{\beta^*_{(j)} \neq 0\}}\leq s_g\quad \text{and}\quad \|\beta^*\|_0 := \sum_{i=1}^p 1_{\{\beta^*_i \neq 0\}} \leq s,
\end{equation}
The focus of this paper is on the estimation of and inference for $\beta^*$ based on $(y, X)$. This problem has great importance in a variety of applications. For example in genome-wide association studies (GWAS) \cite{silver2013pathways}, the genes can be grouped into pathways and it is believed that only a small portion of the pathways contain causal single nucleotide polymorphisms (SNPs), and the number of causal SNPs is much less than the one of non-causal SNPs in a causal pathway. The sparse group Lasso has been applied to identify causal genes or SNPs associated with a certain trait \cite{silver2013pathways}.
Other examples include cancer diagnosis and therapy \cite{vidyasagar2014machine,allahyar2015feral}, classification \cite{rao2015classification}, and climate prediction \cite{chatterjee2012sparse} among many others. The problem can also be viewed as a prototype of various problems in statistics and machine learning, such as the sparse multiple response regression \cite{wang2013block} and multiple task learning \cite{lounici2009taking,lozano2012multi,zhou2017can}. 

The sparse group Lasso \cite{friedman2010note,simon2013sparse,li2015multivariate} provides a classic and straightforward estimator for $\beta^*$:
\begin{equation}\label{eq:sparse-group-lasso}
\hat{\beta} = \argmin_{\beta} \|y - X\beta\|_2^2 + \lambda\|\beta\|_1 + \lambda_g\|\beta\|_{1,2}.
\end{equation}
Here, $\|\beta\|_1 = \sum_{i=1}^p|\beta_i|$ and $\|\beta\|_{1,2} = \sum_{j}\|\beta_{(j)}\|_2$ are $\ell_1$ and $\ell_{1,2}$ convex regularizers to account for element-wise and group-wise sparsity structures, respectively. $\lambda, \lambda_g\geq 0$ are tuning parameters. In the noiseless setting that $\varepsilon=0$, one can apply the constrained $\ell_1+\ell_{1,2}$ minimization instead to estimate $\beta^*$:
\begin{equation}\label{eq:ell_1_minimization}
\begin{split}
\hat{\beta}\quad  =\quad & \argmin\quad \lambda\|\beta\|_1 + \lambda_g \|\beta\|_{1, 2}\\
& \text{subject to} \quad y = X\beta.
\end{split}
\end{equation}
In fact, when $\lambda, \lambda_g$ tend to zero while $\lambda/\lambda_g$ is fixed as a constant, the sparse group Lasso \eqref{eq:sparse-group-lasso} tends to the $\ell_1+\ell_{1,2}$ minimization \eqref{eq:ell_1_minimization}.

When $\beta^*$ is only element-wise sparse, the regular Lasso \cite{tibshirani1996regression}
\begin{equation}\label{eq:Lasso}
\hat{\beta}^{L} = \argmin_{\beta} \|y - X\beta\|_2^2 + \lambda\|\beta\|_1
\end{equation}
can be applied and its theoretical properties have been well studied. See, for example, \cite{bickel2009simultaneous,verzelen2012minimax}. When $\beta^*$ is only group-wise sparse, the group Lasso 
\begin{equation}\label{eq:group-Lasso}
\hat{\beta}^{GL} = \argmin_{\beta} \|y - X\beta\|_2^2 + \lambda_g\|\beta\|_{1,2}
\end{equation}
and its variations have been widely investigated \cite{yuan2006model,lounici2011oracle,bunea2013group}. However, to estimate the simultaneously element-wise and group-wise sparse vector $\beta^*$, despite many empirical successes of sparse group Lasso in practice, the theoretical properties, including optimal rate of convergence and sample complexity, are still unclear so far to the best of our knowledge. 

%%%%%%%%%%%%%%%%%%%
\subsection{Simultaneously Structured Models}
%%%%%%%%%%%%%%%%%%%

More broadly speaking, the simultaneously structured models, i.e., the parameter of interest has multiple structures at the same time, have attracted enormous attention in many fields including statistics, applied mathematics, and machine learning. In addition to the high-dimensional double sparse regression, other simultaneously structured models include sparse principal component analysis \cite{johnstone2009consistency,ma2013sparse}, tensor singular value decomposition \cite{zhang2018tensor,wang2018learning}, simultaneously sparse and low-rank matrix/tensor recovery \cite{oymak2015simultaneously,hao2018sparse}, sparse matrix/tensor SVD \cite{zhang2018optimal}, and sparse phase retrieval \cite{jaganathan2013sparse,shechtman2014gespar,Cai2016PhaseRetrieval}. As shown in \cite{oymak2013noisy,oymak2015simultaneously}, by minimizing multi-objective regularizers with norms associated with these structures (such as $\ell_1$ norm for element-wise sparsity, nuclear norm for low-rankness, and total variation norm for piecewise constant structures), one usually cannot do better than applying an algorithm that only exploits one structure. They particularly illustrated that simultaneously sparse and low-rank structured matrix cannot be well estimated by penalizing $\ell_1$ and nuclear norm regularizers. Instead, non-convex methods were proposed and shown to achieve better performance. 

However based on their results, it remains an open question whether the convex regularization, such as sparse group Lasso or $\ell_1+\ell_{1,2}$ minimization, can achieve good performance in estimation of parameter with two types of sparsity structures, such as the aforementioned high-dimensional double sparse regression. Specifically, as illustrated in Section \ref{sec:noiseless}, a direct application of \cite{oymak2015simultaneously} does not provide a sample complexity lower bound for exact recovery that matches our upper bound.

\subsection{Optimality and Related Literature}

This paper fills the void of statistical limits of sparse group Lasso and provides an affirmative answer to the aforementioned question: by exploiting both element-wise and group-wise sparsity structures, the $\ell_1+\ell_{1,2}$ regularization does provide better performance in high-dimensional double sparse regression. Particularly in the noiseless case, it is shown that $(s, s_g)$-sparse vectors can be exactly recovered and approximately $(s, s_g)$-sparse vectors can be stably estimated with high probability whenever the sample size satisfies $n\gtrsim s_g\log(d/s_g) + s\log(es_gb)$, where $b = \max_{1 \leq i \leq d}b_i$. On the other hand, we prove that exact recovery cannot be achieved by $\ell_1+\ell_{1,2}$ regularization and stable estimation of approximately $(s, s_g)$-sparse vectors is impossible in general unless $n\gtrsim s_g\log(d/s_g) + s\log(es_gb/s)$. We then consider the noisy case and develop the matching upper and lower bounds on the convergence rate for the estimation error. Simulation studies are carried out and the results support our theoretical findings. In addition, statistical inference for the individual coordinates of $\beta^*$ is studied. A confidence interval is constructed based on the debiased sparse group Lasso estimator and its asymptotic property. The results show that by exploring the simultaneously element-wise and group-wise sparsity structures, the debiased sparse group Lasso requires less sample size than the debiased Lasso and debiased group Lasso in the literature \cite{zhang2014confidence,javanmard2014confidence,mitra2016benefit,cai2017confidence}.

The theoretical analysis of sparse group Lasso and $\ell_1+\ell_{1,2}$ minimization is highly non-trivial. First, the regularizer $\lambda\|\cdot\|_1 + \lambda_g\|\cdot\|_{1,2}$ is not decomposable with respect to the support of $\beta^*$ so that the classic techniques of decomposable regularizers \cite{negahban2012unified} and null space property \cite{stojnic2008compressed} may not be suitable here. Despite a substantial body of literature on high-dimensional element-wise sparse vector estimation based on restricted isometry property (RIP) \cite{candes2006robust,candes2007dantzig,cai2013compressed,cai2013sharp,cai2014sparse} and restricted eigenvalue \cite{bickel2009simultaneous}, these techniques cannot provide nearly optimal results for sparse group Lasso here as it is technically difficult to partition general vectors into simultaneously element-wise and group-wise ones that preserves some ordering structures. Departing from the previous literature, our theoretical analysis relies on a novel construction of approximate dual certificate. See Section \ref{sec:proof-sketch} for further details.
Although our results mostly focus on the performance of sparse group Lasso and $\ell_1+\ell_{1,2}$ estimators, the techniques of approximate dual certificate on multi-norm structures here can also be of independent interest.

The statistical properties of sparse group Lasso and related estimators have been studied previously. For example, \cite{chatterjee2012sparse} developed consistency results for estimators with a general tree-structured norm regularizers, of which the sparse group Lasso is a special case. \cite{poignard2018asymptotic} analyzed the asymptotic behaviors of the adaptive sparse group Lasso estimator. \cite{rao2015classification,rao2013sparse} studied the multi-task learning and classification problems based on a variant of sparse group Lasso estimator. \cite{li2015multivariate} studied multivariate linear regression via sparse group Lasso. \cite{ahsen2017error} provided a theoretical framework for developing error bounds of the group Lasso, sparse group Lasso, and group Lasso with tree structured overlapping groups. Specifically, their results imply that the group-wise sparse signal can be exactly recovered with high probability by solving \eqref{eq:ell_1_minimization} if the sample size satisfies $n \gtrsim s_g\left(b + \log d\right)$. 
Different from previous results, this paper focused on both the required sample size and convergence rate of estimation error of sparse group Lasso. To the best of our knowledge, this is the first paper that provides optimal theoretical guarantees for both the sample complexity and estimation error of sparse group Lasso.

\subsection{Organization of the Paper}

The rest of the article is organized as follows. After a brief introduction to notation and preliminaries in Section \ref{sec:notation}, the main theoretical results on constrained $\ell_1+\ell_{1,2}$ minimization in the noiseless setting is presented in Section \ref{sec:noiseless} and the key proof ideas are explained in Section \ref{sec:proof-sketch}. Results for sparse group Lasso in the noisy setting are discussed in Section \ref{sec: noisy}. In particular, the optimal rate of estimation error and statistical inference are studied in Sections \ref{subsec:rate} and \ref{subsec:inference}, respectively. In Section \ref{sec:tuning}, we introduce a practical scheme to select tuning parameters. In Section \ref{sec:numerical}, we provide simulation results in both noiseless and noisy cases to justify our theoretical findings. The proofs of technical results are given in Section \ref{sec:proofs}. All technical lemmas and their proofs can be found in Appendix \ref{sec:technical lemma}.

%%%%%%%%%%%%%%%%%%%%%%%%%%%%%%%%%
\section{$\ell_1+\ell_{1,2}$ Minimization in Noiseless Case}\label{sec:theory}
%%%%%%%%%%%%%%%%%%%%%%%%%%%%%%%%%

\subsection{Notation and Preliminaries}\label{sec:notation}
The following notation will be used throughout the paper. We denote $a\wedge b = \min\{a, b\}, a\vee b = \max\{a, b\}$. Let $\sgn(\cdot)$ be the sign function, i.e., $\sgn(x) = 1, 0, $ or $-1$, if $x>0, x=0$, or $x<0$, respectively. $H_{\alpha}(\cdot)$ is the soft-thresholding function such that $H_{\alpha}(x) = \sgn(x)\cdot \left\{\left(|x|-\alpha\right)\vee 0\right\}$ for any $x\in \mathbb{R}$. We say $a \lesssim b$ and $a \gtrsim b$ if $a \leq Cb$ and $b \leq Ca$ for some uniform constant $C > 0$, respectively. $a \asymp b$ means $a \lesssim b$ and $a \gtrsim b$ both hold. Let the uppercase $C, C_1, C_0, \ldots$ and lowercase $c, c_1, c_0,\ldots$ denote large and small positive constants respectively, whose actual values vary from time to time.  Throughout the paper, we focus on the parameter index set $\{1,\ldots, p\}$ partitioned into $d$ groups. Denote $(1), \ldots, (d) \subseteq \{1,\ldots, p\}$ as the index sets belonging to each group. Additionally, for any group index subset $G\subseteq \{1,\ldots, d\}$, define $(G) = \cup_{j\in G}(j)$, $(G^c) = \cup_{j\notin G}(j)$. For any vector $\gamma$ and index subset $T$, $\gamma_T \in \mathbb{R}^{|T|}$ represents the sub-vector of $\gamma$ with index set $T$. In particular, $\gamma_{(G)}$ represents the sub-vector of $\gamma$ in the union of Groups $j\in G$. Define the $\ell_q$ norm of any vector $\gamma$ as $\|\gamma\|_q = \left(\sum_i|\gamma_i|^q\right)^{1/q}$. For any vector $\beta \in \mathbb{R}^p$ with group structures, we also define the $\ell_{q_1, q_2}$ norm for any $0\leq q_1, q_2\leq \infty$ as
\begin{equation*}
\|\gamma\|_{q_1, q_2} = \left(\sum_{j=1}^d \|\gamma_{(j)}\|_{q_2}^{q_1}\right)^{1/q_1} = \left\{\sum_{j=1}^d \left(\sum_{i \in (j)}|\gamma_{i}|^{q_2}\right)^{q_1/q_2}\right\}^{1/q_1}.
\end{equation*}
In particular, $\|\gamma\|_{0, 2} = \sum_{j=1}^d 1_{\{\gamma_{(j)}\neq 0\}}$ is the number of non-zero groups of $\gamma$,  $\|\gamma\|_{\infty, 2} = \max_j \|\gamma_{(j)}\|_2$ is the maximum $\ell_2$ norm among all groups of $\gamma$, and $\|\gamma\|_{1, 2} = \sum_{j=1}^d \|\gamma_{(j)}\|_2$ is the group-wise $\ell_1$ penalty. With a slight abuse of notation, we simply denote $\|\gamma_T\|_{q_1, q_2}=\|u\|_{q_1,q_2}$ if $u \in \bbR^p$, $u$ restricted on subset $T$ is $\gamma_T$ and $u$ restricted on $T^c$ is $0$.

The focus of this paper is on simultaneously element-wise and group-wise sparse vectors defined as follows.
\begin{Definition}[Simultaneous element-wise and group-wise sparsity]\label{as:sparsity}
	Assume $\beta^*\in \mathbb{R}^p$ is associated with group partition $(1),\ldots, (d)$. For positive integers $s, s_g$ satisfying $s_g \leq d$ and $s_g \leq s \leq \max_{\Omega \subseteq \{1, \dots, d\}, |\Omega| = s_g}\sum_{i \in \Omega}b_{i}$, we say $\beta^*$ is $(s, s_g)$-sparse if
	$$\|\beta^*\|_{0, 2} = \sum_{j=1}^d 1_{\{\beta^*_{(j)}\neq 0\}}\leq s_g, \quad \|\beta^*\|_0 = \sum_i 1_{\{\beta_i^*\neq 0\}}\leq s.$$ 
\end{Definition}

%%%%%%%%%%%%%%%%%%%%%%%%%%%%
\subsection{Noiseless Case and Sample Complexity}\label{sec:noiseless}
%%%%%%%%%%%%%%%%%%%%%%%%%%%%

To analyze the performance of sparse group Lasso and $\ell_1+\ell_{1,2}$ minimization, we first introduce the following assumption on the design matrix $X$. 
\begin{Assumption}[Sub-Gaussian assumption]\label{as:design}
	Suppose all rows of $X$ are i.i.d. centered sub-Gaussian distributed. Specifically, $\mathbb{E}X_{i\cdot}=0, \Var(X_{i\cdot}^\top) = \Sigma$, and for any $\alpha \in \bbR^p$, we have $\mathbb{E}\exp\left(\alpha^\top \Sigma^{-1/2}X_{i\cdot}^\top\right) \leq \exp\left(\kappa^2 \|\alpha\|_2^2/2\right)$ for constant $\kappa>0$. We also assume there exist two constants $C_{\max} \geq c_{\min} > 0$ such that $c_{\min} \leq \sigma_{\min}(\Sigma) \leq \sigma_{\max}(\Sigma) \leq C_{\max}$, where $\sigma_{\max}(\Sigma)$ and $\sigma_{\min}(\Sigma)$ are the largest and smallest eigenvalues of $\Sigma$, respectively. 
\end{Assumption}
Clear, a random matrix $X$ with i.i.d. standard normal entries satisfies this assumption -- this design is referred to as the Gaussian ensemble and has been considered as a benchmark setting in compressed sensing and high-dimensional regression literature \cite{candes2011probabilistic,javanmard2018debiasing}. 

The following theorem shows that the $\ell_1 + \ell_{1, 2}$ minimization achieves the exact recovery with high probability when $\beta^*$ is simultaneously element-wise and group-wise sparse, $X$ is weakly dependent, and Assumption \ref{as:design} holds. The theorem also provides a more general upper bound on estimation error if $\beta^*$ is approximately element-wise and group-wise sparse.
\begin{Theorem}[$\ell_1+\ell_{1,2}$ minimization in noiseless case]\label{th:noiseless}
	Suppose one observes $y = X \beta^*$, where $X$ has the group structure \eqref{eq:X-group} and satisfies Assumption \ref{as:design}, $\beta^*$ is $(s, s_g)$-sparse, and $b = \max_{1 \leq i \leq d}b_i$. Let $T$ be the support of $\beta^*$. Suppose there exist uniform constants $C, c>0$ such that 
	\begin{equation}\label{ineq:sample-complexity}
	n \geq C\left(s_g \log(d/s_g) + s\log(es_gb)\right),
	\end{equation} 
	\begin{equation}\label{ineq:irrepresentable}
	\quad \max_{i\in T^c} \left\|\Sigma_{i, T}\Sigma_{T, T}^{-1} \right\|_2 \leq c/\sqrt{s},
	\end{equation} 
	then the constrained $\ell_1 + \ell_{1,2}$ minimization \eqref{eq:ell_1_minimization} with $\lambda_g = \sqrt{s/s_g}\lambda$ achieves the exact recovery with probability at least $1 - C\exp(-cn/s)$. 
	
	Moreover, if $\beta^*\in \mathbb{R}^p$ is a general vector and $\hat{\beta}$ is the solution to the constrained $\ell_1 + \ell_{1, 2}$ minimization \eqref{eq:ell_1_minimization} with $\lambda_g = \sqrt{s/s_g}\lambda$, then
	\begin{equation}\label{ineq:noiseless_approx}
		\|\hat{\beta} - \beta^*\|_2 \lesssim \min_{S: \substack{\|\beta^*_S\|_0 \leq s, \|\beta^*_S\|_{0, 2} \leq s_g, \\\max_{i \in S^c}\|\Sigma_{i, S}\Sigma_{S, S}^{-1}\|_2 \leq c/\sqrt{s}}}\left(\frac{1}{\sqrt{s}}\|\beta^*_{S^c}\|_1 + \frac{1}{\sqrt{s_g}}\|\beta^*_{S^c}\|_{1, 2}\right).
	\end{equation}
	 with probability at least $1 - C\exp(-cn/s)$. 
\end{Theorem}
\begin{Remark}[Interpretation and comparison]\rm 
	In Theorem \ref{th:noiseless}, the required sample size for achieving exact recovery contains two terms: $s_g\log(d/s_g)$ and $s\log(es_gb)$. Intuitively speaking, $s_g\log(d/s_g)$ corresponds to the complexity of identifying $s_g$ non-zero groups and $s\log(es_gb)$ corresponds to the complexity of estimating $s$ non-zero elements of $\beta$ in $s_g$ known groups.
	
	When $\beta^*$ is only element-wise or group-wise sparse, one can apply respectively the classic $\ell_1$ or $\ell_{1,2}$ minimization to recover $\beta^*$,
	\begin{equation}\label{eq:l_1-mini}
	\hat{\beta}^{\ell_1} = \argmin_{\beta} \|\beta\|_1 \quad \text{subject to}\quad y = X\beta,
	\end{equation}
	\begin{equation}\label{eq:l_12-mini}
	\hat{\beta}^{\ell_{1,2}} = \argmin_{\beta} \|\beta\|_{1,2} \quad \text{subject to}\quad y = X\beta.
	\end{equation}
	The $\ell_1$ minimization and $\ell_{1,2}$ minimization here are respectively the special form of the regular Lasso and group Lasso (if $\lambda, \lambda_g = 0_+$ in \eqref{eq:Lasso} and \eqref{eq:group-Lasso}), respectively. Especially if the group size $b_1 \asymp \cdots \asymp b_d \asymp b$, to ensure exact recovery in the noiseless setting with high probability, \eqref{eq:l_1-mini} requires $n\gtrsim Cs\log(ebd/s)$ \cite{foucart2013mathematical} and group Lasso requires $n\gtrsim s_g (b + \log(ed/s_g))$. The $\ell_1+\ell_{1, 2}$ minimization \eqref{eq:ell_1_minimization} has provable advantages over both regular and group Lasso when $b \gg \log(d) \gg \log(es_gb)$ and $s_gb/\log(es_gb) \gg s \gg s_g$. In particular,  when $s_g = s$, the double sparse regression reduces to the vanilla sparse linear regression, and the upper bound \eqref{ineq:noiseless_approx} matches the classic upper bound for $\ell_1$ minimization \cite{candes2011probabilistic}.

In addition, Condition \eqref{ineq:irrepresentable} is an important technical condition we used in our theoretical analysis. 
\end{Remark}

Next, we consider the sample complexity lower bound. Suppose $b_1 = b_2 = \cdots = b_d$ and $d \geq 2s_g$. Recall that one observes $y = X\beta^*$ without noise and aims to estimate the $(s, s_g)$-sparse vector $\beta^*$ based on $y$ and $X$. As indicated by classic results in compressed sensing \cite{candes2005decoding}, with sufficient computing power, the $\ell_0$ minimization below achieves exact recovery of $\beta^*$
\begin{equation}\label{eq:l_0-minimization}
\hat{\beta}^{\ell_0} = \argmin \|\beta\|_0 \quad \text{subject to}\quad X\beta = y
\end{equation}
as along as $X$ is non-degenerate and $n\geq 2s$. This bound is actually sharp: when $n < 2s$, for any set $T \subseteq \{1, \dots, db\}$ with cardinality $2s$, one can find a vector $\gamma$ such that $\supp(\gamma) \subseteq T$ and $X\gamma = 0$. By choosing an appropriate $T$, we can split the support $\gamma$ to obtain two $(s, s_g)$-sparse vectors $\beta_1, \beta_2$ satisfying $\beta_1 + \beta_2 = \gamma$. Then, $X\beta_1 = X(-\beta_2)$ but there is no way to distinguish $\beta_1$ and $\beta_2$ merely based on $X$ and $y = X\beta_1 = X(-\beta_2)$.

However, the $\ell_0$ minimization \eqref{eq:l_0-minimization} is computational infeasible in practice while a larger sample size is required for applying more practical methods. The following theorem shows that by performing the convex $\ell_1$ regularization, $\ell_{1,2}$ regularization, or any weighted combination of them, one requires at least   $\Omega(s_g\log(d/s_g) + s\log(es_gb/s))$ observations to ensure exact recovery of $(s, s_g)$-sparse vectors.
\begin{Theorem}[Sample complexity lower bound for exact recovery]\label{th:sample_complexity}
	Suppose $b_1 = \cdots = b_d = b$, $d, b \geq 3$. Suppose $X$ is an $n$-by-$(db)$ matrix. If every $(2s, 2s_g)$-sparse vector $\beta \in \bbR^{db}$ is a minimizer of the following programming for some $(\lambda, \lambda_g) \in \{(\lambda, \lambda_g): \lambda, \lambda_g \geq 0, \lambda + \lambda_g > 0\}$: 
	$$\min_z \lambda\|z\|_1 + \lambda_g\|z\|_{1,2} \quad \text{subject to}\quad  Xz = y = X\beta.$$
	In other words, if the $\ell_1+\ell_{1,2}$ minimization exactly recover all $(2s, 2s_g)$-sparse vector $\beta$, then we must have $n \gtrsim s_g\log(d/s_g) + s\log(es_gb/s)$.
\end{Theorem}

The following sample complexity lower bound shows that for arbitrary methods, to ensure stable estimation of all approximately sparse vectors, one requires at least $\Omega(s_g\log(d/s_g) + s\log(es_gb/s))$ observations.
\begin{Theorem}[Sample complexity lower bound for stable estimation]\label{th:lower_noiseless}
	Suppose $b_1 = \cdots = b_d = b$, $b,d \geq 3$. Assume there exists a matrix $X \in \bbR^{n \times (bd)}$, a map $\Delta: \bbR^n \to \bbR^{bd}$ ($\Delta$ may depend on $X$), and a constant $C>0$ satisfying 	
	\begin{equation}\label{eq47}
		\|\beta - \Delta(X\beta)\|_2 \leq C\left(\frac{\|\beta\|_1}{\sqrt{s}} + \frac{\|\beta\|_{1,2}}{\sqrt{s_g}}\right)
	\end{equation}
	for all $\beta \in \bbR^p$ and some $s, s_g$ satisfying $d\geq s_g, s_gb\geq s\geq s_g$. There exists constants $C_0$ and $c_0$ that depend only on $C$ such that whenever $s_g \ge C_0$, we must have
	\begin{equation*}
		n \geq c_0(s_g\log(d/s_g) + s\log(es_gb/s)).
	\end{equation*}
\end{Theorem}

\begin{Remark}[Optimality and comparison with previous results]\rm 
	Theorems \ref{th:sample_complexity} and \ref{th:lower_noiseless} show that the sample complexity upper bound in Theorem \ref{th:noiseless} is rate-optimal under a weak condition: $\log(es_gb) \asymp \log(es_gb) - \log(s)$ or $\log(d) \geq 2s\log(s)/s_g$. Oymak, et al. \cite{oymak2015simultaneously} provided a general analysis for convex regularization of simultaneously structured parameter estimation. Specifically for the high-dimensional double sparse regression, a direct application of their Theorem 3.2 and Corollary 3.1 implies that if $\ell_1+\ell_{1,2}$ minimization can exactly recover $(s, s_g)$-sparse vector $\beta^*$ with a constant probability, one must have $n\gtrsim s$. We can see that Theorem \ref{th:sample_complexity} provides a sharper lower bound on sample complexity.
	
	In addition, by setting $s_g = s$, the lower bound in Theorems \ref{th:sample_complexity} and \ref{th:lower_noiseless} reduces to $n \gtrsim s\log(p/s)$, which matches the optimal sample complexity lower bound for exact recovery of $s$-sparse vectors \cite[Theorem10.11, Proposition 10.7]{foucart2013mathematical}. By setting $s = s_gb$, we obtain a sample complexity lower bound $n \gtrsim s_g(b + \log(d/s_g))$ for (approximate) $s_g$-group-wise sparse vector recovery and stable estimation. To the best of our knowledge, this is the first sample complexity lower bound for group Lasso.
\end{Remark}

\subsection{Proof Sketches}\label{sec:proof-sketch}

We briefly discuss the proof sketches of the main technical results in this section. The detailed proofs are postponed to Section \ref{sec:proofs}.

The proof of Theorem \ref{th:noiseless} is based on a novel dual certificate scheme. The dual certificate \cite{bertsekas2003convex} has been used in the theoretical analysis for various convex optimization methods in high-dimensional problems, such as matrix completion \cite{candes2009exact,gross2011recovering}, compressed sensing \cite{candes2011probabilistic}, robust PCA \cite{candes2011robust}, tensor completion \cite{yuan2016tensor}, etc. The high-dimensional double sparse linear regression exhibits different aspects from these previous works due to the simultaneous sparsity structure. In particular, we can show that if the $u_{et}$ defined below is in the row space of $X$, it can be used as an exact dual certificate for recovery of $(s, s_g)$-sparse vector $\beta^*$: 
\begin{equation}\label{eq:exact-dual-certificate}
u_{et} = v_{et} + w_{et}\in \mathbb{R}^p, \quad  \left\{\begin{array}{ll}
(v_{et})_{(j)} =\sqrt{s/s_g}\beta^*_{(j)}/\|\beta^*_{(j)}\|_2, &  j\in G;\\
\|(v_{et})_{(j)}\|_2< \sqrt{s/s_g}, & j \in G^c; \end{array}\right. \quad 
\left\{\begin{array}{ll}
(w_{et})_T = \sgn(\beta^*_T)\\
\|(w_{et})_{T^c}\|_\infty < 1.
\end{array}\right.
\end{equation}
Here, $T$ and $G$ are the element-wise and group-wise supports of $\beta^*$:
$$T = \{i: \beta_i \neq 0\}\subseteq \{1,\ldots, p\}, \quad G = \{j: \beta_{(j)} \neq 0\}\subseteq \{1,\ldots, d\}.$$
Roughly speaking, $u_{et}$ is the sub-gradient of objective function \eqref{eq:ell_1_minimization} evaluated at $\beta = \beta^*$. If $u_{et}$ is in the row space of $X$, the sub-gradient will be perpendicular to the feasible set of \eqref{eq:ell_1_minimization}, which implies that $\beta^*$ is the unique minimizer of $\ell_1+\ell_{1,2}$ minimization \eqref{eq:ell_1_minimization}. 

For more general vector $\beta^*$ that does not necessarily have a sparse support $T$ or $G$, we consider the following $(s, s_g)$-sparse approximation:
\begin{equation}\label{eq:T}
	\begin{split}
	\beta^{ap} = & \argmin_{S}  \frac{1}{\sqrt{s}}\|\beta^*_{S^c}\|_1 + \frac{1}{\sqrt{s_g}}\|\beta^*_{S^c}\|_{1, 2}\\
	& \text{subject to} \quad \|\beta^*_S\|_0 \leq s. \quad \|\beta^*_S\|_{0, 2} \leq s_g, \quad \max_{i \in S^c}\|\Sigma_{i, S}\Sigma_{S, S}^{-1}\|_2 \leq c/\sqrt{s}.
	\end{split}
\end{equation}
Let $T = \{i: \beta^{ap}_i \neq 0\}$ and $G = \{j: (\beta^{ap})_{(j)}\neq 0\}$ be the element-wise and group-wise supports of $\beta^{ap}$. Define
\begin{equation}\label{eq:exact-dual-certificate_non-sparse}
	\widetilde{u}_0 = \widetilde{v}_0 + \widetilde{w}_0\in\mathbb{R}^p, \quad  \left\{\begin{array}{ll}
	(\widetilde{v}_0)_{(j)} =\sqrt{s/s_g}\beta^*_{T,(j)}/\|\beta^*_{T,(j)}\|_2, &  j\in G;\\
	\|(\widetilde{v}_0)_{(j)}\|_2< \sqrt{s/s_g}, & j \in G^c; \end{array}\right. \quad 
	\left\{\begin{array}{ll}
	(\widetilde{w}_0)_T = \sgn(\beta^*_T)\\
	\|(\widetilde{w}_0)_{T^c}\|_\infty < 1.
	\end{array}\right.
\end{equation}
Here $\beta^*_{T,(j)}\in \mathbb{R}^{b_j}$ is the subvector $\beta^*$ restricted on the $j$-th group with all entries in $T^c$ set to zero. Similarly to the exactly sparse case, if $\widetilde{u}_{0}$ is in the row space of $X$ and the true $\beta^*$ is approximately $(s,s_g)$-sparse, the minimizer of \eqref{eq:ell_1_minimization} will be close to $\beta^*$.

However, it is often difficult to find an exact dual certificate that lies in the row space of $X$ and satisfies stringent conditions in \eqref{eq:exact-dual-certificate} or \eqref{eq:exact-dual-certificate_non-sparse}. We instead propose to analyze via the \emph{approximate dual certificate} defined as \eqref{ineq:dual-certificate-u} in the following lemma.
\begin{Lemma}[Approximate dual certificate for sparse group Lasso]\label{lm:approximate-dual-certificate}
	Suppose $T, G$ are element-wise and group-wise support defined in \eqref{eq:T}. $\widetilde{u}_0$ is defined in \eqref{eq:exact-dual-certificate_non-sparse}. Assume $X$ satisfies $\sigma_{\min}\left(X_T^\top X_T/n\right) \geq c_{\min}/2$. 
	If there exists $u\in \mathbb{R}^p$ in the row span of $X$ satisfying 
	\begin{equation}\label{ineq:dual-certificate-u}
	\begin{split}
	& \|u_T - (\widetilde{u}_0)_T\|_2 \cdot \max_{i\in T^c}\left\|X_T^\top X_i/n\right\|_2\leq c_{\min}/8, \\
	& \|H_{1/2}(u_{(G^c)})\|_{\infty, 2}\leq \sqrt{s_0}/2, \quad \|u_{(G)\backslash T}\|_\infty \leq 1/2,
	\end{split}
	\end{equation}
	Then the conclusion of Theorem \ref{th:noiseless} \eqref{ineq:noiseless_approx} holds with probability at least $1 - 2e^{-cn}$. Here, $H_{1/2}(\cdot)$ is the soft-thresholding operator defined at the beginning of Section \ref{sec:theory}. 
	
	If we additionally assume $\beta^*$ is $(s, s_g)$-sparse, then $\beta^*$ is the unique solution to the sparse group $\ell_1 + \ell_{1, 2}$ minimization \eqref{eq:ell_1_minimization} with probability at least $1 - 2e^{-cn}$.
\end{Lemma}
Lemma \ref{lm:approximate-dual-certificate} shows that the conclusion of Theorem \ref{th:noiseless} holds if there exists an approximate dual certificate $u$ satisfying the condition \eqref{ineq:dual-certificate-u}. The following lemma shows that, under the assumptions in Theorem \ref{th:noiseless}, one can find such an approximate dual certificate with high probability.
\begin{Lemma}\label{lm:dual_certificate_sufficiency}
	Suppose $X$ has group structure \eqref{eq:X-group} and satisfies Assumption \ref{as:design}. Recall $\sigma_{\min}(X_T^\top X_T/n)$ is the least eigenvalue of $X_T^\top X_T/n$. Then $\sigma_{\min}\left(X_T^\top X_T/n\right) \geq 1/2$ and \eqref{ineq:dual-certificate-u} holds with probability at least $1 - Ce^{-cn/s}$, where $T$ is defined in \eqref{eq:T}.
\end{Lemma}

Another key technical tool to the proof of Theorem \ref{th:noiseless} is the following Lemma, which shows that $X$ satisfies the restricted isometry property for all simultaneously element-wise and group-wise sparse vectors with high probability when there are enough samples.
\begin{Lemma}\label{lm:weakRIP}
	If $n \geq C(s_g\log(d/s_g) + s\log(es_gb))$,
\begin{equation}\label{ineq:weak_RIP}
\frac{c_{\min}}{2}\|\gamma\|_2^2 \leq \frac{1}{n}\|X\gamma\|_2^2 \leq (C_{\max} + \frac{c_{\min}}{2})\|\gamma\|_2^2, \quad \forall \gamma \in \{\gamma\in\bbR^p:  \|\gamma\|_0 \leq 2s, \|\gamma\|_{0,2} \leq 2s_g\}
\end{equation}
with probability at least $1 - 2e^{-cn}$. 
\end{Lemma}

	Next we briefly discuss the proof of Theorem \ref{th:sample_complexity}. Consider the quotient space $\bbR^{db}/\text{ker}(X) = \{[\gamma]:= x + \text{ker}(X), \gamma \in \bbR^{db}\}$ and define an associated norm as $\|[\gamma]\| = \inf_{v \in \text{ker}(X)}\{\lambda\|\gamma - v\|_1 + \lambda_g\|\gamma - v\|_{1,2}\}$. We show that there exist $N$ different $(s, s_g)$-sparse vectors $\beta^{(1)}, \dots, \beta^{(N)}$ such that $\log(N) \asymp s\log(es_gb/s) + s_g\log(d/s_g)$ and $\|[\beta^{(i)}]\| = 1, \|[\beta^{(i)}] - [\beta^{(j)}]\| \geq 2/9$ for all $1 \leq i \neq j \leq N$. By a property of the packing number and the fact that dim$(\bbR^{db}/\text{ker}(X)) \leq n$, we must have $N \leq 10^n$. Thus $n \gtrsim \log(N) \asymp s\log(es_gb/s) + s_g\log(d/s_g)$.

    We prove Theorem \ref{th:lower_noiseless} by contradiction. Assume that 
    \begin{equation}\label{eq48}
    	n < c_0\left(s\log(es_gb/s) + s_g\log(d/s_g)\right)
    \end{equation}
    for a sufficiently small constant $c_0$.
    Let $\|\cdot\| = \|\cdot\|_1 + \sqrt{s/s_g}\|\cdot\|_{1,2}$ and $B = \{x \in \bbR^{db}: \|x\| \leq 1\}$ be the unit ball associated with $\|\cdot\|$. Define 
    \begin{equation*}
    d^n(B, \bbR^p) = \inf_{\substack{L^n \text{ is a subspace of } \bbR^p \\\text{ with dim}(\bbR^p/L^n) \leq n} }\left\{\sup_{\beta \in B \cap L^n}\|\beta\|_2\right\},
    \end{equation*}
    We have $d^n(B, \bbR^p) \leq \frac{C}{\sqrt{s}}$ by the assumption of this theorem.
	We can also show that there exists a uniform constant $c > 0$ such that
    \begin{equation*}
    d^n(B, \bbR^p) \geq c\min\left\{\frac{1}{\sqrt{s_0}}, \left[\left(\frac{s_g}{s}\log\left(\frac{c\frac{s}{s_g}d\log(es_gb/s)}{n}\right) + \log(es_gb/s)\right)/n\right]^{1/2}\right\}.
    \end{equation*}
    The previous two inequalities and \eqref{eq48} together imply that 
    \begin{equation*}
    	n \geq c\left(s_g\log\left(\frac{c\frac{s}{s_g}d\log(es_gb/s)}{n}\right) + s\log(es_gb/s)\right) \geq c_0\left(s\log(es_gb/s) + s_g\log(d/s_g)\right) > n.
    \end{equation*}
    This contradiction shows that $n \geq c_0\left(s\log(es_gb/s) + s_g\log(d/s_g)\right)$.
    
%%%%%%%%%%%%%%%%%%%%%%%%%%%%%%
\section{Sparse Group Lasso in Noisy Case}\label{sec: noisy}
%%%%%%%%%%%%%%%%%%%%%%%%%%%%%%

We now turn to the noisy case.

\subsection{Optimal Rate of Estimation Error of Sparse Group Lasso}\label{subsec:rate}

When observations are noisy, we have the following theoretical guarantee for the sparse group Lasso. 
\begin{Theorem}[Upper bound of estimation error]\label{th:noisy}
	Suppose $y = X\beta^* + \varepsilon$, $X$ satisfies Assumption \ref{as:design}, $n \geq C\left(s_g \log(d/s_g) + s\log(es_gb)\right)$ for some uniform constant $C >  0$, $\varepsilon\overset{iid}{\sim} N(0, \sigma^2)$, and $b = \max_{1 \leq i \leq d}b_i$. Then the sparse group Lasso estimator \eqref{eq:sparse-group-lasso} with
	$$\lambda = C\sigma\sqrt{(s\log(es_gb)+s_g\log(ed/s_g))n/s} \quad\text{and}\quad \lambda_g = \sqrt{s/s_g}\lambda$$ 
	satisfies
	\begin{equation*}
	\|\hat{\beta} - \beta^*\|_2 \lesssim \min_{S: \substack{\|\beta^*_S\|_0 \leq s, \|\beta^*_S\|_{0, 2} \leq s_g, \\\max_{i \in S^c}\|\Sigma_{i, S}\Sigma_{S, S}^{-1}\|_2 \leq c/\sqrt{s}}}\left\{\sqrt{\frac{\sigma^2(s_g\log(d/s_g) + s\log(es_gb))}{n}} + \frac{\|\beta^*_{S^c}\|_1}{\sqrt{s}} + \frac{\|\beta^*_{S^c}\|_{1, 2}}{\sqrt{s_g}}\right\}
	\end{equation*}
	with probability at least $1 - C\exp\left(-C\frac{s\log(es_gb) + s_g\log(d/s_g)}{s}\right)$. 
	
	Especially, if $\beta^*$ is exactly $(s, s_g)$-sparse and $\max_{i \in T^c}\|\Sigma_{i, T}\Sigma_{T, T}^{-1}\|_2 \leq c/\sqrt{s}$
	holds, then
	\begin{equation}
    \begin{split}
    \|\hat{\beta} - \beta^*\|_2^2 \lesssim \frac{\sigma^2(s_g\log(d/s_g) + s\log (es_gb))}{n}
    \end{split}
	\end{equation}
	with probability at least $1 - C\exp\left(-C\frac{s\log(es_gb) + s_g\log(d/s_g)}{s}\right)$.
\end{Theorem}

In addition, we focus on the following class of simultaneously element-wise and group-wise sparse vectors,
$$\mathcal{F}_{s, s_g} = \{\beta: \|\beta\|_0 \leq s, \|\beta\|_{0,2} \leq s_g\}.$$
The following minimax lower bound of estimation error holds.
\begin{Theorem}[Lower bound of estimation error]\label{th:lower-noisy}
	Suppose $X$ satisfies Assumption \ref{as:design}, $b_1 = \cdots = b_d = b$, and $d, b \geq 3$. Then we have
	$$\inf_{\hat\beta}\sup_{\beta \in \mathcal{F}_{s, s_g}}\mathbb{E} \|\hat{\beta} - \beta\|_2^2 \gtrsim \frac{\sigma^2(s_g\log(ed/s_g) + s\log(es_gb/s))}{n}.$$
\end{Theorem}
\begin{Remark}\rm 
	Theorems \ref{th:noisy} and \ref{th:lower-noisy} together show that the sparse group Lasso yields the minimax optimal rate of convergence as long as the following condition holds: $\log(es_gb) \asymp \log(es_gb) - \log(s)$ or $\log(d) \gtrsim s\log(s)/s_g$.
\end{Remark}
\begin{Remark}\rm 
	We briefly discuss the main proof ideas of Theorem \ref{th:lower-noisy} here. First, we randomly generate a series of subsets $\Omega^{(i)} \subseteq \{1, \dots, p\}$ as feasible supports of $(s, s_g)$-sparse vectors. 	Then, we prove by a probabilistic argument that there exist $N \gtrsim \left(s_g\log(d/s_g) + s\log(es_gb/s)\right)$ subsets $\{\Omega^{(i)}\}_{i=1}^N$ such that $|\Omega^{(i)} \cap \Omega^{(j)}| < 8s_g\lfloor s/s_g\rfloor/9$ for any $i < j$. Next, we construct a series of candidate $(s, s_g)$-sparse vectors $\beta^{(i)}$ such that $\beta^{(i)}_k = \tau 1_{\{k\in \Omega^{(i)}\}}$. Intuitively speaking, $\{\beta^{(i)}\}_{i=1}^N$ are non-distinguishable based only on observations $(y, X)$ by such a construction. Theorem \ref{th:lower-noisy}  then follows by choosing an appropriate $\tau$ and the generalized Fano's lemma.
\end{Remark}

%%%%%%%%%%%%%%%%%%%%%%%%%%%%%%%
\subsection{Statistical Inference via Debiased Sparse Group Lasso}\label{subsec:inference}
%%%%%%%%%%%%%%%%%%%%%%%%%%%%%%%

We further consider the statistical inference for $\beta^*$ under the double sparse linear regression model. First, let $\hat{\beta}$ be the sparse group Lasso estimator given by \eqref{eq:sparse-group-lasso}. Inspired by the recent advances in inference for high-dimensional linear regression \cite{zhang2014confidence,van2014asymptotically,javanmard2014confidence,cai2017confidence}, we propose the following \emph{debiased sparse group Lasso estimator},
\begin{equation}\label{eq:debiased}
\hat{\beta}^u = \hat{\beta} + \frac{1}{n}\hat{M}X^\top\left(Y - X\hat{\beta}\right).
\end{equation}
Here, $\hat{\Sigma} = \frac{1}{n}\sum_{k=1}^n X_kX_k^\top$  is the sample covariance matrix and $\hat{M} = [\hat{m}_1 ~ \cdots ~ \hat{m}_p]^\top$ is an approximation of the inverse covariance matrix $\Sigma^{-1}$, where $\hat{m}_i$ is the solution to the following convex optimization,
\begin{equation}\label{opt:M}
\begin{split}
\text{minimize} \quad & m^\top \hat{\Sigma} m \\
\text{subject to}\quad &  \|H_{\alpha}(\hat{\Sigma}m - e_i)\|_{\infty, 2}\leq \gamma.
\end{split}
\end{equation}
Here, $H_\alpha$ is the soft-thresholding operator with thresholding level $\alpha$ defined at the beginning of Section \ref{sec:theory} and $e_i$ is the $i$-th vector in the canonical basis of $\mathbb{R}^p$. The following theorem establishes an asymptotic result for debiased sparse group Lasso.
\begin{Theorem}[Asymptotic distribution of debiased sparse group Lasso]\label{th:unbiased}
	Suppose $\beta^* \in \mathbb{R}^p$ is $(s, s_g)$-sparse, $X\in \mathbb{R}^{n\times p}$ satisfies Assumption \ref{as:design}, and $\max_{i \in T^c}\|\Sigma_{i, T}\Sigma_{T, T}^{-1}\|_2 \leq c/\sqrt{s}$. Set $\lambda =  C\sigma\sqrt{\frac{(s\log(es_gb) + s_g\log(d/s_g))n}{s}}$ and $\lambda_g =  \sqrt{\frac{s}{s_g}}\lambda$ in \eqref{eq:sparse-group-lasso},  $\alpha = \frac{\lambda}{n\sigma}, \gamma = \sqrt{\frac{s}{s_g}}\frac{\lambda}{n\sigma}$ in \eqref{opt:M}. 
Then with probability at least $1 - C\exp\left(-C\frac{s\log(es_gb) + s_g\log(d/s_g)}{s}\right)$, the debiased sparse group Lasso estimator $\hat{\beta}^u$ can be decomposed as $\sqrt{n}\left(\hat{\beta}^u - \beta^*\right) = \Delta + w,$ where
	\begin{equation}\label{de-biased1}
	\left\|\Delta\right\|_\infty \leq \frac{C\left(s\log(es_gb) + s_g\log(ed/s_g)\right)}{\sqrt{n}}\sigma, \quad w|X \sim N\left(0,  \sigma^2\hat{M}\hat{\Sigma}\hat{M}^\top\right).
	\end{equation}
	In particular, if $\sqrt{n} \gg s\log(es_gb) + s_g\log(ed/s_g)$, for any $1\leq i\leq p$,
	\begin{equation}\label{de-biased2}
	\frac{\sqrt{n}\left(\hat{\beta}^u_i - \beta^*_i\right)}{\sqrt{\hat{m}_i^\top\hat{\Sigma}\hat{m}_i}} \to N\left(0, \sigma^2\right).
	\end{equation}
\end{Theorem}

\begin{Remark}\rm 
	\eqref{de-biased2} provides a method to construct confidence intervals for $\beta^*$. Specifically if $\hat\sigma$ is a consistent estimator of $\sigma$, such as the scaled sparse group Lasso to be discussed in Section \ref{sec:disc}, 
	$$\left[\hat{\beta}_i^u - \Phi^{-1}(1 - \alpha/2)\hat{\sigma}\sqrt{\frac{\hat{m}_i^\top\hat{\Sigma}\hat{m}_i}{n}},\quad  \hat{\beta}_i^u + \Phi^{-1}(1 - \alpha/2)\hat{\sigma}\sqrt{\frac{\hat{m}_i^\top\hat{\Sigma}\hat{m}_i}{n}}\right]$$ 
	would be an asymptotic $(1 - \alpha)$-confidence interval for $\beta^*_i$. 
	We can see that the debiased sparse group Lasso estimator has the provably advantage on sample complexity ($n \gg (s\log(es_gb) + s_g\log(ed/s_g))^2$) over the ones via debiased Lasso ($n\gg s\log p$, see \cite{zhang2014confidence,javanmard2014confidence,cai2017confidence}) or debiased group Lasso ($n\gg (s_gb + s_g\log p )^2$, see \cite{mitra2016benefit}) for constructing asymptotic confidence intervals of $\beta^*$.
\end{Remark}

%%%%%%%%%%%%%%%%%%%%%%%%%
\section{Simulation Studies}\label{sec:simulation}
%%%%%%%%%%%%%%%%%%%%%%%%%

In this section, we investigate the numerical performance of the sparse group Lasso and $\ell_1+\ell_{1,2}$ minimization for double sparse regression. The results support our theoretical findings in Sections \ref{sec:theory} and \ref{sec: noisy}. We first discuss the practical choice for the tuning parameters used in the proposed algorithms.

%%%%%%%%%%%%%%%%%%
\subsection{Practical Selection of Tuning Parameters}\label{sec:tuning}
%%%%%%%%%%%%%%%%%%

By introducing $\tau$ as a surrogate for $(\lambda_g/\lambda)^2$, we can rewrite the $\ell_1+\ell_{1,2}$ minimization and the sparse group Lasso as
\begin{equation}\label{eq:SGL-minimization}
\begin{split}
\hat{\beta} = \argmin \|\beta\|_1 + \sqrt{\tau}\|\beta\|_{1,2} \quad \text{subject to} \quad y = X\beta,
\end{split}
\end{equation}
\begin{equation}\label{eq:SGL}
\hat\beta = \argmin_{\beta} \|y-X\beta\|_2^2 + \lambda\|\beta\|_1 + \lambda\sqrt{\tau} \|\beta\|_{1,2}.
\end{equation}
As suggested by Theorems \ref{th:noiseless} and \ref{th:noisy}, the theoretical choice of the tuning parameters $(\lambda, \tau)$ relies on $\sigma, s$, and $s_g$ in sparse group Lasso and $\ell_1+\ell_{1,2}$ minimization for double sparse regression. These values, however, are usually unknown in practice. In addition, those theoretical values of tuning parameters may not achieve the best finite-sample numerical performance. We thus introduce in this section a data-driven approach to tuning parameter selection using $K$-fold cross-validation.

We first discuss how to select $\tau$ in the $\ell_1+\ell_{1,2}$ minimization \eqref{eq:SGL-minimization}. Recall $n$ is the sample size, $p$ is the total number of covariates, $d$ is the number of groups, $b_1,\ldots, b_d$ are the number of covariates in each group, and $b = \max_j b_j$. Since the theoretical value $\tau = s/s_g$ and $s/s_g$ must satisfy $1\leq s/s_g\leq b$, for a given integer $L\ge 1$, we introduce a grid
	\begin{equation}\label{eq:S_0}
	S_0 = \{b^{(l-1)/(L-1)}: 1 \leq l \leq L\} 
	\end{equation}
	as a set of candidate values for $\tau$. Here, the grid size $L$ can be set to a typical value of 10, or a larger value if more computing power is available. 	We split the data $\{X_i, y_i\}_{i=1}^n$ into $K$ groups. For $1 \leq k \leq K$, let $J_k \subset \{1, \dots, n\}$ be the index set of the $k$th group and $J_k^c = \{1,\ldots, n\}\backslash J_k$. For each $\tau \in S_0$, we solve
	\begin{equation*}
	\begin{split}
	\hat{\beta}^{(k)}(\tau) = \argmin \|\beta\|_1 + \sqrt{\tau}\|\beta\|_{1,2} \quad \text{subject to} \quad y_{J_k^c} = X_{[J_k^c, :]}\beta
	\end{split}
	\end{equation*}
	and calculate the prediction error 
	\begin{equation*}
	\hat R(\tau) = \sum_{k=1}^{K}\sum_{j \in J_k}\left(y_j - X_{[j,:]}\hat{\beta}^{(k)}(\tau)\right)^2.
	\end{equation*}
	Let $\tau_*$  be the minimizer of the prediction error: $\tau_* = \argmin_{\tau \in S_0}\hat R(\tau).$ Then, the final estimator $\hat{\beta}$ is calculated using \eqref{eq:SGL-minimization} with $\tau_*$. 

	Then we consider the sparse group Lasso \eqref{eq:SGL}, which includes two tuning parameters $(\tau, \lambda)$. We still define $S_0$ in \eqref{eq:S_0} as a grid of candidate values of $\tau$. Following the idea in \cite[Section 3.3]{simon2013sparse}, for each $\tau \in S_0$, we begin with a large value of $\lambda_{\max}(\tau)$ so that $\hat{\beta}$, the outcome of sparse group Lasso \eqref{eq:SGL} with tuning parameters $(\tau, \lambda_{\max}(\tau))$, is zero (this can be achieved by the \texttt{SGL} package\footnote{\url{https://cran.r-project.org/web/packages/SGL/index.html}}). 
	Let  $\lambda_{\min}(\tau)$ be a small fraction of $\lambda_{\max}(\tau)$ (e.g., $\lambda_{\min} = 0.1\lambda_{\max}$ as suggested in \cite[Section 5]{simon2013sparse}). 
	Then we define $\Lambda(\tau) = \left\{\{\lambda_{\min}(\tau)\}^{(L-l)/(L-1)} \cdot \{\lambda_{\max}(\tau)\}^{(l-1)/(L-1)}: l =1,\ldots, L\right\}$. Next, we split the data $\{X_i, y_i\}_{i=1}^n$ into $K$ groups. For $1 \leq k \leq K$, let $J_k \subset \{1, \dots, n\}$ be the index set of the $k$th group and $J_k^c = \{1,\ldots, n\}\backslash J_k$. For each $\tau \in S_0$, $\lambda \in \Lambda(\tau)$, and $k\in \{1,\ldots, K\}$, we solve
\begin{equation*}
\begin{split}
\hat{\beta}^{(k)}(\tau, \lambda) = & \argmin_{\beta} \left\|y_{J_k^c} - X_{[J_k^c,:]}\beta\right\|_2^2 + \lambda\|\beta\|_1 + \lambda\sqrt{\tau}\|\beta\|_{1,2}
\end{split}
\end{equation*}
and calculate the prediction error 
\begin{equation*}
\hat R(\tau, \lambda) = \sum_{k=1}^{K}\sum_{j \in J_k}\left(y_j - X_{[j,:]}\hat{\beta}^{(k)}(\tau, \lambda)\right)^2.
\end{equation*}
Let $(\tau_*, \lambda_*)$ be the minimizer of the prediction error: $(\tau_*, \lambda_*) = \argmin_{\tau \in S_0, \lambda \in \Lambda(\tau)}\hat R(\tau, \lambda).$ The final estimator $\hat{\beta}$ is calculated using \eqref{eq:SGL} with $(\tau_*, \lambda_*)$.

In our simulation studies next, we will examine the performance of this cross-validation scheme with $K = L=10$, $\lambda_{\min}=0.1\lambda_{\max}$.

%%%%%%%%%%%%%%%%%%
\subsection{Numerical Results}\label{sec:numerical}
%%%%%%%%%%%%%%%%%%

\begin{figure}[!h]
	\centering
		\subfigure[$d= 60, b = 20, s_g = 1$]{
		\begin{minipage}[t]{0.48\linewidth}
			\centering
			\includegraphics[height=5cm]{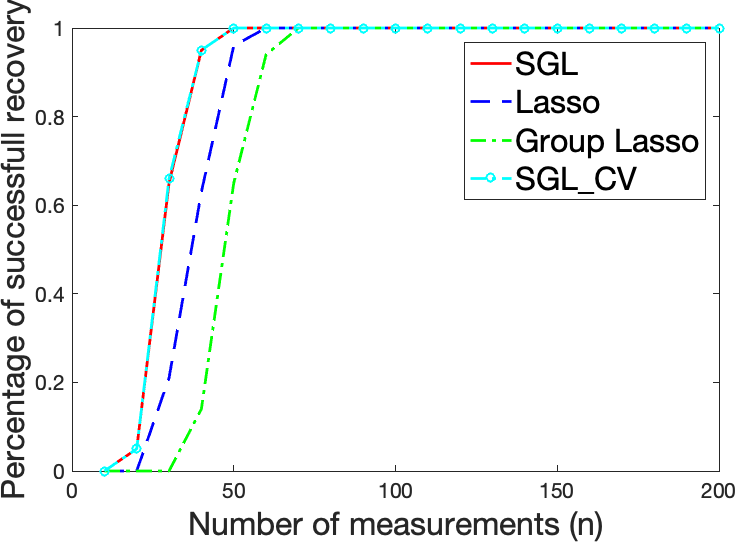}
		\end{minipage}%
	}
        \subfigure[$d = 100, b = 30, s_g = 2$]{
        \begin{minipage}[t]{0.48\linewidth}
        	\centering
        	\includegraphics[height=5cm]{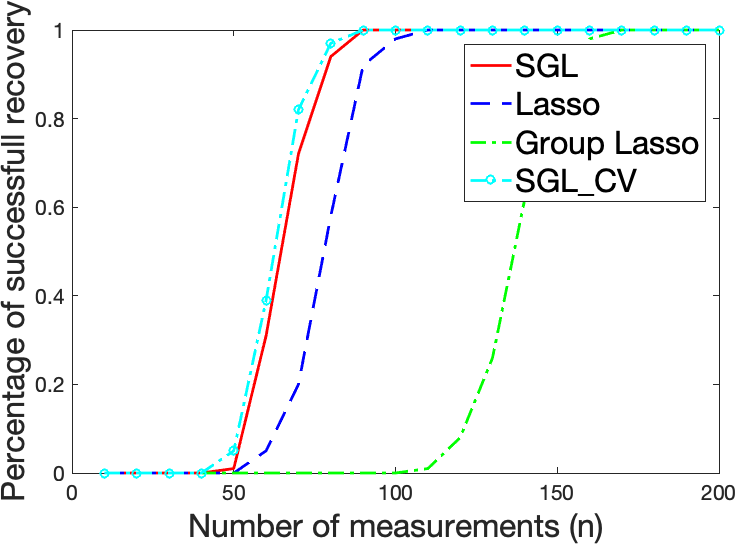}
        \end{minipage}%
    }
	\subfigure[$d= 20, b = 20, s_g = 1$]{
		\begin{minipage}[t]{0.48\linewidth}
			\centering
			\includegraphics[height=5cm]{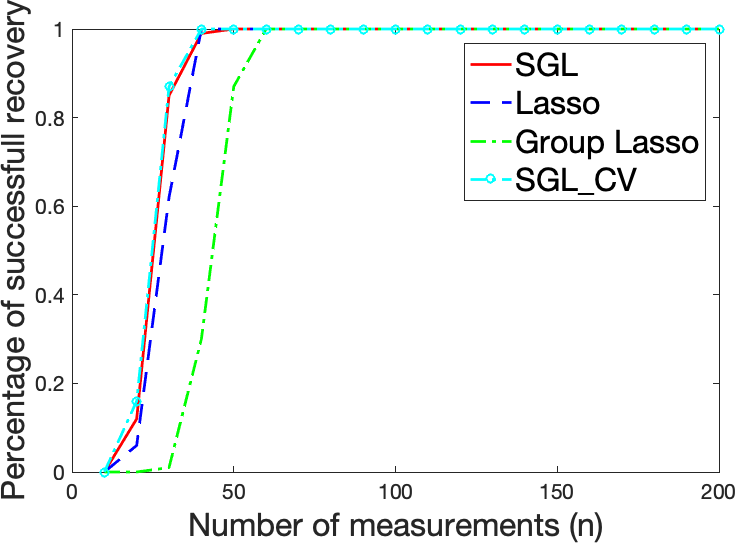}
		\end{minipage}%
	}
	\subfigure[$d = 40, b = 40, s_g = 1$]{
		\begin{minipage}[t]{0.48\linewidth}
			\centering
			\includegraphics[height=5cm]{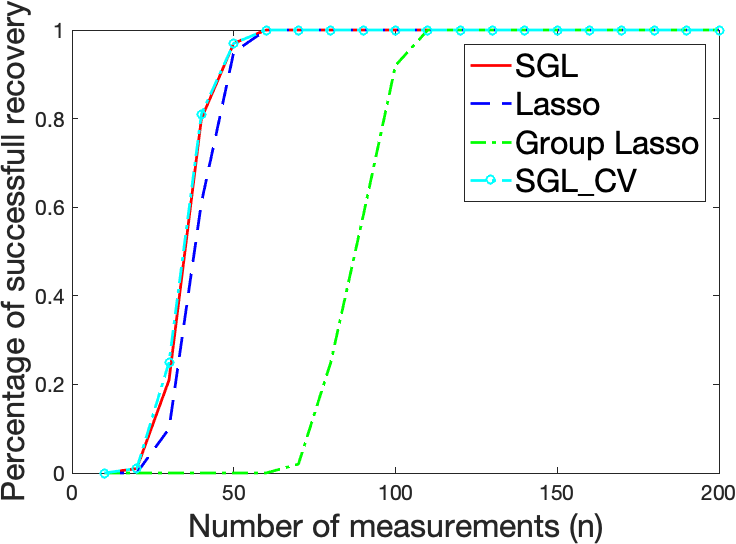}
		\end{minipage}%
	}
	\centering
	\caption{Exact recovery rate in the noiseless case}
	\label{fig1}
\end{figure}

We begin by considering the sample complexity for the exact recovery in the noiseless case. Suppose all group sizes are equal ($b_1 = \cdots = b_d = b$) and the number of observations $n$ varies from 5 to 200. 
We consider four simulation designs with (1) $d=60,b=20,s_g=1$; (2) $d=100,b=30,s_g=2$; (3) $d=b=20,s_g=1$; and (4) $d=b=40,s_g=1$. For each setting, we randomly draw $X \in \bbR^{n \times db}$ with i.i.d. standard normal entries, construct the fixed vector $\beta^*\in \mathbb{R}^{db}$ satisfying 
$$\beta^*_{(j)} =\left\{\begin{array}{ll}
(1,2,3,4,5,0,\dots,0)\in \mathbb{R}^{b} & j = 1, \dots, s_g; \\
0 & j = s_g + 1, \dots, d,
\end{array}\right. $$
and generate $y = X\beta^* = \sum_{j = 1}^{s_g}X_{(j)}\beta^*_{(j)}$.  We implement the $\ell_1+\ell_{1,2}$ minimization  \eqref{eq:ell_1_minimization} with $\lambda_g = \sqrt{s/s_g}\lambda$ (\texttt{SGL}), $\ell_1$ minimization \eqref{eq:l_1-mini} (\texttt{Lasso}), and $\ell_{1,2}$ minimization \eqref{eq:l_12-mini} (\texttt{Group Lasso}), and $\ell_1+\ell_{1,2}$ minimization \eqref{eq:ell_1_minimization} with the tuning parameter $\lambda_g/\lambda$ selected using cross validation discussed in Section \ref{sec:tuning} (\texttt{SGL\_CV}). An exact recovery of $\beta^*$ is considered to be successful if $\|\hat{\beta} - \beta^*\|_2 \leq 10^{-4}$. The successful recovery rate based on 100 replicates is shown in Figure \ref{fig1}. It can be seen that \texttt{SGL} and \texttt{SGL\_CV} have comparable performance and both methods have significantly better performance than \texttt{Lasso} and \texttt{Group Lasso}. This is in line with our theoretical results.

Then we consider the noisy case and focus on average estimation errors of different methods. We generate
\begin{equation*}
	y = X\beta^* + \varepsilon = \sum_{j = 1}^{s_g}X_{(j)}\beta^*_{(j)} + \varepsilon,
\end{equation*}
where $X, \beta^*$ are drawn in the same way as the previous setting and $\varepsilon \overset{iid}{\sim} N(0, 0.1^2)$. We consider four designs: i. $d=60,b=20,s_g=1$; ii. $d=100,b=30,s_g=2$; iii. $d=b=20,s_g=1$; and iv. $d=b=40,s_g=2$. For each case, the number of observations $n$ is chosen from an equally spaced sequence from 5 to 200 and the simulation is replicated for 500 times. We compare the average estimation error of (a) \texttt{SGL\char`_CV1}: sparse group Lasso with theoretical value $\lambda_g = \sqrt{s/s_g}\lambda$ and $\lambda$ selected via cross validation; (b) \texttt{SGL\char`_package}: sparse group Lasso via SGL package\footnote{\url{https://cran.r-project.org/web/packages/SGL/index.html}} in R with the option of automatic tuning parameter selection; (c) \texttt{Lasso}: regular Lasso with tuning parameter selected via cross validation; (d) \texttt{group Lasso}: group Lasso with tuning parameter selected via cross validation; (e) \texttt{SGL\char`_CV2}: sparse group Lasso with both $\lambda$ and $\lambda_g$ selected using the proposed cross validation scheme. We can see the proposed method \texttt{SGL\char`_CV2} achieves smaller estimation error than all other methods, including \texttt{SGL\char`_CV1}, the focus of our theory. These experimental results demonstrate our theory and the applicability of the proposed cross-validation scheme.

\begin{figure}[!h]
	\centering
	\subfigure[$d= 60, b = 20, s_g = 1$]{
		\begin{minipage}[t]{0.48\linewidth}
		\centering
			\includegraphics[height=6cm]{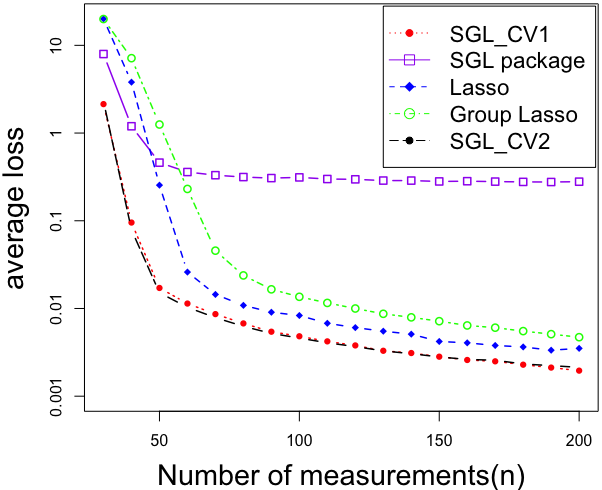}
		\end{minipage}
	}
	\subfigure[$d = 100, b = 30, s_g = 2$]{
		\begin{minipage}[t]{0.48\linewidth}
		\centering
			\includegraphics[height=6cm]{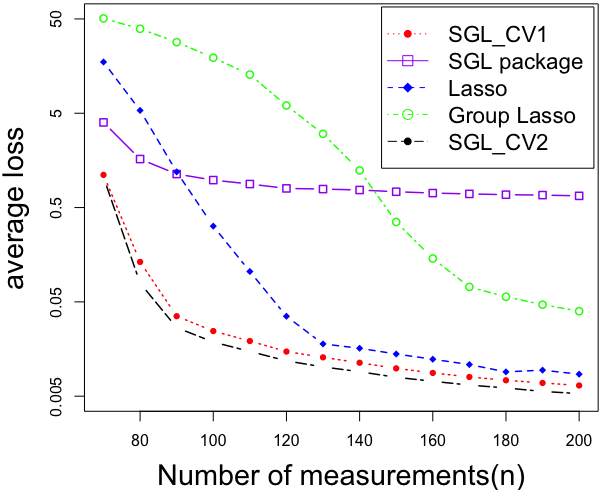}
		\end{minipage}
	}\\
	\subfigure[$d= 20, b = 20, s_g = 1$]{
		\begin{minipage}[t]{0.48\linewidth}
			\centering
			\includegraphics[height=6cm]{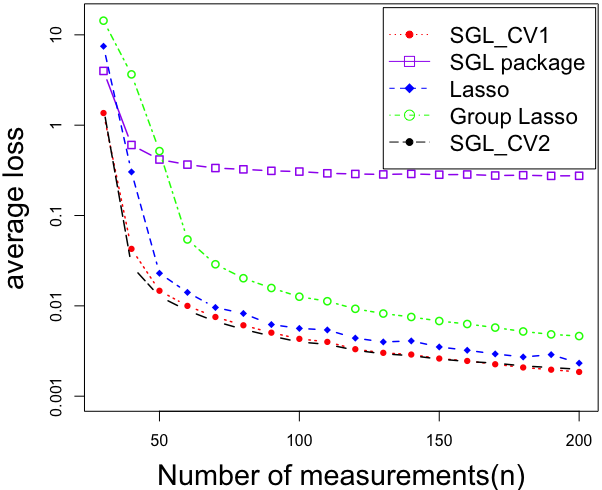}
		\end{minipage}%
	}
	\subfigure[$d = 40, b = 40, s_g = 2$]{
		\begin{minipage}[t]{0.48\linewidth}
			\centering
			\includegraphics[height=6cm]{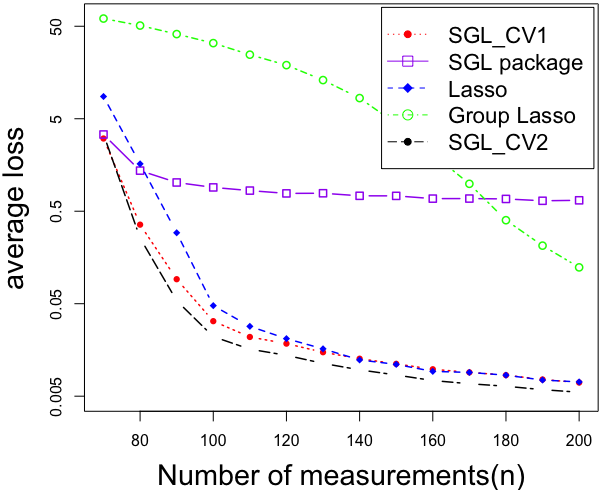}
		\end{minipage}%
	}
	\centering
	\caption{Average estimation error in the noisy case}
	\label{fig4}
\end{figure}

%%%%%%%%%%%%%%%%%%%%%%%%%%
\section{Discussions}\label{sec:disc}
%%%%%%%%%%%%%%%%%%%%%%%%%%

In this paper, we study the high-dimensional double sparse regression and investigate the theoretical properties of the sparse group Lasso and $\ell_1+\ell_{1,2}$ minimization. Particularly, we develop the matching upper and lower bounds on the sample complexity for $\ell_1+\ell_{1,2}$ minimization in the noiseless case. We also prove that the sparse group Lasso achieves minimax optimal rate of convergence in a range of settings in the noisy case. Our results give an affirmative answer to the open question for high-dimensional statistical inference for simultaneously structured model: by introducing both $\ell_1$ and $\ell_{1,2}$ penalties, one can achieve better performance on estimation and statistical inference for simultaneously element-wise and group-wise sparse vectors. 

In addition to $\beta^*$, the estimation and inference for noise level $\sigma$ is another importance task in high-dimensional double sparse regression. Motivated by the recent development of scaled Lasso \cite{sun2012scaled}, one may consider the following scaled sparse group Lasso estimator:
\begin{equation*}
\{\hat{\beta}^{s}, \hat{\sigma}\} = \argmin_{\beta \in \bbR^p, \sigma > 0}\left\{\frac{\|y - X\beta\|_2^2}{\sigma} + n\sigma + \widetilde{\lambda}\|\beta\|_1 + \widetilde{\lambda}_g\|\beta\|_2\right\},
\end{equation*}
where $\widetilde{\lambda}$ and $\widetilde{\lambda}_g$ are tuning parameters that do not rely on $\sigma$. The consistency of $\hat{\sigma}$ can be established based on similar ideas of scaled Lasso in the literature \cite{sun2012scaled,javanmard2014confidence} and the approximate dual certificate in this work. 

Moreover, our technical results can be useful in a variety of other problems with simultaneous sparsity structures. For example, \cite{tibshirani2005sparsity,rinaldo2009properties} considered the estimation of piece-wise constant sparse signals, i.e., both the signal vector and the difference between successive entries of the signal vector are sparse. \cite{jalali2013convex,jalali2019new} discussed the estimation of structured parameters where both the number of non-zero elements and the number of distinct values of the parameter vectors are small. \cite{sprechmann2010collaborative} considered the estimation of matrices with simultaneous sparsity structures within each block and among different blocks. It is interesting to further study the statistical limits, including the sample complexity and minimax optimal rate of  convergence for these problems. In particular, based on the specific sparsity structures of each problem, we can introduce corresponding multi-objective regularizers and the convex regularization methods. The corresponding approximate dual certificates can be proposed, constructed, and analyzed to provide strong theoretical guarantees.

%%%%%%%%%%%%%%%%%%%%%%%%%%%%%%%%%%%%%
\section{Proofs}\label{sec:proofs}
%%%%%%%%%%%%%%%%%%%%%%%%%%%%%%%%%%%%%

We collect the proofs of technical results in this section.

\subsection{Proof of Lemma \ref{lm:approximate-dual-certificate}}
Let $T$ satisfy \eqref{eq:T}. For convenience, we denote $s_0 = s/s_g$ and decompose $u$ as
\begin{equation}\label{eq37}
\begin{split}
u = v + w, \quad & v_i = \left\{
\begin{array}{ll}
u_i - \sqrt{s_0}\beta^*_i/\|\beta^*_{T, (j)}\|_2, & i\in T, i \in (j);\\
u_i, & i\in (G)\backslash T;\\
u_i - H_{1/2}(u_i), & i \in (G^{c}).
\end{array}
\right.\\
& w_{(j)} = \left\{\begin{array}{ll}
\sqrt{s_0}\beta^*_{T, (j)}/\|\beta^*_{T, (j)}\|_2, & j\in G;\\
H_{1/2}(u_{(j)}), & j\notin G.
\end{array}\right.\\
\end{split}
\end{equation}
Note that $|H_{1/2}(x) - x|\leq 1/2$ for any $x\in \mathbb{R}$. Based on the property of \eqref{ineq:dual-certificate-u}, $\|u_{(G)\backslash T}\|_\infty\leq 1/2$, then
\begin{equation}\label{ineq:dual-property-v}
\begin{split}
& \max_{i\in T^c}\left|v_i\right|\leq 1/2, \quad \|v_T - \sgn(\beta^*_T)\|_2 = \|u_T - (\widetilde{u}_0)_T\|_2 \leq \frac{c_{\min}}{8\max_{i\in T^c}\|X_T^\top X_i/n\|_2};
\end{split}
\end{equation}
\begin{equation}\label{ineq:dual-property-w}
\begin{split}
w_{(j)} = \sqrt{s_0}\beta^*_{T, (j)}/\|\beta^*_{T, (j)}\|_2, \text{ if } j\in G; \quad \|w_{(j)}\|_2 \leq \sqrt{s_0}/2, \text{ if } j\notin G.
\end{split}
\end{equation}
Suppose $\hat{\beta}$ is the minimizer to \eqref{eq:ell_1_minimization}, $h = \hat{\beta} - \beta^*$, then based on the sub-differential of $\|\beta\|_1$ and $\|\beta\|_{1,2}$, we have
\begin{equation}\label{ineq:lm-dual-1}
\begin{split}
\mathcal{P}(\hat{\beta}) = & \|\hat{\beta}\|_1 + \sqrt{s_0}\|\hat{\beta}\|_{1,2} = \|\beta^*+h\|_1 + \sqrt{s_0}\|\beta^*+h\|_{1,2} \\
\geq & \|\beta^*_T\|_1 + \sgn(\beta^*_T)^\top h_T + \|h_{T^c}\|_1 + \sqrt{s_0}\left(\|\beta^*_T\|_{1,2} + \sum_{j\in G}\frac{\beta_{T, (j)}^{*\top} h_{(j)}}{\|\beta^*_{T, (j)}\|_2} + \sum_{j\notin G}\|h_{(j)}\|_2\right)\\
& - \|\beta^*_{T^c}\|_1 - \sqrt{s_0}\|\beta^*_{T^c}\|_{1, 2}\\
\geq & \mathcal{P}(\beta^*) + \|h_{T^c}\|_1 + \sqrt{s_0}\|h_{(G^c)}\|_{1, 2} + \sgn(\beta^*_T)^\top h_T + \sum_{j\in G}\frac{\sqrt{s_0}\beta_{T,(j)}^{*\top} h_{(j)}}{\|\beta^*_{T,(j)}\|_2}\\ & - 2\|\beta^*_{T^c}\|_1 - 2\sqrt{s_0}\|\beta^*_{T^c}\|_{1, 2}.
\end{split}
\end{equation}
The last inequality comes from $\|\beta^*\|_1 = \|\beta^*_T\|_1 + \|\beta^*_{T^c}\|_1$ and $\|\beta^*\|_{1,2} \leq \|\beta^*_T\|_{1,2} + \|\beta^*_{T^c}\|_{1,2}$.

In particular, given $Xh = 0$ and that $u$ lies in the row span of $X$, we have $v^\top h + w^\top h = u^\top h = 0$. Therefore,
\begin{equation}\label{ineq:lm-dual-2}
\begin{split}
& \sgn(\beta^*_T)^\top h_T + \sum_{j\in G}\frac{\sqrt{s_0}\beta_{T,(j)}^{*\top} h_{(j)}}{\|\beta^*_{T,(j)}\|_2} = \sgn(\beta^*_T)^\top h_T - v^\top h + \sum_{j\in G}\frac{\sqrt{s_0}\beta_{T,(j)}^{*\top} h_{(j)}}{\|\beta^*_{T,(j)}\|_2} - w^\top h\\
= & - (v_T - \sgn(\beta^*_T))^\top h_T  - v_{T^c}^\top h_{T^c} - \sum_{j\in G} \left(w_{(j)} - \sqrt{s_0}\beta^*_{T,(j)} /\|\beta^*_{T,(j)}\|_2\right)^\top h_{(j)} - (w_{(G^c)})^\top h_{(G^c)}\\ 
\geq & - \|v_T-\sgn(\beta^*_T)\|_2\|h_T\|_2 - \|v_{T^c}\|_\infty \cdot\|h_{T^c}\|_1 \\
& - \max_{j\in G}\left\|w_{(j)} - \sqrt{s_0}\beta^*_{T,(j)} /\|\beta^*_{T,(j)}\|_2\right\|_2\cdot \|h_{(G)}\|_{1, 2} - \|w_{(G^c)}\|_{\infty, 2} \|h_{(G^c)}\|_{1, 2}\\
\overset{\eqref{ineq:dual-property-v} \eqref{ineq:dual-property-w}}{\geq} & -\|v_T - \sgn(\beta^*_T)\|_2\cdot \|h_T\|_2 - \|h_{T^c}\|_1/2 - \sqrt{s_0}\|h_{(G^c)}\|_{1,2}/2.
\end{split}
\end{equation}
Next note that $h = h_T + h_{T^c}$, we must have $X_{T}h_{T} = -X_{T^c}h_{T^c}$, then
\begin{equation}\label{ineq:lm-dual-3}
\begin{split}
\|h_T\|_2 = & \|(X_T^\top X_T/n)^{-1} X_T^\top X_T h_T/n\|_2 \leq \sigma_{\min}^{-1}(X_T^\top X_T/n) \|X_TX_{T^c}h_{T^c}/n\|_2\\
\leq & \frac{2}{c_{\min}}\cdot \max_{i\in T^c}\|X_T^\top  X_{i}/n\|_2 \cdot \|h_{T^c}\|_1.
\end{split}
\end{equation}
Combining \eqref{ineq:dual-property-v}, \eqref{ineq:lm-dual-2}, and \eqref{ineq:lm-dual-3}, one obtains 
$$\sgn(\beta^*_T)^\top h_T + \sum_{j\in G}\frac{\sqrt{s_0}\beta_{T,(j)}^{*\top} h_{(j)}}{\|\beta^*_{T,(j)}\|_2} \geq - 3/4\cdot \|h_{T^c}\|_1  -\sqrt{s_0}\|h_{(G^c)}\|_{1,2}/2.$$ 
Plug this inequality to \eqref{ineq:lm-dual-1}, we finally have
$$\mathcal{P}(\hat{\beta}) \geq \mathcal{P}(\beta^*) + \|h_{T^c}\|_1 /4 + \sqrt{s_0}\|h_{(G^c)}\|_{1,2}/2 -  2\|\beta^*_{T^c}\|_1 - 2\sqrt{s_0}\|\beta^*_{T^c}\|_{1, 2}. $$
Since $\hat{\beta}$ is the minimizer to \eqref{eq:ell_1_minimization}, we must have $\mathcal{P}(\hat{\beta}) \leq \mathcal{P}(\beta^*)$, then 
\begin{equation}\label{ineq:noiseless_h}
	\|h_{T^c}\|_1 /4 + \sqrt{s_0}\|h_{(G^c)}\|_{1,2}/2 \leq 2\|\beta^*_{T^c}\|_1 + 2\sqrt{s_0}\|\beta^*_{T^c}\|_{1, 2}.
\end{equation}
If $\beta^*$ is $(s, s_g)$-sparse, immediately we have $h_{T^c} =0$. Then $0 = X_T^\top X h = (X_T^\top X_T) h_T$. By $\sigma_{\min}(X_T^\top X_T/n) \geq c_{\min}/2$, we know $X_T^\top X_T/n$ is non-singular, then $h_T = 0$. 

Now, we consider the general case. Without loss of generality, suppose $G = \{1, \dots, g\}$, where $g \leq s_g$. Denote $T_1$ as the indices of the $s$ largest entries of $h_{(G) \backslash T}$, $T_2$ as the indices of the $s$ largest entries of $h_{(G) \backslash [T \cup T_1]}$, and so on. For $s_g + 1 \leq i \leq d$, denote $S_{i, 1}$ as the indices of the $\lfloor s/s_g\rfloor$ largest entries of $h_{(i)}$, $S_{i, 2}$ as the indices of the $\lfloor s/s_g\rfloor$ largest entries of $h_{(i)\backslash S_{i, 1}}$, and so on. Let $\widetilde{S}_1, \dots, \widetilde{S}_{\sum_{i = g + 1}^{d}\lceil b_i/\lfloor s/s_g\rfloor\rceil}$ be an arrangement of $S_{i, j} (1 \leq j \leq \lceil b_i/\lfloor s/s_g\rfloor\rceil, g + 1 \leq i \leq d)$ such that $\|h_{\widetilde{S}_1}\|_2^2 \geq \dots \geq \|h_{\widetilde{S}_{\sum_{i = g + 1}^{d}\lceil b_i/\lfloor s/s_g\rfloor\rceil}}\|_2^2$. Let $R_1 = \cup_{i = 1}^{s_g}\widetilde{S}_i$, $R_2 = \cup_{i = s_g + 1}^{2s_g}\widetilde{S}_i$, and so on. Then $(T_1, T_2, \dots, R_1, R_2, \dots)$ is a partition of $T^c$, and $|T_i|, |R_j| \leq s, |g(T_i)|, |g(R_j)| \leq s_g$, where $g(S) = \{i_1, \dots, i_k\}$ if $S \subseteq \cup_{j = 1}^k(i_j)$ and $S \cap (i_j)$ are not empty for all $1 \leq j \leq k$. Let $\widetilde{T} = T \cup T_1 \cup R_1$. 
If \eqref{ineq:weak_RIP} holds, then 
\begin{equation}\label{eq27}
\frac{c_{\min}}{2}\|h_{\widetilde{T}}\|_2^2 \leq \frac{1}{n}\|X_{\widetilde{T}}h_{\widetilde{T}}\|_2^2 = \frac{1}{n}\langle X_{\widetilde{T}}h_{\widetilde{T}}, Xh\rangle - \frac{1}{n}\langle X_{\widetilde{T}}h_{\widetilde{T}}, X_{\widetilde{T}^c}h_{\widetilde{T}^c}\rangle.
\end{equation}
Since $Xh = 0$, we have
\begin{equation}\label{eq26}
	\langle X_{\widetilde{T}}h_{\widetilde{T}}, Xh\rangle = 0.
\end{equation}
Now, we consider $\left|\langle X_{\widetilde{T}}h_{\widetilde{T}}, X_{\widetilde{T}^c}h_{\widetilde{T}^c}\rangle\right|$. By triangle inequality, 
\begin{equation*}
	\left|\langle X_{\widetilde{T}}h_{\widetilde{T}}, X_{\widetilde{T}^c}h_{\widetilde{T}^c}\rangle\right| \leq \left|\langle X_{T}h_{T}, X_{\widetilde{T}^c}h_{\widetilde{T}^c}\rangle\right| + \left|\langle X_{T_1}h_{T_1}, X_{\widetilde{T}^c}h_{\widetilde{T}^c}\rangle\right| + \left|\langle X_{R_1}h_{R_1}, X_{\widetilde{T}^c}h_{\widetilde{T}^c}\rangle\right|.
\end{equation*}
The triangle inequality shows that
\begin{equation*}
\left|\langle X_{T}h_{T}, X_{\widetilde{T}^c}h_{\widetilde{T}^c}\rangle\right| \leq \sum_{i \geq 2}\left|\langle X_{T}h_{T}, X_{T_i}h_{T_i}\rangle\right| + \sum_{j \geq 2}\left|\langle X_{T}h_{T}, X_{R_j}h_{R_j}\rangle\right|.
\end{equation*}
Combine the parallelogram identity and \eqref{ineq:weak_RIP} together, we have
\begin{equation*}
\left|\langle X_{T}h_{T}, X_{T_i}h_{T_i}\rangle\right| \leq C_{\max}n\|h_T\|_2\|h_{T_i}\|_2, \quad \left|\langle X_{T}h_{T}, X_{R_j}h_{R_j}\rangle\right| \leq C_{\max}n\|h_T\|_2\|h_{R_j}\|_2.
\end{equation*}
Thus,
\begin{equation}\label{eq18}
\left|\langle X_{T}h_{T}, X_{\widetilde{T}^c}h_{\widetilde{T}^c}\rangle\right| \leq C_{\max}n\|h_T\|_2(\sum_{i \geq 2}\|h_{T_i}\|_2 + \sum_{j \geq 2}\|h_{R_j}\|_2). 
\end{equation}
By (3.10) in \cite{candes2007dantzig}, we have
\begin{equation}\label{eq19}
\sum_{i \geq 2}\|h_{T_i}\|_2 \leq s^{-1/2}\|h_{(G)\backslash T}\|_1,
\end{equation}
and 
\begin{equation*}
    \begin{split}
    \sum_{j \geq 2}\|h_{R_j}\|_2 =& \sum_{j \geq 2}\left(\sum_{i=(j-1)s_g + 1}^{js_g}\|h_{\widetilde{S}_i}\|_2^2\right)^{1/2} \leq \sum_{j \geq 2}\sqrt{s_g}\|h_{\widetilde{S}_{(j-1)s_g}}\|_2  \leq \sum_{j \geq 2}\sqrt{s_g}\sum_{i=(j-2)s_g+1}^{(j-1)s_g}\|h_{\widetilde{S}_{i}}\|_2/s_g  \\=& s_g^{-1/2}\sum_{k}\|h_{\widetilde{S}_k}\|_2 = s_g^{-1/2}\sum_{i=g + 1}^d\sum_{j}\|h_{S_{i, j}}\|_2.
    \end{split}
\end{equation*}
For all $g + 1 \leq i \leq d$, apply (3.10) in \cite{candes2007dantzig} again, 
\begin{equation*}
\sum_{j \geq 2}\|h_{S_{i, j}}\|_2 \leq (\lfloor s/s_g\rfloor)^{-1/2}\|h_{(i)}\|_1 \leq \sqrt{2}(s/s_g)^{-1/2}\|h_{(i)}\|_1.
\end{equation*}
Moreover, by the definition of $S_{i,1}$,
\begin{equation*}
\sum_{i=g + 1}^d\|h_{S_{i, 1}}\|_2 \leq \sum_{i=g + 1}^d\|h_{(i)}\|_2 = \|h_{(G^c)}\|_{1, 2}.
\end{equation*}
Therefore,
\begin{equation}\label{eq20}
\begin{split}
\sum_{j \geq 2}\|h_{R_j}\|_2 \leq& s_g^{-1/2}\left(\sum_{i=g + 1}^d\sqrt{2}(s/s_g)^{-1/2}\|h_{(i)}\|_1\right) + s_g^{-1/2}\|h_{(G^c)}\|_{1, 2}\\ =& \sqrt{2}s^{-1/2}\|h_{(G^c)}\|_1 + s_g^{-1/2}\|h_{(G^c)}\|_{1, 2}.
\end{split}
\end{equation}
Combine \eqref{eq18}, \eqref{eq19} and \eqref{eq20} together, if \eqref{ineq:weak_RIP} holds, we have
\begin{equation*}
\begin{split}
\left|\langle X_{T}h_{T}, X_{\widetilde{T}^c}h_{\widetilde{T}^c}\rangle\right| \leq& C_{\max}n\|h_T\|_2(s^{-1/2}\|h_{(G)\backslash T}\|_1 + \sqrt{2}s^{-1/2}\|h_{(G^c)}\|_1 + s_g^{-1/2}\|h_{(G^c)}\|_{1, 2})\\
\leq& C_{\max}n\|h_T\|_2(\sqrt{2}s^{-1/2}\|h_{T^c}\|_1 + s_g^{-1/2}\|h_{(G^c)}\|_{1, 2}).
\end{split}
\end{equation*}
Similarly, if \eqref{ineq:weak_RIP} holds, then
$\left|\langle X_{T_1}h_{T_1}, X_{\widetilde{T}^c}h_{\widetilde{T}^c}\rangle\right| \leq C_{\max}n\|h_{T_1}\|_2(\sqrt{2}s^{-1/2}\|h_{T^c}\|_1 + s_g^{-1/2}\|h_{(G^c)}\|_{1, 2})$ and $\left|\langle X_{R_1}h_{R_1}, X_{\widetilde{T}^c}h_{\widetilde{T}^c}\rangle\right| \leq C_{\max}n\|h_{R_1}\|_2(\sqrt{2}s^{-1/2}\|h_{T^c}\|_1 + s_g^{-1/2}\|h_{(G^c)}\|_{1, 2})$. 
Thus, with probability at least $1 - 2e^{-cn}$,
\begin{equation}\label{eq22}
\begin{split}
\left|\langle X_{\widetilde{T}}h_{\widetilde{T}}, X_{\widetilde{T}^c}h_{\widetilde{T}^c}\rangle\right| \leq& C_{\max}n\left(\|h_T\|_2 + \|h_{T_1}\|_2 + \|h_{R_1}\|_2\right)(\sqrt{2}s^{-1/2}\|h_{T^c}\|_1 + s_g^{-1/2}\|h_{(G^c)}\|_{1, 2})\\
\leq& \sqrt{3}C_{\max}n\|h_{\widetilde{T}}\|_2(\sqrt{2}s^{-1/2}\|h_{T^c}\|_1 + s_g^{-1/2}\|h_{(G^c)}\|_{1, 2}).
\end{split}
\end{equation}
The last inequality holds due to Cauchy-Schwarz inequality. Combine \eqref{eq27}, \eqref{eq26}, \eqref{eq22} and Lemma \ref{lm:weakRIP} together, we know that with probability at least $1 - 2e^{-cn}$,
\begin{equation*}
    \frac{c_{\min}}{2}\|h_{\widetilde{T}}\|_2^2 \leq \sqrt{3}C_{\max}\|h_{\widetilde{T}}\|_2(\sqrt{2}s^{-1/2}\|h_{T^c}\|_1 + s_g^{-1/2}\|h_{(G^c)}\|_{1, 2}),
\end{equation*}
i.e., with probability at least $1 - 2e^{-cn}$,
\begin{equation*}
	\|h_{\widetilde{T}}\|_2 \leq 2\sqrt{3}\frac{C_{\max}}{c_{\min}}(\sqrt{2}s^{-1/2}\|h_{T^c}\|_1 + s_g^{-1/2}\|h_{(G^c)}\|_{1, 2}).
\end{equation*}
Finally, by \eqref{ineq:noiseless_h}, \eqref{eq19}, \eqref{eq20} and the previous inequality, with probability at least $1 -  2e^{-cn}$, 
\begin{equation*}
\begin{split}
\|h\|_2 \leq& \|h_{\widetilde{T}}\|_2 + \sum_{i \geq 2}\|h_{T_i}\|_2 + \sum_{j \geq 2}\|h_{R_j}\|_2\\ \leq& 2\sqrt{3}\frac{C_{\max}}{c_{\min}}(\sqrt{2}s^{-1/2}\|h_{T^c}\|_1 + s_g^{-1/2}\|h_{(G^c)}\|_{1, 2}) + \sqrt{2}s^{-1/2}\|h_{T^c}\|_2 + s_g^{-1/2}\|h_{(G^c)}\|_{1, 2}\\ \leq& C\left(\frac{1}{\sqrt{s}}\|\beta^*_{T^c}\|_1 + \frac{1}{\sqrt{s_g}}\|\beta^*_{T^c}\|_{1, 2}\right).
\end{split}
\end{equation*}

In summary, we have finished the proof of this lemma. \quad $\square$

\subsection{Proof of Lemma \ref{lm:dual_certificate_sufficiency}}

Let $T$ satisfy \eqref{eq:T}. Given $\|\beta^*_T\|_{0,2} \leq s_g$, without loss of generally we assume that 
$$\beta^*_{T,(s_g+1)},\cdots, \beta^*_{T,(d)} =0.$$
We also denote $T_{(j)}$ as the support of $\beta^*_{T,(j)}$. First by Lemma \ref{lm:sub-Gaussian-concentration} Part 3 with 
$$v \in \mathbb{R}^{p}, v_k = \left\{\begin{array}{ll}
1, & k=i;\\
0, & k\neq i;
\end{array}\right.\quad U \in \mathbb{R}^{p \times |T|} = \mathbb{R}^{(\sum_{i = 1}^{d}b_i) \times |T|}, U_{[T, :]} = I; U_{[T^c, :]} = 0,$$ and notice that $x\log(eu/x) \geq \log(eu)$ for all $1 \leq x \leq u$,
we have
\begin{equation}\label{ineq:thm1-1}
\begin{split}
& \bbP\left(\max_{i\in T^c} \left\|X_T^\top X_i/n\right\|_2 \geq 1/2\right) \leq \sum_{i\in T^c} \bbP\left(\left\|X_T^\top X_i/n\right\|_2 \geq 1/2\right)\\
\leq& \sum_{i\in T^c} \bbP\left(\left\|X_T^\top X_i/n - \mathbb{E} X_T^\top X_i/n\right\|_2 + \|\mathbb{E}X_T^\top X_i/n\|_2 \geq 1/2\right)\\
\leq& \sum_{i\in T^c} \bbP\left(\left\|X_T^\top X_i/n - \mathbb{E} X_T^\top X_i/n\right\|_2 \geq 1/2 - \|\Sigma_{T, T}\|\|\Sigma_{i, T}\Sigma_{T, T}^{-1}\|_2\right)\\
\leq& \sum_{i\in T^c} \bbP\left(\left\|X_T^\top X_i/n - \mathbb{E} X_T^\top X_i/n\right\|_2 \geq 1/4\right)\\
\leq & db \cdot C\exp\left(Cs - n\right) \leq  C\exp\left(\log(d) + \log(b) + Cs - n\right)\\
 \leq& C\exp\left(s_g\log(ed/s_g) + s\log(es_gb/s) + Cs - n\right) \leq C\exp(-cn)
\end{split}
\end{equation}
provided that $n\geq C\left(s\log(es_gb/s) + s_g\log(d/s_g)\right)$ for some large constant $C>0$. Note that the fourth inequality comes from the facts that $\|\Sigma_{T, T}\| \leq \|\Sigma\| \leq C_{\max}$ and $\|\Sigma_{i, T}\Sigma_{T, T}^{-1}\|_2 \leq c/\sqrt{s} \leq 1/(4C_{\max})$. By Lemma \ref{lm:X-random-matrix-properties} Part 1, we also know
\begin{equation*}
\begin{split}
\bbP\left(\sigma_{\min}(X_T^\top X_T/n) \leq c_{\min}/2 \right) \leq & \bbP\left(\|X_T^\top X_T/n - \Sigma_{T, T}\| \geq c_{\min}/2\right)\\
\leq& \bbP\left(\|X_T^\top X_T\Sigma_{T, T}^{-1}/n - I_{|T|}\|\|\Sigma_{T, T}\| \geq c_{\min}/2\right)\\
\leq& \bbP\left(\|X_T^\top X_T\Sigma_{T, T}^{-1}/n - I_{|T|}\| \geq c_{\min}/(2C_{\max})\right)\\
\leq & C\exp\left(Cs - cn\right) \leq C\exp(-cn).
\end{split}
\end{equation*}
Next, we apply the well-regarded golfing scheme \cite{gross2011recovering,candes2011probabilistic} to find an approximate dual certificate $u$ that satisfies \eqref{ineq:dual-certificate-u}. Let
\begin{equation}
u_0\in\mathbb{R}^p, \quad (u_0)_{(j)} = \left\{\begin{array}{ll}
\sqrt{s/s_g}\beta^*_{T,(j)}/\|\beta^*_{T,(j)}\|_2 + \sgn(\beta^*_{T,(j)}), &  j\in G;\\
0, & j \in G^c.
\end{array}\right. 
\end{equation}
Immediately we have $(u_0)_T = (\widetilde{u}_0)_T$.
We divide $n$ rows of $X$ into non-overlapping batches, say $X_{[I_1, :]}, X_{[I_2, :]}, \ldots,$ with $|I_l| = n_l$. Here, $n_1, n_2, \ldots$ will be specified a little while later. Consider the following sequences
\begin{equation}
\begin{split}
& \alpha_0 = u_0, \\ 
& \gamma_l = X_{[I_l, :]}^\top X_{[I_l, T]}\Sigma_{T, T}^{-1}/n_l\cdot (\alpha_{l-1})_T, \quad \alpha_l = \alpha_{l-1}-\gamma_l,\quad l=1,2, \ldots, l_{\max}.
\end{split}
\end{equation}
Finally the approximate dual certificate is defined as
\begin{equation}\label{eq:approximate dual certificate}
u = \sum_{l=1}^{l_{\max}} \gamma_l = \sum_{l=1}^{l_{\max}} X_{[I_l, :]}^\top X_{[I_l, T]}\Sigma_{T, T}^{-1}/n_l \cdot (\alpha_{l-1})_T.
\end{equation}
From the inductive definition we can see 
$$(\alpha_l)_T = (I - X_{[I_l, T]}^\top X_{[I_l, T]}\Sigma_{T, T}^{-1}/n_l)(\alpha_{l-1})_T, \quad (\gamma_l)_{T^c} = X_{[I_l, T^c]}^\top X_{[I_l, T]}\Sigma_{T, T}^{-1}/n_l \cdot (\alpha_{l-1})_T, \quad l=1,2,\ldots. $$
Next, we apply the random matrix results (Lemmas \ref{lm:X-random-matrix-properties} and \ref{lm:sub-Gaussian-concentration}) and obtain the following tail probabilities.
\begin{itemize}
	\item if $n_l \geq Cst_l$ for large constant $C>0$ and $t_l\geq C$, by Part 1 of Lemma \ref{lm:X-random-matrix-properties},
\begin{equation}\label{ineq:thm-1-3}
\begin{split}
& \bbP\left(\|X_{[I_l, T]}^\top X_{[I_l, T]}\Sigma_{T, T}^{-1}/n_l - I_{|T|}\| \geq C\sqrt{st_l/n_l}\right) \\
\leq & C\exp\left(Cs - n_l\min\left\{\frac{st_l}{n_l}, \left(\frac{st_l}{n_l}\right)^{1/2}\right\}\right) \leq C\exp\left(-cst_l\right);
\end{split}
\end{equation}
\item Suppose $q_{l-1} = (\alpha_{l-1})_T \in \mathbb{R}^{|T|}$ is independent of $X_{[I_l, :]}$. If $n_l \geq \frac{C (s_0\log (es_gb/s) + \log d)}{\min\{s_0\delta_l^2, \sqrt{s_0}\delta_l\}}$ for $\delta_l \geq C\max_{i \in T^c}\|\Sigma_{i, T}\Sigma_{T, T}^{-1}\|_2 \geq C(\max_{i \in T^c}\|\Sigma_{i, T}\Sigma_{T, T}^{-1}\|_2)\|\Sigma_{T,T}^{-1}q_{l-1}\|_2/\|q_{l-1}\|_2$, by Lemma \ref{lm:X-random-matrix-properties} Part 2,
\begin{equation}\label{ineq:thm-1-4}
\begin{split}
& \bbP\left(\max_{j\in G^c}\left\|H_{\|q_{l-1}\|_2\delta_l}\left(X_{[I_l, (j)]}^\top X_{[I_l, T]}\Sigma_{T, T}^{-1}/n_l\cdot q_{l-1}\right)\right\|_2 \geq \sqrt{s_0}\|q_{l-1}\|_2\delta_l \right)\\ 
\leq & \sum_{j\in G^c}\bbP\left(\left\|H_{\|q_{l-1}\|_2\delta_l}\left(X_{[I_l, (j)]}^\top X_{[I_l, T]}(\Sigma_{T, T}^{-1}q_{l-1})/n_l\right)\right\|_2 \geq \sqrt{s_0}\|q_{l-1}\|_2\delta_l\right)\\
\leq& d\cdot \binom{b}{\lceil s_0\rceil}\exp\left(Cs_0 - cn_l\min\left\{\frac{s_0\|q_{l - 1}\|_2^2\delta_l^2}{\kappa^4\|\Sigma_{T, T}^{-1}q_{l-1}\|_2^2}, \frac{\sqrt{s_0}\|q_{l - 1}\|_2\delta_l}{\kappa^2\|\Sigma_{T, T}^{-1}q_{l-1}\|_2}\right\}\right)\\
&+ d\cdot \binom{b}{\lfloor s_0\rfloor}\exp\left(Cs_0 - cn_l\min\left\{\frac{s_0\|q_{l - 1}\|_2^2\delta_l^2}{\kappa^4\|\Sigma_{T, T}^{-1}q_{l-1}\|_2^2}, \frac{\sqrt{s_0}\|q_{l - 1}\|_2\delta_l}{\kappa^2\|\Sigma_{T, T}^{-1}q_{l-1}\|_2}\right\}\right)\\
\leq & 2d\cdot \left(\frac{eb}{\lfloor s_0\rfloor}\right)^{\lceil s_0\rceil}\exp\left(Cs_0 - cn_l \min\{s_0\delta_l^2, \sqrt{s_0}\delta_l\} \right)\\
\leq & C\exp\left(\log(d) + Cs_0\log(2es_gb/s) + Cs_0 - cn_l \min\{s_0\delta_l^2, \sqrt{s_0}\delta_l\} \right)\\
\leq & C\exp(-cn_l\min\{s_0\delta_l^2, \sqrt{s_0}\delta_l\});
\end{split}
\end{equation}
The third inequality comes from $\|\Sigma_{T, T}^{-1}\| \leq \frac{1}{c_{\min}}$.
\item Suppose $q_{l-1} = (\alpha_{l-1})_T \in \mathbb{R}^{|T|}$ is fixed. If $n_l\min\{\theta_l^2, \theta_l\} \geq C\log(es_gb), \theta_l \geq 2\max_{i \in T^c}\|\Sigma_{i, T}\Sigma_{T, T}^{-1}\|_2$, by Lemma \ref{lm:sub-Gaussian-concentration} Part 2,
\begin{equation}\label{ineq:thm-1-5}
\begin{split}
& \bbP\left(\|X_{[I_l, (G)\backslash T]}^\top X_{[I_l, T]}\Sigma_{T, T}^{-1}/n_l\cdot q_{l-1}\|_\infty \geq \theta_l\|q_{l-1}\|_2\right) \\
\leq & \sum_{i\in (G)\backslash T} \bbP\left(|X_{[I_l, i]}^\top X_{[I_l, T]}/n_l \cdot (\Sigma_{T, T}^{-1}q_{l-1})| \geq \theta_l\|q_{l-1}\|_2\right) \\
\leq& \sum_{i\in (G)\backslash T} \bbP\left(|X_{[I_l, i]}^\top X_{[I_l, T]}/n_l \cdot (\Sigma_{T, T}^{-1}q_{l-1}) - \Sigma_{i, T}\Sigma_{T, T}^{-1}q_{l - 1}| \geq \theta_l\|q_{l-1}\|_2 - |\Sigma_{i, T}\Sigma_{T, T}^{-1}q_{l - 1}|\right) \\
\leq& \sum_{i\in (G)\backslash T} \bbP\left(|X_{[I_l, i]}^\top X_{[I_l, T]}/n_l \cdot (\Sigma_{T, T}^{-1}q_{l-1}) - \Sigma_{i, T}\Sigma_{T, T}^{-1}q_{l - 1}| \geq \theta_l\|q_{l-1}\|_2 - \|\Sigma_{i, T}\Sigma_{T, T}^{-1}\|_2\|q_{l - 1}\|_2\right) \\
\leq& \sum_{i\in (G)\backslash T} \bbP\left(|X_{[I_l, i]}^\top X_{[I_l, T]}/n_l \cdot (\Sigma_{T, T}^{-1}q_{l-1}) - \Sigma_{i, T}\Sigma_{T, T}^{-1}q_{l - 1}| \geq \frac{1}{2}\theta_l\|q_{l-1}\|_2\right) \\
\leq& \sum_{i\in (G)\backslash T} \bbP\left(|X_{[I_l, i]}^\top X_{[I_l, T]}/n_l \cdot (\Sigma_{T, T}^{-1}q_{l-1}) - \Sigma_{i, T}\Sigma_{T, T}^{-1}q_{l - 1}| \geq \frac{c_{\min}}{2}\theta_l\|\Sigma_{T, T}^{-1}q_{l-1}\|_2\right) \\
\leq & s_gb\cdot C \exp\left(-cn_l\min\{\theta_l^2, \theta_l\}\right) = C\exp(\log(s_gb) - cn_l\min\{\theta_l^2, \theta_l\}) \\
\leq & C\exp\left(-cn_l\min\{\theta_l^2, \theta_l\}\right).
\end{split}
\end{equation}
\end{itemize}
Then we specify $\{n_l, t_l, \delta_l, \theta_l\}_{l\geq 1}$ as follows,
\begin{itemize}
	\item $n_1 = n_2 \geq C(s\log (es_gb) + s_g\log (d/s_g))$, $t_1 = t_2 = cn_1/(s\log (es)) \geq C$, $\delta_1 = \delta_2 = 1/(16\sqrt{s})$, $\theta_1 = \theta_2 = 1/(16\sqrt{s})$; 
	\item $n_3 = \cdots = n_{l_{\max}} \asymp \frac{n_1}{l_{\max}-2} \geq C(s\log (es_gb) + s_g\log (d/s_g))/\log(es)$, $t_3=\cdots = t_{l_{\max}} =cn_3/s \geq C$, $\delta_3 = \cdots = \delta_{l_{\max}} = \log(es)/(16\sqrt{s}) \geq \max\{(\log(es)/s)^{1/2}/16, \log(es)\sqrt{s_0}/(16s)\}$, $\theta_3 = \cdots = \theta_{l_{\max}} = (\log(es)/s)^{1/2}/16$, with $l_{\max} = \lceil C\log(es)\rceil + 2$.
\end{itemize}
We can see the following events happen
\begin{equation}\label{ineq:thm1-6}
\begin{split}
& \|X_{[I_l, T]}^\top X_{[I_l, T]}\Sigma_{T, T}^{-1}/n_l - I_{|T|}\| \leq C\sqrt{st_l/n_l} \leq \sqrt{1/\log(es)}, \quad l=1,2;\\
& \|X_{[I_l, T]}^\top X_{[I_l, T]}\Sigma_{T, T}^{-1}/n_l - I_{|T|}\| \leq C\sqrt{st_l/n_l} \leq 1/2, \quad l=3,\ldots, l_{\max};
\end{split}
\end{equation}
\begin{equation}
\begin{split}
& \max_{j\in G^c}\left\|H_{\|q_{l-1}\|_2/(16\sqrt{s})}\left(X^\top_{[I_l, (j)]}X_{[I_l, T]}\Sigma_{T, T}^{-1}/n_l\cdot q_{l-1}\right)\right\|_2 \leq \sqrt{s_0}\|q_{l-1}\|_2/(16\sqrt{s}), \quad l=1, 2;\\
& \max_{j\in G^c}\left\|H_{\|q_{l-1}\|_2\cdot \log(es)/(16\sqrt{s})}\left(X^\top_{[I_l, (j)]}X_{[I_l, T]}\Sigma_{T, T}^{-1}/n_l\cdot q_{l-1}\right)\right\|_2 \\
\leq & \sqrt{s_0}\|q_{l-1}\|_2\log(es)/(16\sqrt{s}), \quad l=3,\ldots, l_{\max};\\
\end{split}
\end{equation}
\begin{equation}
\begin{split}
& \left\|X^\top_{[I_l, (G)\backslash T]} X_{[I_l, T]}\Sigma_{T, T}^{-1}/n_l\cdot q_{l-1}\right\|_\infty \leq \|q_{l-1}\|_2/(16\sqrt{s}),\quad l=1, 2\\
& \left\|X^\top_{[I_l, (G)\backslash T]} X_{[I_l, T]}\Sigma_{T, T}^{-1}/n_l\cdot q_{l-1}\right\|_\infty \leq \|q_{l-1}\|_2\cdot (\log(es)/s)^{1/2}/16,\quad l=3, \ldots, l_{\max}.
\end{split}
\end{equation}
with probability at least $1 - C\log(es)\exp(-c\frac{n}{\log(es)}) - C\log(es)\exp\left(-c\frac{n}{s_g}\right) - C\log(es)\exp\left(-c\frac{n}{s}\right)$. By triangle inequality, $u_0$ satisfies 
\begin{equation}\label{ineq:thm1-7}
\begin{split}
\|u_0\|_2 \leq & \sqrt{s/s_g}\left(\sum_{j\in G}\left\|\frac{\beta^*_{T,(j)}}{\|\beta^*_{T,(j)}\|_2}\right\|_2^2\right)^{1/2} + \|\sgn(\beta^*_T)\|_2 \leq 2\sqrt{s}.
\end{split}
\end{equation}
When $\max_{i \in T^c}\|X_T^\top X_i/n\|_2 \leq \frac{1}{2}$ and \eqref{ineq:thm1-6}-\eqref{ineq:thm1-7} hold, we have
\begin{equation}\label{ineq:q-l-2-norm}
\begin{split}
 \|q_0\|_2 &\leq 2\sqrt{s},\\
 \|q_1\|_2 &= \left\|(I_{|T|} - X_{I_1, T}^\top X_{I_1, T}\Sigma_{T, T}^{-1}/n_1)q_0\right\| \leq \|I_{|T|} - X_{I_1, T}^\top X_{I_1, T}\Sigma_{T, T}^{-1}/n_1\|\cdot \|q_0\|_2\\ &\leq 2\sqrt{s/\log(es)};\\
\text{similarly, } & \|q_2\|_2 \leq \|q_1\|_2/\sqrt{\log(es)} \leq 2\sqrt{s}/(\log(es));\\
& \|q_l\|_2 \leq \|q_{l-1}\|_2/2 \leq \cdots \leq \|q_2\|/2^{l-2} \leq 2^{3-l}\sqrt{s}/(\log(es)), \quad l\geq 3.
\end{split}
\end{equation}
For large constant $C>0$,$\|q_{l_{\max}}\|_2 \leq 2^{3- C\log(es)}\sqrt{s}/\log(es) \leq c_{\min}/8$. Notice that
\begin{equation*}
	u_T = (\sum_{l = 1}^{l_{\max}}\gamma_l)_T = (\sum_{l = 1}^{l_{\max}}(\alpha_{l - 1} - \alpha_{l}))_T = (\alpha_0 - \alpha_{l_{\max}})_T = (\widetilde{u}_0)_T - (q_{l_{\max}})_T, 
\end{equation*}
we know that
\begin{equation*}
    \|u_T - (\widetilde{u}_0)_T\|_2\cdot \max_{i\in T^c}\|X_T^\top X_i/n\|_2 = \|q_{l_{\max}}\|_2\cdot \max_{i\in T^c}\|X_T^\top X_i/n\|_2 \leq \frac{c_{\min}}{8}\cdot\frac{1}{2} < \frac{c_{\min}}{8}.
\end{equation*}
In addition,
\begin{equation*}
\begin{split}
\|u_{(G)\backslash T}\|_\infty \leq & \sum_{l=1}^{l_{\max}} \left\|X_{[I_l, (G)\backslash T]}^\top X_{[I_l, T]}\Sigma_{T, T}^{-1}/n_l\cdot (\alpha_{l-1})_T\right\|_\infty \\
\leq & \|q_0\|_2/(16\sqrt{s}) + \|q_1\|_2/(16\sqrt{s}) + \sum_{l=3}^{l_{\max}} \|q_{l-1}\|_2\cdot (\log(es)/s)^{1/2}/16 \\
\leq & 1/8 + 1/8 + \sum_{l=3}^\infty 2^{4-l}/16 \leq 1/2.
\end{split}
\end{equation*}
Since
\begin{equation*}
	\begin{split}
	    &\|q_0\|_2/(16\sqrt{s}) + \|q_1\|_2/(16\sqrt{s}) + \sum_{l=3}^{l_{\max}} \|q_{l-1}\|_2\cdot \log(es)/(16\sqrt{s})\\
	    \leq& \frac{1}{8} + \frac{1}{8} + \sum_{l=3}^{l_{\max}}2^{4-l}\sqrt{s}/(\log(es))\cdot \log(es)/(16\sqrt{s}) \leq \frac{1}{2},
	\end{split}
\end{equation*}
\begin{equation*}
\begin{split}
\left\|H_{1/2}(u_{(G^c)})\right\|_{\infty, 2}
\leq & \|H_{\|q_0\|_2/(16\sqrt{s}) + \|q_1\|_2/(16\sqrt{s}) + \sum_{l=3}^{l_{\max}} \|q_{l-1}\|_2\cdot \log(es)/(16\sqrt{s})}(u_{(G^c)})\|_{\infty, 2}\\
\leq & \sum_{l=1}^{2} \left\|H_{\|q_{l - 1}\|_2/(16\sqrt{s})}\left(X^\top_{[I_l, (G^c)]}X_{[I_l, T^c]}q_{l-1}\right)\right\|_{\infty, 2}\\ &+ \sum_{l=3}^{l_{\max}} \left\|H_{\|q_{l - 1}\|_2\cdot \log(es)/(16\sqrt{s})}\left(X^\top_{[I_l, (G^c)]}X_{[I_l, T^c]}q_{l-1}\right)\right\|_{\infty, 2}\\
\leq& \sum_{l = 1}^{2}\sqrt{s_0}\|q_{l-1}\|_2/(16\sqrt{s}) + \sum_{l=3}^{l_{\max}} \sqrt{s_0}\|q_{l-1}\|_2\cdot \log(es)/(16\sqrt{s})\\
\leq & \sqrt{s_0}/2.
\end{split}
\end{equation*}
Thus, the construction of $u$ satisfies all required condition in Lemma \ref{lm:approximate-dual-certificate} with probability at least $1 - C\exp\left(-c\frac{n}{s}\right)$. This has finished the proof of this lemma. 
\quad $\square$

\subsection{Proof of Lemma \ref{lm:weakRIP}}
Let $g(S)$ be the group support of set $S$, that is, $g(S) = \{i_1, \dots, i_k\}$ if $S \subset \cup_{j = 1}^k(i_j)$ and $S \cap (i_j)$ are not empty for all $1 \leq j \leq k$. Lemma \ref{lm:X-random-matrix-properties} Part 1 and the union bound show that 
\begin{equation*}
\begin{split}
&\bbP\left(\exists \gamma \in \bbR^p, \|\gamma\|_0 \leq 2s, \|\gamma\|_{0,2} \leq 2s_g, \frac{1}{n}\|X\gamma\|_2^2 \notin \left[\frac{c_{\min}}{2}\|\gamma\|_2^2, (C_{\min} + \frac{c_{\min}}{2})\|\gamma\|_2^2\right] \right)\\
=&\bbP\left(\exists x \in \bbR^{2s \wedge p}, S \subseteq {1, \cdots, p}, |S| = 2s \wedge p, |g(S)| \leq 2s_g, \frac{1}{n}\|X_{S}x\|_2^2 \notin \left[\frac{c_{\min}}{2}\|\gamma\|_2^2, (C_{\min} + \frac{c_{\min}}{2})\|\gamma\|_2^2\right]\right)\\
\leq& \sum_{S \subseteq \{1, \dots, p\}, |S| = 2s \wedge p, |g(S)| \leq 2s_g}\bbP\left(\forall x \in \bbR^{2s \wedge p}, \frac{1}{n}\|X_{S}x\|_2^2 \notin \left[\frac{c_{\min}}{2}\|\gamma\|_2^2, (C_{\min} + \frac{c_{\min}}{2})\|\gamma\|_2^2\right]\right)\\
\leq& \sum_{S \subseteq \{1, \dots, p\}, |S| = 2s \wedge p, |g(S)| \leq 2s_g}\bbP\left(\|\frac{1}{n}X_{S}^\top X_{S} - \Sigma_{S, S}\| \geq \frac{c_{\min}}{2}\right)\\
\leq& \sum_{S \subseteq \{1, \dots, p\}, |S| = 2s \wedge p, |g(S)| \leq 2s_g}\bbP\left(\|\frac{1}{n}X_{S}^\top X_{S}\Sigma_{S, S}^{-1} - I_{|S|}\| \geq \frac{c_{\min}}{2C_{\max}}\right)\\
\leq& \left[\binom{d}{2s_g} \vee 1\right]\binom{2s_gb}{2s}\cdot 2\exp\left(Cs - cn\right)\\
\leq& \left(\frac{ed}{2s_g}\right)^{2s_g}\left(\frac{e\cdot 2s_gb}{2s}\right)^{2s}\cdot 2\exp\left(Cs - cn\right)
\\ \leq& 2\exp\left(2s\log(es_gb/s) + 2s_g\log(ed/s_g) + Cs - cn\right)\\ \leq& 2e^{-cn}.
\end{split}
\end{equation*}
\quad $\square$

\subsection{Proof of Theorem \ref{th:sample_complexity}}
If $d \geq 3s_g$ and $b \geq 3s/s_g$, by \eqref{eq25}, we can find $\Omega^{(1)}, \dots, \Omega^{(N)} \subset \{1, \dots, db\}$ such that $|\Omega^{(i)}| = s_g\lfloor s/s_g\rfloor$, $|\Omega_{(k)}^{(i)}| = \lfloor s/s_g\rfloor1_{\{\Omega_{(k)}^{(i)} \text{ is not empty}\}}$ for all $1 \leq i \leq N$, $1 \leq k \leq d$, and
\begin{equation}\label{eq33}
	\left|\Omega^{(i)} \cap \Omega^{(j)}\right| \leq 8s_g\lfloor s/s_g\rfloor/9, \quad 1 \leq i \neq j \leq N,
\end{equation}
\begin{equation}\label{eq34}
	\left|\left\{k|\left|\Omega_{(k)}^{(i)} \cap \Omega_{(k)}^{(j)}\right| \geq 2\lfloor s/s_g\rfloor/3\right\}\right| \leq 2s_g/3, \quad 1 \leq i \neq j \leq N,
\end{equation}
where $N = \left\lfloor \left(\frac{d}{2\sqrt{2}s_g}\right)^{s_g/3}\left(\frac{b}{2\sqrt{2}\lfloor s/s_g\rfloor}\right)^{s/9}\right\rfloor$. 
For any $1 \leq j \leq db, 1 \leq i \leq N$, define
\begin{equation*}
	\beta^{(i)}_j = \left\{ \begin{array}{ll}
\frac{1}{\lambda s_g\lfloor s/s_g\rfloor + \lambda_g s_g\sqrt{\lfloor s/s_g\rfloor}}, & j \in \Omega^{(i)}\\
0 & j \notin \Omega^{(i)},
\end{array} \right.
\end{equation*}
then $\|\beta^{(i)}\|_0 \leq s, \|\beta^{(i)}\|_{0, 2} \leq s_g$.
We consider the quotient space
\begin{equation*}
	\bbR^{db}/\text{ker}(X) = \left\{[x] := x + \text{ker}(X), x \in \bbR^{db}\right\}.
\end{equation*}
Then the dimension of $\bbR^{db}/\text{ker}(X)$ is rank$(X) \leq n$. Define the norm $\|[x]\| = \inf_{v \in \text{ker}(X)}\{\lambda\|x - v\|_1 + \lambda_g\|x - v\|_{1, 2}\}$. For any vector $x \in \bbR^{db}$ satisfying $\|x\|_0 \leq 2s, \|x\|_{0, 2} \leq 2s_g$, note that $x - v$ with $v \in \text{ker}(X)$ satisfies $X(x - v) = Xx$, by our assumption, we have $\|[x]\| = \lambda\|x\|_1 + \lambda_g\|x\|_{1, 2}$. Thus $\|[\beta^{(1)}]\| = \cdots = \|[\beta^{(N)}]\| = 1$. Moreover, by \eqref{eq33} and \eqref{eq34},
\begin{equation*}
	\begin{split}
	\|\beta^{(i)} - \beta^{(j)}\|_1 =& \frac{1}{\lambda s_g\lfloor s/s_g\rfloor + \lambda_g s_g\sqrt{\lfloor s/s_g\rfloor}}\left(|\Omega^{(i)}| + |\Omega^{(j)}| - 2|\Omega^{(i)} \cap \Omega^{(j)}|\right)\\ \geq& \frac{2s_g\lfloor s/s_g\rfloor}{9(\lambda s_g\lfloor s/s_g\rfloor + \lambda_g s_g\sqrt{\lfloor s/s_g\rfloor})},
	\end{split}
\end{equation*}
and
\begin{equation*}
    \begin{split}
    \|\beta^{(i)} - \beta^{(j)}\|_{1,2}
    =& \sum_{k = 1}^{d}\|\beta^{(i)}_{(k)} - \beta^{(j)}_{(k)}\|_{2}\\
    \geq& \sum_{k \in S_{i,j}}\|\beta^{(i)}_{(k)} - \beta^{(j)}_{(k)}\|_{2}\\
    \geq& \frac{1}{\lambda s_g\lfloor s/s_g\rfloor + \lambda_g s_g\sqrt{\lfloor s/s_g\rfloor}}\sqrt{\frac{2\lfloor s/s_g\rfloor}{3}}\cdot |S_{i,j}|\\
    \geq& \frac{1}{\lambda s_g\lfloor s/s_g\rfloor + \lambda_g s_g\sqrt{\lfloor s/s_g\rfloor}}\sqrt{\frac{2\lfloor s/s_g\rfloor}{3}}\cdot \frac{s_g}{3},
    \end{split}
\end{equation*}
where $S_{i,j} = \left\{k|\Omega_{(k)}^{(i)}, \Omega_{(k)}^{(j)} \text{ are not empty sets}, \left|\Omega_{(k)}^{(i)} \cap \Omega_{(k)}^{(j)}\right| < 2\lfloor s/s_g\rfloor/3\right\}$.\\
Since $\beta^{(i)} - \beta^{(j)}$ is $(2s, 2s_g)$-sparse,
\begin{equation*}
	\begin{split}
	\left\|[\beta^{(i)}] - [\beta^{(j)}]\right\| =& \left\|[\beta^{(i)} - \beta^{(j)}]\right\| = \lambda\left\|\beta^{(i)} - \beta^{(j)}\right\|_1 + \lambda_g\left\|\beta^{(i)} - \beta^{(j)}\right\|_{1, 2} \geq 2/9.
	\end{split}
\end{equation*}
By \cite[Proposition C.3]{foucart2013mathematical}, we have $N \leq 10^{\text{rank}(X)} \leq 10^n$. Therefore we have
\begin{equation*}
	\left\lfloor \left(\frac{d}{2\sqrt{2}s_g}\right)^{s_g/3}\left(\frac{b}{2\sqrt{2}\lfloor s/s_g\rfloor}\right)^{s/9}\right\rfloor \leq 10^n, 
\end{equation*}
which means that $n \geq c(s_g\log(d/s_g) + s\log(es_gb/s))$.

If $d < 3s_g$ or $b < 3s/s_g$, let $s_g'=[s_g/3] \vee 1 \geq s_g/5, s' = [s/15] \vee s_g'$, then $d \geq 3s_g'$ and $b \geq 3s'/s_g'$. Since all $(2s, 2s_g)$-sparse vectors can be exactly recovered by the $\ell_1 + \ell_{1,2}$ minimization and $s' \leq s, s_g' \leq s_g$, we know that the $\ell_1 + \ell_{1,2}$ minimization exactly recover all $(2s', 2s_g')$-sparse vectors. Therefore, we have 
\begin{equation}\label{ineq1}
	\begin{split}
	n \geq& c(s_g'\log(d/s_g') + s'\log(es_g'b/s')) \geq c\bigg(\frac{s_g}{5}\cdot \log\bigg(\frac{d}{s_g}\bigg) + \frac{s}{15}\cdot \log\bigg(\frac{eb(s_g/5)}{s/15} \vee eb\bigg)\bigg)\\ \geq& c'(s_g\log(d/s_g) + s\log(es_gb/s)).
	\end{split}
\end{equation}
$\square$

\subsection{Proof of Theorem \ref{th:lower_noiseless}}
We would like prove Theorem \ref{th:lower_noiseless} by contradiction. Let
\begin{equation*}
	c = \min\left\{\frac{1}{8}, c', \sqrt{\frac{c'}{256}}\right\}, \quad c_0 = \min\left\{\frac{c}{2e}, \frac{c^2}{2C^2}, 16c^2\right\}, \quad C_0 = \max\left\{\frac{C^2}{c^2}, \frac{1}{32c^2}\right\},
\end{equation*}
where $c'$ is a uniform constant such that $n \geq c'(s\log(es_gb/s) + s_g\log(d/s_g))$ if the conditions in Theorem \ref{th:sample_complexity} are satisfied.
Assume for contradiction that
\begin{equation}\label{eq44}
n < c_0(s\log(es_gb/s) + s_g\log(d/s_g)).
\end{equation}
Let $s_0 = s/s_g$, define the norm $\|\cdot\| = \|\cdot\|_1 + \sqrt{s_0}\|\cdot\|_{1,2}$. Let $B = \{x \in \bbR^p| \|x\| \leq 1\}$,
\begin{equation*}
d^n(B, \bbR^p) = \inf_{\substack{L^n \text{ is a subspace of } \bbR^p \\\text{ with dim}(\bbR^p/L^n) \leq n} }\left\{\sup_{\beta \in B \cap L^n}\|\beta\|_2\right\}.
\end{equation*}
By \cite[Theorem 10.4]{foucart2013mathematical}, we have
\begin{equation}\label{eq43}
d^n(B, \bbR^p) \leq \sup_{\beta \in B}\|\beta - \Delta(X\beta)\|_2 \leq \frac{C}{\sqrt{s}}\sup_{\beta \in B}\left(\|\beta\|_1 + \sqrt{s_0}\|\beta\|_{1,2}\right) = \frac{C}{\sqrt{s}}.
\end{equation}
If
\begin{equation}\label{eq42}
d^n(B, \bbR^p) \geq c\min\left\{\frac{1}{\sqrt{s_0}}, \left[\left(\frac{s_g}{s}\log\left(\frac{c\frac{s}{s_g}d\log(es_gb/s)}{n}\right) + \log(es_gb/s)\right)/n\right]^{1/2}\right\},
\end{equation}
since
\begin{equation*}
	\frac{C}{\sqrt{s}} \leq \frac{c\sqrt{C_0}}{\sqrt{s}} \leq \frac{c\sqrt{s_g}}{\sqrt{s}} = \frac{c}{\sqrt{s_0}},
\end{equation*}
\eqref{eq43} and \eqref{eq42} together imply that
\begin{equation}\label{eq46}
n \geq 	\frac{c^2}{C^2}\left(s_g\log\left(\frac{c\frac{s}{s_g}d\log(es_gb/s)}{n}\right) + s\log(es_gb/s)\right).
\end{equation}
By \eqref{eq44}, 
\begin{equation}\label{eq45}
\begin{split}
\frac{c\frac{s}{s_g}d\log(es_gb/s)}{n} >& \frac{c\frac{s}{s_g}d\log(es_gb/s)}{c_0(s\log(es_gb/s) + s_g\log(d/s_g))}\\\geq& 2e\frac{\frac{s}{s_g}d\log(es_gb/s)}{s\log(es_gb/s) + s_g\log(d/s_g)}\\ \geq& \min\left\{e\frac{\frac{s}{s_g}d\log(es_gb/s)}{s\log(es_gb/s)}, \frac{e\frac{s}{s_g}d\log(es_gb/s)}{s_g\log(ed/s_g)}\right\}\\
\geq & \min\left\{\frac{ed}{s_g}, \frac{\frac{ed}{s_g}}{\log(\frac{ed}{s_g})}\right\} \geq \left(\frac{ed}{s_g}\right)^{1/2}.
\end{split}
\end{equation}
In the last inequality, we used $x^{1/2} \geq \log(x)/2$ for all $x \geq 1$.\\
Combine \eqref{eq46} and \eqref{eq45} together, we have
\begin{equation*}
n \geq \frac{c^2}{2C^2}\left(s\log(es_gb/s) + s_g\log(d/s_g)\right) \geq c_0\left(s\log(es_gb/s) + s_g\log(d/s_g)\right) > n,
\end{equation*}
contradiction!

Thus, we only need to prove \eqref{eq42} based on \eqref{eq44}. We still use the proof of contradiction. If
\begin{equation*}
	d^n(B, \bbR^p) < c\min\left\{\frac{1}{\sqrt{s_0}}, \left[\left(\frac{s_g}{s}\log\left(\frac{c\frac{s}{s_g}d\log(es_gb/s)}{n}\right) + \log(es_gb/s)\right)/n\right]^{1/2}\right\} := \mu,
\end{equation*}
then there exists a subspace $L^n$ of $\bbR^p$ with dim$(\bbR^p/L^n) \leq n$ such that for all $v \in L^n\backslash \{0\}$, 
\begin{equation*}
\|v\|_2 < \mu\left(\|v\|_1 + \sqrt{s_0}\|v\|_{1,2}\right).
\end{equation*}
Let $B \in \bbR^{n \times p}$ satisfying ker$(B) = L^n$. Let $s' = \lfloor \frac{1}{32\mu^2}\rfloor, s_g' = \lfloor s'/s_0\rfloor$, by \eqref{eq44} and \eqref{eq45},
\begin{equation*}
\frac{1}{8}s_0^{-1/2} \geq cs_0^{-1/2} \geq \mu \geq c\min\left\{\sqrt{\frac{C_0}{s}}, \left(\frac{\frac{s_g}{2s}\log(d/s_g) + \log(es_gb/s)}{c_0(s_g\log(d/s_g) + s\log(es_gb/s))}\right)^{1/2}\right\} \geq \frac{1}{4\sqrt{2}}s^{-1/2}, 
\end{equation*}
which means that 
\begin{equation*}
1 \leq s' \leq s, \quad 1 \leq s_g' \leq s_g.
\end{equation*}
Moreover, we have $\frac{1}{64\mu^2} < s' \leq \frac{1}{32\mu^2}$. For any $(2s', 2s_g')$-sparse $\beta$ with support set $T$ and group support set $G$, and $v \in \text{ker}(A)$, by Cauchy-Schwarz inequality,
\begin{equation*}
\begin{split}
\|v_T\|_1 + \sqrt{s_0}\|v_{(G)}\|_{1,2} \leq& \sqrt{2s'}\|v_T\|_2 + \sqrt{s_0}\sqrt{2s_g'}\|v_T\|_2 \leq 2\sqrt{2s'}\|v_T\|_2\\ <& 2\sqrt{2}\frac{1}{4\sqrt{2}\mu}\mu\left(\|v\|_1 + \sqrt{s_0}\|v\|_{1,2}\right)
= \frac{1}{2}\left(\|v\|_1 + \sqrt{s_0}\|v\|_{1,2}\right),
\end{split}
\end{equation*}
i.e.,
\begin{equation*}
\|v_T\|_1 + \sqrt{s_0}\|v_{(G)}\|_{1,2} < \|v_{T^c}\|_1 + \sqrt{s_0}\|v_{(G^c)}\|_{1,2}.
\end{equation*}
Based on Cauchy-Schwarz inequality and the sub-differential of $\|\beta\|_1$ and $\|\beta\|_{1,2}$, we have
\begin{equation*}
\begin{split}
&\|\beta + v\|_1 + \sqrt{s_0}\|\beta + v\|_{1,2}\\
\geq& \|\beta\|_1 + \text{sgn}(\beta)^\top v_T + \|v_{T^c}\|_1 + \sqrt{s_0}\left(\|\beta\|_{1,2} + \sum_{j \in G}\frac{\beta^\top_{(j)}v_{(j)}}{\|\beta_{(j)}\|_2} + \sum_{j \in G^c}\|v_j\|_2\right)\\
\geq& \|\beta\|_1 - \|v_T\|_1 + \|v_{T^c}\|_1 + \sqrt{s_0}\left(\|\beta\|_{1,2} - \|v_{(G)}\|_2 + \|v_{(G^c)}\|_2\right)\\
>& \|\beta\|_1 + \sqrt{s_0}\|\beta\|_{1,2}.
\end{split}
\end{equation*}
By Theorem \ref{th:sample_complexity}, 
\begin{equation*}
n \geq c'(s'\log(es_g'b/s') + s_g'\log(d/s_g')) \geq c's'\left\{\log\left(\frac{es_gb}{2s}\right) + \frac{1}{2s_0}\log(s_0d/s')\right\} \geq c's'\log\left(\frac{es_gb}{2s}\right).
\end{equation*}
Thus
\begin{equation*}
\begin{split}
n \geq& c's'\left(\log\left(\frac{es_gb}{2s}\right) + \frac{s_g}{s}\log\left(\frac{c'\frac{s}{s_g}d\log(es_gb/s)}{n}\right)\right)\\ >& \frac{c'}{64\mu^2}\left(\frac{1}{4}\log(es_gb/s) + \frac{s_g}{s}\log\left(\frac{c\frac{s}{s_g}d\log(es_gb/s)}{n}\right)\right)\\ \geq& n
\end{split}
\end{equation*}
provided that $c = \min\left\{\frac{1}{8}, c', \sqrt{\frac{c'}{256}}\right\}$, 
contradiction! This means that \eqref{eq42} holds if \eqref{eq44} is true.

Therefore, we have finished the proof of Theorem \ref{th:lower_noiseless}.
\quad $\square$

\subsection{Proof of Theorem \ref{th:noisy}}

Let $\lambda = C\sigma\sqrt{\frac{s\log(es_gb) + s_g\log(d/s_g)}{s}n}, \lambda_g = \sqrt{s/s_g}\lambda$. By \eqref{ineq:particular_lemma3} in Lemma \ref{lm:bernstein-sub-gaussian} and \eqref{guassian tail}, one has
\begin{equation}\label{eq6}
\begin{split}
    &\bbP\left(\left\|H_{\frac{1}{10}\lambda}(X^\top \varepsilon)\right\|_{\infty, 2} \geq \frac{1}{10}\lambda_g\right)\\ \leq& \bbP\left(\exists 1 \leq j \leq d, \left\|H_{\frac{1}{10}\lambda}(X_{(j)}^\top \varepsilon)\right\|_{2} \geq \frac{1}{10}\lambda_g, \|\varepsilon\|_2 \geq 5\sqrt{n\sigma^2}\right) + \bbP\left(\|\varepsilon\|_2 \geq 5\sqrt{n\sigma^2}\right)\\
    \leq& \bbP\left(\exists 1 \leq j \leq d, \left\|H_{\frac{1}{10}\lambda}(X_{(j)}^\top \varepsilon)\right\|_{2} \geq \frac{1}{10}\lambda_g\bigg| \|\varepsilon\|_2 \geq 5\sqrt{n\sigma^2}\right) + \bbP\left(\|\varepsilon\|_2 \geq 5\sqrt{n\sigma^2}\right)\\
    \leq& d\exp\left(-C\frac{s\log(es_gb) + s_g\log(d/s_g)}{s_g}\right) + e^{-n}\\
    =& \exp\left(\log(s_g) + \log(d/s_g) - C\frac{s\log(es_gb) + s_g\log(d/s_g)}{s_g}\right) + e^{-n}\\
    \leq& \exp\left(-C\frac{s\log(es_gb) + s_g\log(d/s_g)}{s_g}\right) + e^{-n}.
\end{split}
\end{equation}
By the definition of $\hat{\beta}$ and KKT condition, we have
\begin{equation*}
	X^\top(y - X\hat\beta) + \lambda z_1 + \lambda_g z_2 = 0, 
\end{equation*}
where 
\begin{equation*}
	\left\{\begin{array}{ll}
	(z_1)_i  = \text{sgn}(\hat\beta_i), & \hat\beta_i \neq 0;\\
	|(z_1)_i| \leq 1, & \hat\beta_i = 0; \end{array}\right. \quad
	\left\{\begin{array}{ll}
	(z_2)_{(j)}  = \frac{\hat{\beta}_{(j)}}{\|\hat{\beta}_{(j)}\|_2}, & \hat{\beta}_{(j)} \neq 0;\\
	\|(z_2)_{(j)}\|_2 \leq 1, & \hat{\beta}_{(j)} = 0. \end{array}\right. 
\end{equation*}
Therefore,
\begin{equation*}
	\|H_{\lambda}(X^\top(X\hat{\beta} - y))\|_{\infty, 2} \leq \lambda_g.
\end{equation*}
\eqref{eq6}, Lemma \ref{lm:soft-thresholding} Part 1 and the previous inequality together imply that 
\begin{equation}\label{eq7}
	\bbP\left(\left\|H_{(1 + \frac{1}{10})\lambda}(X^\top Xh)\right\|_{\infty, 2} \leq (1 + \frac{1}{10})\lambda_g\right) \geq 1 - \exp\left(-C\frac{s\log(es_gb) + s_g\log(d/s_g)}{s_g}\right) - e^{-n},
\end{equation}
where $h = \hat{\beta} - \beta^*$.
By the definition of $\hat{\beta}$, we have
\begin{equation*}
	\|y - X\hat\beta\|_2^2 + \lambda\|\hat\beta\|_1 + \lambda_g\|\hat{\beta}\|_{1, 2} \leq \|y - X\beta^*\|_2^2 + \lambda\|\beta^*\|_1 + \lambda_g\|\beta^*\|_{1, 2}.
\end{equation*}
\eqref{ineq:lm-dual-1} and the previous inequality show that 
\begin{equation}\label{eq15}
	\begin{split}
	&\|Xh\|_2^2 + \lambda\|h_{T^c}\|_1 + \lambda_g\|h_{(G^c)}\|_{1, 2}\\ \leq& 2\langle Xh, \varepsilon\rangle - \lambda\cdot \text{sgn}(\beta^*_T)^\top h_T - \lambda_g\sum_{j \in G}\frac{\beta_{T,(j)}^{*\top}h_{(j)}}{\|\beta^*_{T,(j)}\|_2} + 2\lambda\|\beta^*_{T^c}\|_1 + 2\lambda_g\|\beta^*_{T^c}\|_{1,2}. 
	\end{split}
\end{equation}

First, we consider $\langle Xh, \varepsilon\rangle$. Denote $P = X_T(X_T^\top X_T)^{-1}X_T^\top$, since $Xh = X_Th_T + X_{T^c}h_{T^c}$ and $(I_n - P)X_T = 0$,
\begin{equation}\label{eq12}
	\begin{split}
	&\left|\langle Xh, \varepsilon\rangle\right| \leq \left|\langle PXh, \varepsilon\rangle\right| + \left|\langle (I_n - P)Xh, \varepsilon\rangle\right|\\ =& \left|\langle X_T^\top Xh, (X_T^\top X_T)^{-1}X_T^\top\varepsilon\rangle\right| + \left|\langle (I_n - P)X_{T^c}h_{T^c}, \varepsilon\rangle\right|\\
	=& \left|\langle X_T^\top Xh, (X_T^\top X_T)^{-1}X_T^\top\varepsilon\rangle\right| + \left|\langle X_{T^c}h_{T^c}, (I_n - P)\varepsilon\rangle\right|.
	\end{split}
\end{equation}
Therefore, to give an upper bound of $\left|\langle Xh, \varepsilon\rangle\right|$, we only need to bound $\left|\langle X_T^\top Xh, (X_T^\top X_T)^{-1}X_T^\top\varepsilon\rangle\right|$ and $\left|\langle X_{T^c}h_{T^c}, (I_n  - P)\varepsilon\rangle\right|$, respectively.
By Part 1 of Lemma \ref{lm:X-random-matrix-properties} and also notice that $c_{\min} \leq \sigma_{\min}(\Sigma) \leq \sigma_{\max}(\Sigma) \leq C_{\max}$,
\begin{equation}\label{eq8}
	\begin{split}
	\bbP\left(\left\|\left(\frac{1}{n}X_T^\top X_T\right)^{-1}\right\| \geq \frac{2}{c_{\min}}\right) \leq& \bbP\left(\|\frac{1}{n}X_T^\top X_T - \Sigma_{T, T}\| \geq \frac{c_{\min}}{2}\right)\\ \leq& \bbP\left(\|\frac{1}{n}X_T^\top X_T\Sigma_{T, T}^{-1} - I_s\| \geq \frac{c_{\min}}{2C_{\max}}\right)\\ \leq& 2\exp\left(Cs - cn\right) \leq 2\exp\left(-cn\right).
	\end{split}
\end{equation}
\eqref{eq8}, Lemma \ref{lm:infinity norm} and Cauchy-Schwarz inequality together imply that with probability at least $1 - \exp\left(-C\frac{s\log(es_gb)+s_g\log(d/s_g)}{s}\right) - 2\exp(-cn)$, 
\begin{equation*}
	\|(X_T^\top X_T)^{-1}X_T^\top\varepsilon\|_1 \leq \frac{2}{c_{\min}}\frac{\sqrt{s}}{n}\|X_T^\top\varepsilon\|_2 \leq \frac{2}{c_{\min}}\frac{s}{n}\|X_T^\top\varepsilon\|_{\infty} \leq C\frac{s}{n}\sqrt{n\frac{s\log(es_gb)+s_g\log(d/s_g)}{s}\sigma^2} \leq C\frac{s}{n}\lambda,
\end{equation*}
\begin{equation*}
	\|(X_T^\top X_T)^{-1}X_T^\top\varepsilon\|_{1, 2} \leq \sqrt{s_g}\|(X_T^\top X_T)^{-1}X_T^\top\varepsilon\|_2 \leq \frac{2}{c_{\min}}\frac{\sqrt{s_g}}{n}\|X_T^\top\varepsilon\|_2 \leq C\frac{\sqrt{s\cdot s_g}}{n}\lambda.
\end{equation*}
Combine Lemma \ref{lm:soft-thresholding} Part 2, \eqref{eq7} and the previous two inequalities together, with probability at least $1 - 2\exp\left(-C\frac{s\log(es_gb)+s_g\log(d/s_g)}{s}\right) - 3e^{-cn}$, 
\begin{equation}\label{eq9}
	\begin{split}
	\left|\langle X_T^\top Xh, (X_T^\top X_T)^{-1}X_T^\top\varepsilon\rangle\right| \leq &\frac{11}{10}\lambda\|(X_T^\top X_T)^{-1}X_T^\top\varepsilon\|_1 + \frac{11}{10}\lambda_g\|(X_T^\top X_T)^{-1}X_T^\top\varepsilon\|_{1, 2}\\
	\leq& C\frac{s}{n}\lambda^2.
	\end{split}
\end{equation}
Similarly to the proof of \eqref{eq6}, also notice that $\|(I_n - P)\varepsilon\|_2 \leq \|\varepsilon\|_2$ and $X_{(G^c)}$ is independent of $I_n - P$, we have
\begin{equation*}
	\begin{split}
	&\bbP\left(\left\|H_{\frac{1}{10}\lambda}\left(X_{(G^c)}^\top(I_n - P) \varepsilon\right)\right\|_{\infty, 2} \geq \frac{1}{10}\lambda_g\right)\\
	\leq& \bbP\left(\exists j \in G^c, \left\|H_{\frac{1}{10}\lambda}\left(X_{(j)}^\top(I_n - P)\varepsilon\right)\right\|_2 \geq \frac{1}{10}\lambda_g\bigg|\|(I_n - P)\varepsilon\|_2 \geq 5\sqrt{n\sigma^2}\right)\\ &+ \bbP\left(\|\varepsilon\|_2 \geq 5\sqrt{n\sigma^2}\right)\\
	\leq& \exp\left(-C\frac{s\log(es_gb) + s_g\log(d/s_g)}{s_g}\right) + e^{-n}.
	\end{split}
\end{equation*}
By Lemma \ref{lm:soft-thresholding} Part 2 and \eqref{eq6}, with probability at least $1 - \exp\left(-C\frac{s\log(es_gb) + s_g\log(d/s_g)}{s_g}\right) - e^{-n}$,
\begin{equation*}
	\left|\langle X_{(G^c)}h_{(G^c)}, (I_n - P)\varepsilon\rangle\right| = \left|\langle h_{(G^c)}, X_{(G^c)}^\top(I_n - P)\varepsilon\rangle\right| \leq \frac{1}{10}\lambda\|h_{(G^c)}\|_1 + \frac{1}{10}\lambda_g\|h_{(G^c)}\|_{1, 2}.
\end{equation*}
Notice that $X_{T^c\backslash (G^c)}$ and $I_n - P$ are independent and $|T^c\backslash (G^c)| \leq |G| \leq s_gb$, by Lemma \ref{lm:infinity norm}, with probability at least $1  - \exp(-C\frac{s\log(es_gb) + s_g\log(d/s_g)}{s}) - e^{-n}$,
\begin{equation*}
	\begin{split}
	\left|\langle X_{T^c\backslash (G^c)}h_{T^c\backslash(G^c)}, (I_n - P)\varepsilon\rangle\right| &\leq \|h_{T^c\backslash(G^c)}\|_{1}\|X_{T^c\backslash (G^c)}^\top(I_n - P)\varepsilon\|_{\infty}\\ &\leq C\sqrt{n\frac{s\log(es_gb) + s_g\log(d/s_g)}{s}\sigma^2}\|h_{T^c\backslash(G^c)}\|_1\\ &\leq \frac{1}{10}\lambda\|h_{T^c\backslash(G^c)}\|_1.
	\end{split}
\end{equation*}
Combine the previous two inequalities together, we have
\begin{equation}\label{eq10}
	\begin{split}
	\left|\langle X_{T^c}h_{T^c}, (I_n - P)\varepsilon\rangle\right| \leq& \left|\langle X_{(G^c)}h_{(G^c)}, (I_n - P)\varepsilon\rangle\right| + \left|\langle X_{T^c\backslash (G^c)}h_{T^c\backslash(G^c)}, (I_n - P)\varepsilon\rangle\right|\\ \leq& \frac{1}{10}\lambda\|h_{T^c}\|_1 + \frac{1}{10}\lambda_g\|h_{(G^c)}\|_{1, 2}
	\end{split}	
\end{equation}
with probability $1 - C\exp\left(-C\frac{s\log(es_gb) + s_g\log(d/s_g)}{s}\right) - Ce^{-cn}$. 
\begin{comment}
	By\eqref{ineq:thm1-1} and \eqref{eq36}, with probability at least $1 - C(es_gb)^{-C}$,
	\begin{equation}\label{eq11}
	\begin{split}
	\left|\langle h_{T^c}, X_{T^c}X_T(X_T^\top X_T)^{-1}X_T^\top\varepsilon\rangle\right| \leq& \|h_{T^c}\|_1\|X_{T^c}^\top X_T(X_T^\top X_T)^{-1}X_T^\top\varepsilon\|_{\infty}\\
	=& \|h_{T^c}\|_1\max_{i\in T^c}|X_{i}^\top X_T(X_T^\top X_T)^{-1}X_T^\top\varepsilon|\\
	\leq& \|h_{T^c}\|_1\max_{i\in T^c}\|X_{i}^\top X_T\|_2\|(X_T^\top X_T)^{-1}X_T^\top\varepsilon\|_2\\
	\leq& \|h_{T^c}\|_1\cdot \frac{n}{2} \\
	\leq& \|h_{T^c}\|_1\cdot C\sqrt{\frac{s\log(es_gb)}{n}\sigma^2} \leq \frac{1}{10}\lambda\|h_{T^c}\|_1.
	\end{split}
	\end{equation}
\end{comment}
Combine \eqref{eq12}, \eqref{eq9} and \eqref{eq10} together, we know that with probability at least $1 - C\exp\left(-C\frac{s\log(es_gb) + s_g\log(d/s_g)}{s}\right) - Ce^{-cn}$,
\begin{equation}\label{ineq: inner product}
	|\langle Xh, \varepsilon\rangle| \leq C\frac{s}{n}\lambda^2 + \frac{1}{10}\lambda\|h_{T^c}\|_1 + \frac{1}{10}\lambda_g\|h_{(G^c)}\|_{1, 2}.
\end{equation}

Moreover, by the proof of Theorem \ref{th:noiseless}, with probability at least $1 - C\exp(-cn/s)$, there exists an approximate dual certificate $u \in \bbR^p$ in the row span of $X$ satisfying \eqref{ineq:dual-certificate-u}, and $\|v_T - \text{sgn}(\beta^*_T)\|_2 \leq \frac{1}{8}$, where $v$ is defined in \eqref{eq37}. Similarly to \eqref{ineq:lm-dual-2}, we have
\begin{equation*}
\begin{split}
&\sgn(\beta^*_T)^\top h_T + \sum_{j\in G}\frac{\sqrt{s_0}\beta_{T,(j)}^{*\top} h_{(j)}}{\|\beta^*_{T,(j)}\|_2}\\
\geq&  -\|v_T - \sgn(\beta^*_T)\|_2\cdot \|h_T\|_2 - \|h_{T^c}\|_1/2 - \sqrt{s_0}\|h_{(G^c)}\|_{1,2}/2 + \langle h, u\rangle\\
\geq& -\frac{c_{\min}}{8}\cdot \|h_T\|_2 - \|h_{T^c}\|_1/2 - \sqrt{s_0}\|h_{(G^c)}\|_{1,2}/2 + \langle h, u\rangle.
\end{split}
\end{equation*}
By Lemma \ref{lm:approximate dual certificate}, with probability at least $1 - Ce^{-cn/s}$, $u = X^\top w$ with  $\|w\|_2 \leq C\sqrt{s/n}$. Therefore, with probability at least $1 - Ce^{-cn/s}$, 
\begin{equation*}
	|\langle h, u\rangle| = |\langle Xh, w\rangle| \leq \|Xh\|_2\|w\|_2 \leq C\sqrt{s/n}\|Xh\|_2.
\end{equation*}
The two previous inequalities together imply that
\begin{equation}\label{eq13}
	\sgn(\beta^*_T)^\top h_T + \sum_{j\in G}\frac{\sqrt{s_0}\beta_{T,(j)}^{*\top} h_{(j)}}{\|\beta^*_{T,(j)}\|_2} \geq  -\frac{c_{\min}}{8}\cdot \|h_T\|_2 - \|h_{T^c}\|_1/2 - \sqrt{s_0}\|h_{(G^c)}\|_{1,2}/2 -C\sqrt{s/n}\|Xh\|_2
\end{equation}
with probability at least $1 - Ce^{-cn/s}$.

Combine \eqref{eq15}, \eqref{ineq: inner product} and \eqref{eq13} together, with probability at least $1 - Ce^{-C\frac{s\log(es_gb) + s_g\log(d/s_g)}{s}} - Ce^{-cn/s}$, 
\begin{equation}\label{eq16}
\begin{split}
	&\|Xh\|_2^2 + \frac{3}{10}\lambda\|h_{T^c}\|_1 + \frac{3}{10}\lambda_g\|h_{(G^c)}\|_{1, 2}\\ \leq& C\frac{s}{n}\lambda^2 + \frac{c_{\min}}{8}\lambda\|h_T\|_2 + C\sqrt{s/n}\lambda\|Xh\|_2 + 2\lambda\|\beta^*_{T^c}\|_1 + 2\lambda_g\|\beta^*_{T^c}\|_{1,2}.
\end{split}
\end{equation}
By \eqref{ineq:thm1-1}, \eqref{eq7} and \eqref{eq8}, with probability at least $1 - \exp\left(-C\frac{s\log(es_gb) + s_g\log(d/s_g)}{s_g}\right) - Ce^{-cn}$,
\begin{equation}\label{eq17}
	\begin{split}
	\|h_T\|_2 \leq& \|(X_T^\top X_T)^{-1}\|\|X_T^\top X_Th_T\|_2\\
	\leq& \frac{2}{c_{\min}n}\|X_T^\top Xh - X_T^\top X_{T^c}h_{T^c}\|_2\\
	\leq& \frac{2}{c_{\min}n}\left(\|X_T^\top Xh\|_2 + \|X_T^\top X_{T^c}h_{T^c}\|_2\right)\\
	\leq& \frac{2}{c_{\min}n}\left(\|H_{\frac{11}{10}\lambda}(X_T^\top Xh)\|_2 + \frac{11}{10}\sqrt{s}\lambda + n\sum_{i \in T^c}\|X_T^\top X_i/n\|_2|h_i|\right)\\
	\leq& \frac{2}{c_{\min}n}\left(\sqrt{s_g}\|H_{\frac{11}{10}\lambda}(X_T^\top Xh)\|_{\infty,2} + \frac{11}{10}\sqrt{s}\lambda + n\max_{i \in T^c}\|X_T^\top X_i/n\|_2\|h_{T^c}\|_1\right)\\
	\leq& \frac{2}{c_{\min}n}\left(\sqrt{s_g}\frac{11}{10}\lambda_g + \frac{11}{10}\sqrt{s}\lambda + \frac{n}{2}\|h_{T^c}\|_1\right)\\
	\leq& \frac{5}{c_{\min}}\frac{\sqrt{s}}{n}\lambda + \frac{1}{c_{\min}}\|h_{T^c}\|_1.	
	\end{split}
\end{equation}
The fourth inequality comes from $\|x\|_2 \leq \|H_{\alpha}(x)\|_2 + \sqrt{s}\alpha$ for $x \in \bbR^s$; the fifth inequality holds since $\|X_T^\top Xh\|_{0, 2} \leq s_g$.\\
\eqref{eq16} and \eqref{eq17} together imply that
\begin{equation*}
	\begin{split}
	\|Xh\|_2^2 + \frac{7}{40}\lambda\|h_{T^c}\|_1 + \frac{3}{10}\lambda_g\|h_{(G^c)}\|_{1, 2} \leq& C\frac{s}{n}\lambda^2 + C\sqrt{s/n}\lambda\|Xh\|_2 + 2\lambda\|\beta^*_{T^c}\|_1 + 2\lambda_g\|\beta^*_{T^c}\|_{1,2}
	\end{split}
\end{equation*}
with probability at least $1 - C\exp(-C\frac{s\log(es_gb) + s_g\log(d/s_g)}{s}) - Ce^{-cn/s}$.
Also notice that 
\begin{equation*}
	C\sqrt{s/n}\lambda\|Xh\|_2 \leq \|Xh\|_2^2 + C\frac{s}{n}\lambda^2,
\end{equation*}
with probability at least $1 - C\exp(-C\frac{s\log(es_gb) + s_g\log(d/s_g)}{s}) - Ce^{-cn/s}$,
\begin{equation}\label{ineq:h_{T^c}}
	\|h_{T^c}\|_1 + \sqrt{s_0}\|h_{(G^c)}\|_{1,2} \leq C\left(\frac{s}{n}\lambda + \|\beta^*_{T^c}\|_1 + \sqrt{s_0}\|\beta^*_{T^c}\|_{1,2}\right).
\end{equation}

From the proof of Lemma \ref{lm:approximate-dual-certificate}, we know that \eqref{eq27} and \eqref{eq22} hold with probability at least $1 - 2e^{-cn}$.
By Lemma \ref{lm:soft-thresholding} Part 2 and \eqref{eq7}, with probability at least $1 - \exp\left(-C\frac{s\log(es_gb) + s_g\log(d/s_g)}{s_g}\right) - e^{-n}$,
\begin{equation}\label{eq21}
\begin{split}
\left|\langle X_{\widetilde{T}}h_{\widetilde{T}}, Xh\rangle\right| =& \left|\langle h_{\widetilde{T}}, X_{\widetilde{T}}^\top Xh\rangle\right| \leq \frac{11}{10}\left(\lambda\|h_{\widetilde{T}}\|_1 + \lambda_g\|h_{\widetilde{T}}\|_{1, 2}\right)\\ \leq& \frac{11}{10}\left(\lambda\cdot\sqrt{3s}\|h_{\widetilde{T}}\|_2 + \lambda_g\sqrt{2s_g}\|h_{\widetilde{T}}\|_2\right) \leq 4\lambda\sqrt{s}\|h_{\widetilde{T}}\|_2.
\end{split}
\end{equation}
The second inequality is due to $\|h_{\widetilde{T}}\|_0 \leq 3s, \|h_{\widetilde{T}}\|_{0,2} \leq 2s_g$ and Cauchy-Schwarz inequality.

Combine \eqref{eq27}, \eqref{eq22}, \eqref{ineq:h_{T^c}} and \eqref{eq21} together, with probability  at least $1 - C\exp(-C\frac{s\log(es_gb) + s_g\log(d/s_g)}{s}) - Ce^{-cn/s}$, we have
\begin{equation*}
	\begin{split}
	\frac{c_{\min}}{2}\|h_{\widetilde{T}}\|_2^2 \leq& \frac{1}{n}4\lambda\sqrt{s}\|h_{\widetilde{T}}\|_2 + \sqrt{3}C_{\max}\|h_{\widetilde{T}}\|_2(\sqrt{2}s^{-1/2}\|h_{T^c}\|_1 + s_g^{-1/2}\|h_{(G^c)}\|_{1, 2})\\
	\leq& \frac{1}{n}4\lambda\sqrt{s}\|h_{\widetilde{T}}\|_2 + \sqrt{3}C_{\max}\|h_{\widetilde{T}}\|_2\cdot \frac{C}{\sqrt{s}}\left(\frac{s}{n}\lambda + \|\beta^*_{T^c}\|_1 + \sqrt{s_0}\|\beta^*_{T^c}\|_{1,2}\right)\\
	\leq& C\left(\frac{\sqrt{s}}{n}\lambda + \frac{1}{\sqrt{s}}\|\beta^*_{T^c}\|_1 + \frac{1}{\sqrt{s_g}}\|\beta^*_{T^c}\|_{1,2}\right)\|h_{\widetilde{T}}\|_2.
	\end{split}
\end{equation*}
Therefore, with probability at least $1 - C\exp(-C\frac{s\log(es_gb) + s_g\log(d/s_g)}{s}) - Ce^{-cn/s}$, 
\begin{equation}\label{ineq:T}
	\|h_{\widetilde{T}}\|_2 \leq C\left(\frac{\sqrt{s}}{n}\lambda + \frac{1}{\sqrt{s}}\|\beta^*_{T^c}\|_1 + \frac{1}{\sqrt{s_g}}\|\beta^*_{T^c}\|_{1,2}\right).
\end{equation}

By \eqref{eq19}, \eqref{eq20}, \eqref{ineq:h_{T^c}} and the previous inequality, also notice that $e^{-cn/s} \leq e^{-C\frac{s\log(es_gb) + s_g\log(d/s_g)}{s}}$, with probability at least $1 - C\exp(-C\frac{s\log(es_gb) + s_g\log(d/s_g)}{s})$, 
\begin{equation}\label{eq39}
	\begin{split}
	\|h\|_2 \leq& \|h_{\widetilde{T}}\|_2 + \sum_{i \geq 2}\|h_{T_i}\|_2 + \sum_{j \geq 2}\|h_{R_j}\|_2 \leq \|h_{\widetilde{T}}\|_2 + \sqrt{2}s^{-1/2}\|h_{T^c}\|_2 + s_g^{-1/2}\|h_{(G^c)}\|_{1, 2}\\ \leq& C\left(\frac{\sqrt{s}}{n}\lambda + \frac{1}{\sqrt{s}}\|\beta^*_{T^c}\|_1 + \frac{1}{\sqrt{s_g}}\|\beta^*_{T^c}\|_{1,2}\right),
	\end{split}
\end{equation}
i.e., with probability at least $1 - C\exp(-C\frac{s\log(es_gb) + s_g\log(d/s_g)}{s})$, 
\begin{equation*}
\|h\|_2 \leq C\left(\sqrt{\frac{\sigma^2(s_g\log(d/s_g) + s\log(es_gb))}{n}} + \frac{1}{\sqrt{s}}\|\beta^*_{T^c}\|_1 + \frac{1}{\sqrt{s_g}}\|\beta^*_{T^c}\|_{1,2}\right).
\end{equation*}

Moreover, if $\beta^*$ is $(s, s_g)$-sparse, then $\|\beta^*_{T^c}\|_1 = \|\beta^*_{T^c}\|_{1,2} = 0$. Therefore, with probability at least $1 - C\exp(-C\frac{s\log(es_gb) + s_g\log(d/s_g)}{s})$, 
\begin{equation*}
	\|h\|_2^2 \leq \frac{C\sigma^2(s_g\log(d/s_g) + s\log(es_gb))}{n}.
\end{equation*}
\quad $\square$

\subsection{Proof of Theorem \ref{th:lower-noisy}}
First, we consider the case that $d \geq 3s_g$ and $b \geq 3s/s_g$.
Let $\omega^{(1)}, \dots, \omega^{(N)}$ be uniformly randomly vectors from 
\begin{equation*}
	A = \{\omega \in \{0, 1\}^{db}|\sum_j 1_{\{\omega_{(j)} \neq 0\}} = s_g, \|\omega_{(j)}\|_0 = \lfloor s/s_g\rfloor \text{ if } \omega_{(j)} \neq 0\}. 
\end{equation*}
Denote $\Omega^{(i)} = \{j|\omega_{j}^{(i)} \neq 0\}$, $\Omega_{(k)}^{(i)} = \{j|j \in (k), \omega_{j}^{(i)} \neq 0\}$ and $\beta^{(i)} = \tau \omega^{(i)}$,  for all $1 \leq i \leq N, 1 \leq k \leq d$, where $\tau$ is a parameter that will be specified later. Obviously, $\|\beta^{(i)}\|_0 = s_g\lfloor s/s_g\rfloor \leq s$, therefore $\|\beta^{(i)} - \beta^{(j)}\|_2^2 \leq 2s_g\lfloor s/s_g\rfloor\tau^2 \leq 2s\tau^2$. 

Moreover, if $|\Omega^{(i)} \cap \Omega^{(j)}| \geq 8s_g\lfloor s/s_g\rfloor/9$, then we must have 
\begin{equation*}
	\left|\left\{k|\omega_{(k)}^{(i)}, \omega_{(k)}^{(j)} \neq 0, \left|\Omega_{(k)}^{(i)} \cap \Omega_{(k)}^{(j)}\right| \geq 2\lfloor s/s_g\rfloor/3\right\}\right| \geq 2s_g/3,
\end{equation*}
otherwise $|\Omega^{(i)} \cap \Omega^{(j)}| \leq \frac{2s_g}{3}\lfloor s/s_g\rfloor + \frac{s_g}{3}2\lfloor s/s_g\rfloor/3 \leq 8s_g\lfloor s/s_g\rfloor/9$, which is a contradiction.

Therefore,
\begin{equation}\label{eq38}
	\begin{split}
	&\bbP\left(\|\beta^{(i)} - \beta^{(j)}\|_2^2 \leq 2s_g\lfloor s/s_g\rfloor\tau^2/9\right)\\ =& \bbP\left(|\Omega^{(i)} \cap \Omega^{(j)}| \geq 8s_g\lfloor s/s_g\rfloor/9\right)\\ \leq& \bbP\left(\left|\left\{k|\omega_{(k)}^{(i)}, \omega_{(k)}^{(j)} \neq 0, \left|\Omega_{(k)}^{(i)} \cap \Omega_{(k)}^{(j)}\right| \geq 2\lfloor s/s_g\rfloor/3\right\}\right| \geq 2s_g/3\right)\\
	\leq& \frac{\sum_{l = \lceil 2s_g/3\rceil}^{s_g}\binom{s_g}{l}\left[\sum_{t = \lceil 2\lfloor s/s_g\rfloor/3\rceil}^{\lfloor s/s_g\rfloor}\binom{\lfloor s/s_g\rfloor}{t}\binom{b - \lfloor s/s_g\rfloor}{\lfloor s/s_g\rfloor - t}\right]^l\binom{b}{\lfloor s/s_g\rfloor}^{s_g - l}\binom{d - l}{s_g - l}}{\binom{d}{s_g}\binom{b}{\lfloor s/s_g\rfloor}^{s_g}}\\
	=& \sum_{l = \lceil 2s_g/3\rceil}^{s_g}\binom{s_g}{l}\frac{\binom{d - l}{s_g - l}}{\binom{d}{s_g}}\cdot \left[\sum_{t = \lceil 2\lfloor s/s_g\rfloor/3\rceil}^{\lfloor s/s_g\rfloor}\binom{\lfloor s/s_g\rfloor}{t}\frac{\binom{b - \lfloor s/s_g\rfloor}{\lfloor s/s_g\rfloor - t}}{\binom{b}{\lfloor s/s_g\rfloor}}\right]^l.
	\end{split}
\end{equation}
Note that 
\begin{equation*}
    \frac{\binom{d - l}{s_g - l}}{\binom{d}{s_g}} = \frac{\frac{(d - l)\cdots (d - s_g + 1)}{(s_g - l)!}}{\frac{d(d - 1)\cdots (d - s_g + 1)}{s_g!}} = \frac{s_g(s_g - 1)\cdots (s_g - l + 1)}{d(d - 1)\cdots (d - l + 1)}\leq \left(\frac{s_g}{d}\right)^l,
\end{equation*}
The inequality holds since $\frac{s_g - i}{d - i} \leq \frac{s_g}{d}$ for all $1 \leq i \leq s_g$.\\
Similarly, for $1 \leq t \leq \lfloor s/s_g\rfloor$,
\begin{equation*}
	\frac{\binom{b - \lfloor s/s_g\rfloor}{\lfloor s/s_g\rfloor - t}}{\binom{b}{\lfloor s/s_g\rfloor}} \leq \frac{\binom{b - t}{\lfloor s/s_g\rfloor - t}}{\binom{b}{\lfloor s/s_g\rfloor}} \leq \left(\frac{\lfloor s/s_g\rfloor}{b}\right)^t.
\end{equation*}
Combine \eqref{eq38} and the previous two inequalities together, we have
\begin{equation}\label{ineq: random}
	\begin{split}
	    &\bbP\left(\|\beta^{(i)} - \beta^{(j)}\|_2^2 \leq 2s_g\lfloor s/s_g\rfloor\tau^2/9\right)\\ \leq& \sum_{l = \lceil 2s_g/3\rceil}^{s_g}\binom{s_g}{l}\left(\frac{s_g}{d}\right)^l\cdot \left[\sum_{t = \lceil 2\lfloor s/s_g\rfloor/3\rceil}^{\lfloor s/s_g\rfloor}\binom{\lfloor s/s_g\rfloor}{t}\left(\frac{\lfloor s/s_g\rfloor}{b}\right)^t\right]^l\\
	    \leq& \sum_{l = \lceil 2s_g/3\rceil}^{s_g}\binom{s_g}{l}\left(\frac{s_g}{d}\right)^l\cdot \left[\sum_{t = \lceil 2\lfloor s/s_g\rfloor/3\rceil}^{\lfloor s/s_g\rfloor}\binom{\lfloor s/s_g\rfloor}{t}\left(\frac{\lfloor s/s_g\rfloor}{b}\right)^{2\lfloor s/s_g\rfloor/3}\right]^l\\
	    \leq& \sum_{l = \lceil 2s_g/3\rceil}^{s_g}\binom{s_g}{l}\left(\frac{s_g}{d}\right)^l\cdot \left[2^{\lfloor s/s_g\rfloor}\left(\frac{\lfloor s/s_g\rfloor}{b}\right)^{2\lfloor s/s_g\rfloor/3}\right]^l\\
	    \leq& \sum_{l = \lceil 2s_g/3\rceil}^{s_g}\binom{s_g}{l}\left(\frac{s_g}{d}\right)^{2s_g/3}\cdot \left[\left(\frac{2\sqrt{2}\lfloor s/s_g\rfloor}{b}\right)^{2\lfloor s/s_g\rfloor/3}\right]^{2s_g/3}\\
	    \leq& \left(\frac{2\sqrt{2}s_g}{d}\right)^{2s_g/3}\cdot \left(\frac{2\sqrt{2}\lfloor s/s_g\rfloor}{b}\right)^{2s/9}.
	\end{split}
\end{equation}
Set $N = \left\lfloor\left(\frac{d}{2\sqrt{2}s_g}\right)^{s_g/3}\left(\frac{b}{2\sqrt{2}\lfloor s/s_g \rfloor}\right)^{s/9}\right\rfloor$, then
\begin{equation*}
    \begin{split}
    &\bbP\left(\forall 1 \leq i \neq j \leq N, \|\beta^{(i)} - \beta^{(j)}\|_2^2 > 2s_g\lfloor s/s_g\rfloor\tau^2/9\right)\\	\geq& 1 - \frac{N(N - 1)}{2}\left(\frac{2\sqrt{2}s_g}{d}\right)^{2s_g/3}\cdot \left(\frac{2\sqrt{2}\lfloor s/s_g\rfloor}{b}\right)^{2s/9}\\ >& 0.
    \end{split}
\end{equation*}

i.e., the probability that $\beta^{(1)}, \dots, \beta^{(N)}; \Omega^{(1)}, \cdots, \Omega^{(N)}$ satisfy
\begin{equation}\label{eq24}
	\frac{s}{9}\tau^2 < 2s_g\lfloor s/s_g\rfloor\tau^2/9 < \min_{i \neq j}\|\beta^{(i)} - \beta^{(j)}\|_2^2 \leq 2s\tau^2,
\end{equation}
\begin{equation}\label{eq25}
	|\Omega^{(i)} \cap \Omega^{(j)}| < 8s_g\lfloor s/s_g\rfloor/9, \quad \forall 1 \leq i < j \leq N
\end{equation}
is positive. For convenience, we fix $\beta^{(1)}, \dots, \beta^{(N)}$ to be the vectors satisfying \eqref{eq24}.

Denote $y^{(i)} = X\beta^{(i)} + \varepsilon$ for all $1 \leq i \leq N$. We consider the Kullback-Leibler divergence between different distribution pairs:
\begin{equation*}
	D_{KL}\left((y^{(i)}, X), (y^{(j)}, X)\right) = \mathbb{E}_{(y^{(j)}, X)}\left[\log\left(\frac{p(y^{(i)}, X)}{p(y^{(j)}, X)}\right)\right],
\end{equation*}
where $p(y^{(i)}, X)$ is the probability density of $(y^{(i)}, X)$. Conditioning on $X$, we have
\begin{equation*}
	\mathbb{E}_{(y^{(j)}, X)}\left[\log\left(\frac{p(y^{(i)}, X)}{p(y^{(j)}, X)}\right)|X\right] = \frac{\|X(\beta^{(i)} - \beta^{(j)})\|_2^2}{2\sigma^2}.
\end{equation*}
Thus for $1 \leq i \neq j \leq N$,
\begin{equation}\label{ineq:KL}
	\begin{split}
	D_{KL}\left((y^{(i)}, X), (y^{(j)}, X)\right) =& \mathbb{E}_X \frac{\|X(\beta^{(i)} - \beta^{(j)})\|_2^2}{2\sigma^2} = \frac{n(\beta^{(i)} - \beta^{(j)})^\top\Sigma (\beta^{(i)} - \beta^{(j)})}{2\sigma^2}\\ \leq& \frac{3n\|\beta^{(i)} - \beta^{(j)}\|_2^2}{4\sigma^2} \leq \frac{3ns\tau^2}{2\sigma^2}.
	\end{split}
\end{equation}
In the first inequality, we used $\sigma_{\max}(\Sigma) \leq \frac{3}{2}$.\\
By generalized Fano's Lemma,
\begin{equation*}
	\inf_{\hat{\beta}}\sup_{\beta \in \mathbb{F}_{s, s_g}} \mathbb{E}\|\hat{\beta} - \beta\|_2 \geq \frac{\sqrt{s\tau^2/9}}{2}\left(1 - \frac{\frac{3ns\tau^2}{2\sigma^2} + \log 2}{\log N}\right).
\end{equation*}
Since $\log N \asymp s_g\log(\frac{d}{s_g}) + s\log\left(\frac{es_gb}{s}\right)$, by setting $\tau = c\sqrt{\frac{\sigma^2\left(s_g\log(\frac{d}{s_g}) + s\log\left(\frac{es_gb}{s}\right)\right)}{ns}}$, we have
\begin{equation*}
	\inf_{\hat{\beta}}\sup_{\beta \in \mathbb{F}_{s, s_g}} \mathbb{E}\|\hat{\beta} - \beta\|_2^2 \geq \left(\inf_{\hat{\beta}}\sup_{\beta \in \mathbb{F}_{s, s_g}} \mathbb{E}\|\hat{\beta} - \beta\|_2\right)^2 \geq c\frac{\sigma^2\left(s_g\log(d/s_g) + s\log(es_gb/s)\right)}{n}. 
\end{equation*}

If $d < 3s_g$ or $b < 3s/s_g$, let $s_g'=[s_g/3] \vee 1 \geq s_g/5, s' = [s/15] \vee s_g'$, then $d \geq 3s_g'$ and $b \geq 3s'/s_g'$. Similarly to \eqref{ineq1}, we have
\begin{equation*}
    \begin{split}
    \inf_{\hat{\beta}}\sup_{\beta \in \mathbb{F}_{s, s_g}} \mathbb{E}\|\hat{\beta} - \beta\|_2^2 \geq& \inf_{\hat{\beta}}\sup_{\beta \in \mathbb{F}_{s', s_g'}} \mathbb{E}\|\hat{\beta} - \beta\|_2^2 \geq c\frac{\sigma^2\left(s_g'\log(d/s_g') + s'\log(es_g'b/s')\right)}{n}\\ \geq& c'\frac{\sigma^2\left(s_g\log(d/s_g) + s\log(es_gb/s)\right)}{n}.
    \end{split}
\end{equation*}
\quad $\square$
\subsection{Proof of Theorem \ref{th:unbiased}}

The proof of Theorem \ref{th:unbiased} relies on the following key lemma, which shows that $\Sigma^{-1}$ is in the feasible set of the optimization problem \eqref{opt:M} with high probability by choosing appropriate $\alpha$ and $\gamma$.
\begin{Lemma}\label{lm:optimization}
	By setting $\alpha = C\sqrt{\frac{s\log(es_gb) + s_g\log(d/s_g)}{sn}}, \gamma = \sqrt{\frac{s}{s_g}}\alpha$ in \eqref{opt:M}, we have 
	\begin{equation*}
	\bbP\left(\max_{1 \leq i \leq p}\|H_{\alpha}(e_i - \frac{1}{n}X^\top X\Sigma^{-1}e_i)\|_{\infty, 2} \leq \gamma\right) \geq 1 - 4\exp\left(- C\frac{s\log(es_gb) + s_g\log(d/s_g)}{s_g}\right).
	\end{equation*}
\end{Lemma}

Note that $Y = X\beta^* + \varepsilon$, we have
\begin{equation*}
\sqrt{n}(\hat{\beta}^u - \beta^*) = \sqrt{n}\left(\hat{\beta} - \beta^* + \frac{1}{n}\hat{M}X^\top\left(Y - X\hat{\beta}\right)\right) = \sqrt{n}\left(I - \frac{1}{n}\hat{M}X^\top X\right)(\hat{\beta} - \beta^*) + \frac{1}{\sqrt{n}}\hat{M}X^\top\varepsilon.
\end{equation*}
Since $\varepsilon_i \stackrel{i.i.d.}{\sim} N(0, \sigma^2)$, we know that 
\begin{equation*}
	\frac{1}{\sqrt{n}}\hat{M}X^\top\varepsilon|X \sim N\left(0, \hat{M}\hat{\Sigma}\hat{M}^\top\right).
\end{equation*}
Denote $h = \hat{\beta} - \beta^*$. Since $\beta^*$ is $(s, s_g)$-sparse, by \eqref{ineq:h_{T^c}}, \eqref{eq39} and Cauchy-Schwarz inequality, with probability at least $1 - C\exp(-C\frac{s\log(es_gb) + s_g\log(d/s_g)}{s})$, 
\begin{equation*}
\|h\|_1 \leq \|h_{T}\|_1 + \|h_{T^c}\|_1 \leq \sqrt{s}\|h_{T}\|_2 + \|h_{T^c}\|_1 \leq \sqrt{s}\|h\|_2 + \|h_{T^c}\|_1 \leq C\frac{s}{n}\lambda.
\end{equation*}
\begin{equation*}
	\|h\|_{1,2} \leq \|h_{(G)}\|_{1,2} + \|h_{(G^c)}\|_{1,2} \leq \sqrt{s_g}\|h_{(G)}\|_2 + \|h_{(G^c)}\|_{1,2} \leq \sqrt{s_g}\|h\|_2 + \|h_{(G^c)}\|_{1,2} \leq C\frac{\sqrt{s\cdot s_g}}{n}\lambda.
\end{equation*}
In addition, Lemma \ref{lm:optimization} shows that $\Sigma^{-1}$ is in the feasible set of \eqref{opt:M} with probability at least $1 - C\exp(-C\frac{s\log(es_gb) + s_g\log(d/s_g)}{s})$. By the definition of $\hat{M}$,
\begin{equation}\label{eq30}
\max_{i} \|H_{\alpha}(e_i - \hat{\Sigma}\hat{M}^\top e_i)\|_{\infty, 2} = \max_{i} \|H_{\alpha}(e_i - \hat{\Sigma}\hat{m}_i)\|_{\infty, 2} \leq \gamma.
\end{equation}
Combining these facts, by Lemma \ref{lm:soft-thresholding} Part 2, we must have
\begin{equation*}
    \begin{split}
    \left\|(I - \frac{1}{n}\hat{M}XX^\top)(\hat{\beta} - \beta^*)\right\|_\infty =& \max_{i}\left|\langle e_i - \hat{\Sigma}\hat{M}^\top e_i, h\rangle\right|\\ \leq& \alpha\|h\|_1 + \gamma\|h\|_{1,2}\\ \leq& C\frac{s}{n}\alpha\lambda + C\frac{\sqrt{s\cdot s_g}}{n}\gamma\lambda\\ =& \frac{C(s\log(es_gb) + s_g\log(d/s_g))}{n}\sigma
    \end{split}
\end{equation*}
with probability at least $1 - C\exp(-C\frac{s\log(es_gb) + s_g\log(d/s_g)}{s})$.
This has finished the proof of \eqref{de-biased1}.
 
 Next, we consider $\hat{m}_i^\top\hat{\Sigma}\hat{m}_i$. By \eqref{eq30} and Lemma \ref{lm:soft-thresholding} Part 2, we have
 \begin{equation*}
     1 - \langle e_i, \hat{\Sigma}\hat{m}_i\rangle = \langle e_i, e_i - \hat{\Sigma}\hat{m}_i\rangle \leq \alpha\|e_i\|_1 + \gamma\|e_i\|_{1,2} = \alpha + \gamma.
 \end{equation*}
Therefore, for any $c \geq 0$,
\begin{equation*}
	\hat{m}_i^\top\hat{\Sigma}\hat{m}_i \geq \hat{m}_i^\top\hat{\Sigma}\hat{m}_i + c(1 - \alpha - \gamma) - c\langle e_i, \hat{\Sigma}\hat{m}_i\rangle \geq \min_{m}\left\{m^\top\hat{\Sigma}m + c(1 - \alpha - \gamma) - c\langle e_i, \hat{\Sigma}m\rangle\right\}.
\end{equation*}
Since $m = ce_i/2$ achieves the minimum of the right hand side, we have
\begin{equation*}
	\hat{m}_i^\top\hat{\Sigma}\hat{m}_i \geq c(1 - \alpha - \gamma) - \frac{c^2}{4}\hat{\Sigma}_{i,i}.
\end{equation*}
If $\hat{\Sigma}_{ii} > 0$ for all $1 \leq i \leq p$, by setting $c = 2(1 - \alpha - \gamma)/\hat{\Sigma}_{i,i}$, we have
\begin{equation}\label{eq31}
	\hat{m}_i^\top\hat{\Sigma}\hat{m}_i \geq \frac{(1 - \alpha - \gamma)^2}{\hat{\Sigma}_{i,i}}, \quad \forall 1 \leq i \leq p.
\end{equation}
Moreover, by Lemma \ref{lm:sub-Gaussian-concentration} Part 2 with $u = v = e_i$, we have
\begin{equation*}
	\bbP\left(\left|\hat{\Sigma}_{i,i} - \Sigma_{i, i}\right| \geq \frac{c_{\min}}{2}\right) \leq 2\exp\left(-cn\right).
\end{equation*}
By the union bound,
\begin{equation*}
	\begin{split}
	    &\bbP\left(\exists 1 \leq i \leq p, \left|\hat{\Sigma}_{i,i} - \Sigma_{i, i}\right| \geq \frac{c_{\min}}{2}\right)\\
	    \leq& \sum_{i = 1}^{p}\bbP\left(\left|\hat{\Sigma}_{i,i} - \Sigma_{i, i}\right| \geq \frac{c_{\min}}{2}\right)\\
	    \leq& db\cdot 2\exp\left(-cn\right)\\ \leq& 2\exp\left(-cn\right).
	\end{split}
\end{equation*}
Therefore, with probability at least $1 - 2\exp\left(-cn\right)$, 
\begin{equation*}
	\frac{c_{\min}}{2} \leq \hat{\Sigma}_{i,i} \leq C_{\max} + \frac{c_{\min}}{2}, \quad \forall 1 \leq i \leq p.
\end{equation*}
\eqref{eq31} and the previous inequality together imply that with probability at least $1 - 2\exp\left(-cn\right)$, 
\begin{equation*}
    \hat{m}_i^\top\hat{\Sigma}\hat{m}_i \geq \frac{1}{2C_{\max}}, \quad \forall 1 \leq i \leq p.
\end{equation*}
\eqref{de-biased1} and the previous inequality together imply \eqref{de-biased2}.
\quad $\square$

\bibliographystyle{ieeetr}
\bibliography{reference}

\appendix

\section{Technical Lemmas}\label{sec:technical lemma}
We collect all additional technical lemmas and their proofs in this section.

\begin{Lemma}[Bernstein-type Inequality for Soft-thresholded Sub-Gaussian Vectors]\label{lm:bernstein-sub-gaussian}
	~\\ Suppose the rows of $X\in \mathbb{R}^{n \times p}$ are independent sub-Gaussian vectors satisfying Assumption \ref{as:design}. $w\in \mathbb{R}^{n}$ is a fixed vector, $\Omega$ is a subset of $\{1,\ldots, p\}$ with $|\Omega| = r$. Then
	\begin{equation}\label{ineq:bernstein-sub-gaussian}
	\bbP\left(\left\|\sum_{k=1}^n w_k X_{k, \Omega}\right\|_2 \geq \sqrt{C_{\max}}\kappa\|w\|_2\cdot\left(\sqrt{r} + \sqrt{2t}\right)\right) \leq \exp(-t).
	\end{equation}
	For any fixed vector $w\in \mathbb{R}^n$ and fixed index subset $\Omega\subseteq \{1,\ldots, p\}$ with $|\Omega| = r$, 
	\begin{equation}\label{ineq:soft_thresholding_concentration}
	\begin{split}
	\bbP\left(\left\|H_{(\delta\|w\|_2)}(w^\top X_{\Omega})\right\|_2 \geq t\|w\|_2\right) \leq & \binom{r}{\lfloor(t/\delta)^2\rfloor\wedge r} \cdot \exp\left(-(t/(\kappa\sqrt{C_{\max}}) - (t/\delta)\wedge \sqrt{r})_+^2/2\right)\\
	& + \binom{r}{\lceil(t/\delta)^2\rceil}\cdot \exp\left(-(t/(\kappa\sqrt{C_{\max}}) - \sqrt{\lceil (t/\delta)^2 \rceil})_+^2/2\right).
	\end{split}
	\end{equation}
	In particular, for any $b \geq r$, if $\bar{\lambda} = C\|w\|_2\sqrt{\frac{s\log(es_gb) + s_g\log(d/s_g)}{s}}$, $\bar{\lambda}_g = \sqrt{s/s_g}\bar{\lambda}$, we have
	\begin{equation}\label{ineq:particular_lemma3}
	\begin{split}
	\bbP\left(\left\| H_{\bar{\lambda}}(w^\top X_{\Omega})\right\|_2\geq \bar{\lambda}_g \right) \leq \exp\left(-C\frac{s\log(es_gb) + s_g\log(d/s_g)}{s_g}\right).
	\end{split}
	\end{equation}
\end{Lemma}

{\bf\noindent Proof of Lemma \ref{lm:bernstein-sub-gaussian}.} We only need to focus on the case where $\|w\|_2=1$. Let $W_{\Omega} = X_{\Omega}\Sigma_{\Omega, \Omega}^{-1/2}$, immediately we know that $W_{1, \Omega},\ldots, W_{n, \Omega}$ are isotropic sub-Gaussian distributed. Then for any fixed $w$, $w^\top W_{\Omega}$ is also an isotropic sub-Gaussian vector such that for any $\alpha\in \mathbb{R}^{r}$,
\begin{equation*}
    \begin{split}
    \mathbb{E}\exp\left(w^\top W_{\Omega} \alpha\right) =& \mathbb{E}\exp\left(w^\top X_{\Omega}\Sigma_{\Omega, \Omega}^{-1/2}\alpha\right) = \mathbb{E}\exp\left(w^\top X\Sigma^{-1/2}(\Sigma^{1/2})_{\cdot,\Omega}\Sigma_{\Omega, \Omega}^{-1/2}\alpha\right)\\ \leq& \exp\left(\kappa^2\|(\Sigma^{1/2})_{\cdot,\Omega}\Sigma_{\Omega, \Omega}^{-1/2}\alpha\|_2^2/2\right) = \exp\left(\kappa^2\|\alpha\|_2^2/2\right).
    \end{split}
\end{equation*}
The last equation holds since $(\Sigma^{1/2})_{\Omega,\cdot}(\Sigma^{1/2})_{\cdot,\Omega} = (\Sigma^{1/2}\Sigma^{1/2})_{\Omega,\Omega} = \Sigma_{\Omega,\Omega}$.\\
By the tail inequality of sub-Gaussian quadratic form (\cite[Theorem 2.1]{hsu2012tail}), 
\begin{equation*}
\bbP\left(\left\|w^\top W_{\Omega}\right\|_2^2 \geq \kappa^2\left(r + 2\sqrt{rt} + 2t\right)\right) \leq \exp(-t).
\end{equation*}
By taking square-root of the previous inequality, we have
\begin{equation*}
	\bbP\left(\left\|w^\top W_{\Omega}\right\|_2 \geq \kappa\|w\|_2\cdot\left(\sqrt{r} + \sqrt{2t}\right)\right) \leq \exp(-t).
\end{equation*}
Also note that 
\begin{equation*}
	\left\|w^\top X_{\Omega}\right\|_2 = \left\|w^\top W_{\Omega}\Sigma_{\Omega, \Omega}^{1/2}\right\|_2 \leq \left\|\Sigma_{\Omega, \Omega}^{1/2}\right\|\left\|w^\top W_{\Omega}\right\|_2 \leq \|\Sigma\|^{1/2}\left\|w^\top W_{\Omega}\right\|_2 \leq \sqrt{C_{\max}}\left\|w^\top W_{\Omega}\right\|_2,
\end{equation*}
we obtain \eqref{ineq:bernstein-sub-gaussian}. 

For the second part of proof, note that
\begin{equation*}
\begin{split}
& \bbP\left(\left\|H_\delta(w^\top X_{\Omega})\right\| \geq t\right)\\
\leq & \bbP\left(\exists \Lambda\subseteq \Omega, \text{ such that all entries of $|w^\top X_{\Lambda}| \geq \delta$ and $\|w^\top X_{\Lambda}\|_2\geq t$}\right)\\
\leq & \bbP\left(\exists \Lambda\subseteq \Omega, \sqrt{|\Lambda|}\delta \leq t, \|w^\top X_{\Lambda}\|_2\geq t\right)\\
& + \bbP\left(\exists \Lambda\subseteq \Omega, \sqrt{|\Lambda|}\delta>t, \text{ all entries of $|w^\top X_{\Lambda}|\geq \delta$}\right)\\
\leq & \sum_{\substack{\Lambda\subseteq \Omega\\|\Lambda| = \lfloor(t/\delta)^2\rfloor\wedge r}} \bbP\left(\|w^\top X_{\Lambda}\|_2\geq t\right) + \sum_{\substack{\Lambda\subseteq \Omega\\|\Lambda| = \lceil (t/\delta)^2 \rceil}} \bbP\left(\|w^\top X_\Lambda\|_2\geq t\right).
\end{split}
\end{equation*}
By the first part of this lemma,
\begin{equation*}
\bbP\left(\|w^\top X_{\Lambda}\|_2\geq \sqrt{C_{\max}}\kappa\|w\|_2 t\right) \leq \exp\left(- \left(t - \sqrt{|\Lambda|}\right)_+^2/2\right).
\end{equation*}
Plug in this to the previous inequality, one has
\begin{equation*}
\begin{split}
\bbP\left(\left\|H_\delta(w^\top X_{\Omega})\right\| \geq t\right) \leq & \binom{r}{\lfloor(t/\delta)^2\rfloor\wedge r} \cdot \exp\left(-(t/(\kappa\sqrt{C_{\max}}) - (t/\delta)\wedge \sqrt{r})_+^2/2\right) \\
& + \binom{r}{\lceil(t/\delta)^2\rceil}\cdot \exp\left(-(t/(\kappa\sqrt{C_{\max}}) - \sqrt{\lceil (t/\delta)^2 \rceil})_+^2/2\right).
\end{split}
\end{equation*}

Specifically, if $\delta = C\sqrt{\frac{s\log(es_gb) + s_g\log(d/s_g)}{s}}$, $t = \sqrt{s/s_g}\delta$, 
\begin{equation*}
	t/(\kappa\sqrt{C_{\max}}) - \sqrt{\lceil (t/\delta)^2 \rceil} \geq C\sqrt{\frac{s\log(es_gb) + s_g\log(d/s_g)}{s_g}} - \sqrt{2\frac{s}{s_g}} \geq C\sqrt{\frac{s\log(es_gb) + s_g\log(d/s_g)}{s_g}}.
\end{equation*}
Therefore, \eqref{ineq:soft_thresholding_concentration} shows that
\begin{equation*}
	\begin{split}
	\bbP\left(\left\| H_{\bar{\lambda}}(w^\top X_{\Omega})\right\|_2\geq \bar{\lambda}_g \right) \leq& r^{\lfloor(t/\delta)^2\rfloor}\exp\left(-C\frac{s\log(es_gb) + s_g\log(d/s_g)}{s_g}\right)\\ &+ r^{\lceil(t/\delta)^2\rceil}\exp\left(-C\frac{s\log(es_gb) + s_g\log(d/s_g)}{s_g}\right)\\ \leq& 2r^{2s/s_g}\exp\left(-C\frac{s\log(es_gb) + s_g\log(d/s_g)}{s_g}\right)\\ \leq& \exp\left(\log 2 + \frac{2s\log(eb)}{s_g} - C\frac{s\log(es_gb) + s_g\log(d/s_g)}{s_g}\right) \\ \leq& \exp\left(-C\frac{s\log(es_gb) + s_g\log(d/s_g)}{s_g}\right).
	\end{split}
\end{equation*}
\quad $\square$

\begin{Lemma}[sub-Gaussian quadratic form concentrations]\label{lm:sub-Gaussian-concentration}
	Suppose $Z \in \mathbb{R}^p$ is a sub-Gaussian vector satisfying Assumption \ref{as:design}. 
	\begin{enumerate}[leftmargin=*]
	\item For any fixed $u, v\in \mathbb{R}^{p}, u, v \neq 0$, $u^\top Z Z^\top v$ is sub-exponential such that for every $t > 0$,
	\begin{equation}\label{ineq:quadratic-concentration-1}
	\mathbb{P}\left(\left|u^\top Z Z^\top v - \mathbb{E}u^\top Z Z^\top v \right| \geq t\|u\|_2\|v\|_2\right) \leq C\exp(-ct/\kappa^2).
	\end{equation}

	\item In addition, suppose $X = [X_1^\top,\ldots, X_n^\top]^\top \in \mathbb{R}^{n\times p}$ is a random matrix with independent random sub-Gaussian rows satisfying Assumption \ref{as:design},
	\begin{equation}\label{ineq:quadratic-concentration-2}
	\bbP\left(\left|\frac{1}{n}\sum_{k=1}^n u^\top X_kX_k^\top v - u^\top \Sigma v\right| \geq t\|u\|_2\|v\|_2\right)\leq 2\exp\left(-cn\min\left\{\frac{t^2}{\kappa^4}, \frac{t}{\kappa^2}\right\}\right).
	\end{equation}
	
	\item More generally, for any fixed matrix $U\in \mathbb{R}^{p\times r}$, the following concentration inequality in spectral norm holds,
	\begin{equation}\label{ineq:quadratic-concentration-3}
	\bbP\left(\left\|\frac{1}{n}\sum_{k=1}^n U^\top X_kX_k^\top v - U^\top \Sigma v\right\|_2 \geq t\|U\|\|v\|_2\right) \leq 2\exp\left(Cr -cn\min\left\{\frac{t^2}{\kappa^4}, \frac{t}{\kappa^2}\right\}\right).
	\end{equation}
	\end{enumerate}
\end{Lemma}

{\bf\noindent Proof of Lemma \ref{lm:sub-Gaussian-concentration}.} Since we can rescale $u$ and $v$ without essentially changing the problem, without loss of generality we assume $\|u\|_2 = \|v\|_2=1$. Let $A = uv^\top$, then $u^\top ZZ^\top v = Z^\top uv^\top Z = Z^\top AZ$. By Assumption \ref{as:design}, $\mathbb{E} Z = 0$ and $\|\langle Z, e_i\rangle\|_{\psi_2} \leq C\kappa$. By Hanson-Wright inequality (\cite[Theorem 1.1]{rudelson2013hanson}),
\begin{equation*}
	\begin{split}
	\mathbb{P}\left(\left|u^\top Z Z^\top v - \mathbb{E}u^\top Z Z^\top v \right| \geq t\right) =& \bbP\left(|Z^\top AZ - \mathbb{E}Z^\top AZ| \geq t\right)\\ \leq& 2\exp\left[-c\min\left(\frac{t^2}{\kappa^4\|A\|_{HS}^2}, \frac{t}{\kappa^2\|A\|}\right)\right]\\
	\leq& 2\exp\left[-c\min\left(\frac{t^2}{\kappa^4}, \frac{t}{\kappa^2}\right)\right],
	\end{split}
\end{equation*}
where 
\begin{equation*}
	\|A\|_{HS} = \left(\sum_{i,j}|a_{i,j}|^2\right)^{1/2} = \left(\sum_{i, j}|u_iv_j|^2\right)^{1/2} = \|u\|_2\|v\|_2 = 1,
\end{equation*}
\begin{equation*}
	\|A\| = \max_{\|x\|_2 \leq 1}\|Ax\|_2 = \max_{\|x\|_2 \leq 1}\|uv^\top x\|_2 = \|u\|_2\max_{\|x\|_2 \leq 1}|v^\top x| = \|u\|_2\|v\|_2 = 1.
\end{equation*}
Therefore, for every $t \geq \kappa^2$, 
\begin{equation*}
	\mathbb{P}\left(\left|u^\top Z Z^\top v - \mathbb{E}u^\top Z Z^\top v \right| \geq t\right) \leq 2\exp\left(-ct/\kappa^2\right).
\end{equation*}
Thus, there exists a constant $c < \log 2$, for every $t \geq 0$,
\begin{equation*}
    \mathbb{P}\left(\left|u^\top Z Z^\top v - \mathbb{E}u^\top Z Z^\top v \right| \geq t\right) \leq 2\exp\left(-ct/\kappa^2\right).
\end{equation*}
Notice that $\mathbb{E}u^\top X_kX_k^\top v = u^\top\Sigma v$ for all $1 \leq k \leq n$, by Bernstein-type concentration inequality (c.f., \cite[Proposition 5.16]{vershynin2010introduction}),
\begin{equation*}
\begin{split}
\bbP\left(\left|\frac{1}{n}\sum_{k=1}^n u^\top X_kX_k^\top v - u^\top\Sigma v\right| \geq t\right)\leq 2\exp\left(-cn\min\left\{\frac{t^2}{\kappa^4}, \frac{t}{\kappa^2}\right\}\right).
\end{split}
\end{equation*}
This has finished the proof of \eqref{ineq:quadratic-concentration-2}. 

Finally, we consider \eqref{ineq:quadratic-concentration-3}, which can be done by an $\varepsilon$-net argument and the result in \eqref{ineq:quadratic-concentration-2}. 
For any $w \in \bbR^r, \|w\|_2 =1$, set $u = Uw$ in \eqref{ineq:quadratic-concentration-2}, we have
\begin{equation*}
	\bbP\left(\left|\frac{1}{n}\sum_{k=1}^n w^\top U^\top X_kX_k^\top v - w^\top U^\top \Sigma v\right| \geq \frac{t}{2}\|Uw\|_2\|v\|_2\right) \leq 2\exp\left(-cn\min\left\{\frac{t^2}{\kappa^4}, \frac{t}{\kappa^2}\right\}\right).
\end{equation*}
By \cite[Lemma 5.3]{vershynin2010introduction}, we can find a $\frac{1}{2}$-net $\mathcal{N}_{\frac{1}{2}}$ of $S^{r - 1} = \{x|x \in \bbR^r, \|x\|_2 = 1\}$ with $|\mathcal{N}_{\frac{1}{2}}| \leq 5^r$. 
By the union bound,
\begin{equation}\label{eq4}
	\begin{split}
	&\bbP\left(\forall w \in \mathcal{N}_{\frac{1}{2}}, \left|\frac{1}{n}\sum_{k=1}^n w^\top U^\top X_kX_k^\top v - w^\top U^\top \Sigma v\right| \geq \frac{t}{2}\|Uw\|_2\|v\|_2\right)\\ \leq& 5^r\cdot 2\exp\left(-cn\min\left\{\frac{t^2}{\kappa^4}, \frac{t}{\kappa^2}\right\}\right).
	\end{split}
\end{equation}
For any $g \in \bbR^r, g \neq 0$, set $x = \frac{g}{\|g\|_2} \in \argmax_{w \in \bbR^r, \|w\|_2 = 1}|w^\top g|$, we can find $y \in \mathcal{N}_{\frac{1}{2}}$ such that $\|x - y\|_2 \leq \frac{1}{2}$. By triangle inequality,
\begin{equation*}
	\|g\|_2 - |y^\top g| = |x^\top g| - |y^\top g| \leq |x^\top g - y^\top g| \leq \|x - y\|_2\|g\|_2 \leq \frac{1}{2}\|g\|_2.
\end{equation*}
Therefore, 
\begin{equation*}
    \sup_{w \in \bbR^r, \|w\|_2 = 1}\left|\frac{1}{n}\sum_{k=1}^n w^\top U^\top X_kX_k^\top v - w^\top U^\top \Sigma v\right| \leq 2\sup_{w \in \mathcal{N}_{\frac{1}{2}}}\left|\frac{1}{n}\sum_{k=1}^n w^\top U^\top X_kX_k^\top v - w^\top U^\top \Sigma v\right|.	
\end{equation*}
The \eqref{eq4} and the previous inequality together, also notice that $\|U\| = \sup_{w \in \bbR^r, \|w\|_2 = 1}\|Uw\|_2$, we have
\begin{equation}\label{eq3}
	\begin{split}
	&\bbP\left(\sup_{w \in \bbR^r, \|w\|_2 = 1}\left|\frac{1}{n}\sum_{k=1}^n w^\top U^\top X_kX_k^\top v - w^\top U^\top \Sigma v\right| \geq t\|U\|\|v\|_2\right)\\
	\leq& \bbP\left(\sup_{w \in \mathcal{N}_{\frac{1}{2}}}\left|\frac{1}{n}\sum_{k=1}^n w^\top U^\top X_kX_k^\top v - w^\top U^\top \Sigma v\right| \geq \frac{t}{2}\|U\|\|v\|_2\right)\\
	\leq& \bbP\left(\forall w \in \mathcal{N}_{\frac{1}{2}}, \left|\frac{1}{n}\sum_{k=1}^n w^\top U^\top X_kX_k^\top v - w^\top U^\top \Sigma v\right| \geq \frac{t}{2}\|Uw\|_2\|v\|_2\right)\\
	\leq& 5^r\cdot 2\exp\left(-cn\min\left\{\frac{t^2}{\kappa^4}, \frac{t}{\kappa^2}\right\}\right).
	\end{split}
\end{equation}
Finally, note that 
\begin{equation*}
	\left\|\frac{1}{n}\sum_{k=1}^n U^\top X_kX_k^\top v - U^\top \Sigma v\right\|_2 = \sup_{w \in \bbR^r, \|w\|_2 = 1}\left|\frac{1}{n}\sum_{k=1}^n w^\top U^\top X_kX_k^\top v - w^\top U^\top \Sigma v\right|,
\end{equation*}
we have proved \eqref{ineq:quadratic-concentration-3}.
 \quad $\square$

We collect the random matrix properties of $X$ in the following lemma. These properties will be extensively used in the main content of the paper.
\begin{Lemma}\label{lm:X-random-matrix-properties}
	Suppose $X = [X_1^\top,\ldots, X_n^\top]^\top \in \mathbb{R}^{n\times p}$ is a random matrix with independent random sub-Gaussian rows satisfying Assumption \ref{as:design}.	
	\begin{enumerate}
	\item Suppose $T \subseteq\{1,\ldots, p\}$ is with cardinality $s$. Then,
	\begin{equation}\label{ineq:X_TX_T-tail-bound}
	\bbP\left(\left\|\frac{1}{n}X_T^\top X_T\Sigma_{T, T}^{-1} - I_s\right\| \geq t\right) \leq 2\exp\left(Cs - cn\min\left\{\frac{t^2}{\kappa^4}, \frac{t}{\kappa^2}\right\}\right);
	\end{equation}
	\item For any fixed vector $\alpha\in \mathbb{R}^s$, $\delta>0$, and fixed index subset $\Omega \subseteq T^c$ satisfying $|\Omega| = r$, $t \geq \delta \geq C(\max_{i \in T^c}\|\Sigma_{i, T}\Sigma_{T, T}^{-1}\|_2)\|\alpha\|_2$, 
	\begin{equation}\label{ineq:H_lambda-tail-bound}
	\begin{split}
	& \bbP\left(\left\|H_\delta(\alpha^\top X_T^\top X_\Omega/n)\right\|_2 \geq t\right)\\
	\leq & \binom{r}{\lfloor(t/\delta)^2\rfloor\wedge r}\exp\left(C\lfloor (t/\delta)^2 \rfloor\wedge r - cn\min\left\{\frac{t^2}{\kappa^4\|\alpha\|_2^2}, \frac{t}{\kappa^2\|\alpha\|_2}\right\}\right)\\
	&  + \binom{r}{\lceil(t/\delta)^2\rceil}_+\exp\left(C\lceil (t/\delta)^2 \rceil - cn\min\left\{\frac{t^2}{\kappa^4\|\alpha\|_2^2}, \frac{t}{\kappa^2\|\alpha\|_2}\right\}\right).;
	\end{split}
	\end{equation}
	Here, $H_{\lambda}(\cdot)$ is the soft-thresholding estimator at level $\lambda$.
	\end{enumerate}
\end{Lemma}
{\bf\noindent Proof of Lemma \ref{lm:X-random-matrix-properties}.} 
\begin{enumerate}[leftmargin=*]
	\item The first statement is via $\varepsilon$-net. Denote $W_T = X_T\Sigma_{T, T}^{-1/2}$, then the rows of $W_T$ are independent isotropic sub-Gaussian distributed.
	For any fixed vector $x \in S^{s - 1} = \{x: x \in \bbR^s, \|x\|_2 = 1\}$, by \cite[Lemma 5.5]{vershynin2010introduction}, $Z_i = \langle (W_{T})_{i\cdot}^\top, x\rangle$ are independent sub-Gaussian random variables with $\mathbb{E}Z_i^2 = 1$ and $\|Z_i\|_{\psi_2} \leq C\kappa$. Therefore, by Remark 5.18 and Lemma 5.14 in \cite{vershynin2010introduction},$\|Z_i^2 - 1\|_{\psi_1} \leq 2\|Z_ i^2\|_{\psi_1} \leq 4\|Z_i\|_{\psi_2}^2 \leq C\kappa^2$. Bernstein-type inequality shows that
	\begin{equation*}
		\bbP\left(\left|\frac{1}{n}\|W_Tx\|_2^2 - 1\right| \geq \frac{t}{2}\right) = \bbP\left(\left|\frac{1}{n}\sum_{i = 1}^{n}(Z_i^2 - 1)\right| \geq \frac{t}{2}\right) \leq 2\exp\left(-cn\min\left\{\frac{t^2}{\kappa^4}, \frac{t}{\kappa^2}\right\}\right).
	\end{equation*}
	
	By \cite[Lemma 5.2]{vershynin2010introduction}, we can find a $\frac{1}{4}$-net $\mathcal{N}_{\frac{1}{4}}$ of $S^{s - 1} = \{x: x \in \bbR^s, \|x\|_2 = 1\}$ with $|\mathcal{N}_{\frac{1}{4}}| \leq 9^s$. The union bound tells us
	\begin{equation}\label{eq1}
		\bbP\left(\max_{x \in \mathcal{N}_{\frac{1}{4}}}\left|\frac{1}{n}\|W_Tx\|_2^2 - 1\right| \geq \frac{t}{2}\right) \leq 9^s\cdot 2\exp\left(-cn\min\left\{\frac{t^2}{\kappa^4}, \frac{t}{\kappa^2}\right\}\right).
	\end{equation}
	By \cite[Lemma 5.4]{vershynin2010introduction},
	\begin{equation}\label{eq2}
	    \|\frac{1}{n}W_T^\top W_T - I_s\| \leq 2\max_{x \in \mathcal{N}_{\frac{1}{4}}}\left|\langle(\frac{1}{n}W_T^\top W_T - I_s)x, x\rangle\right| = 2\max_{x \in \mathcal{N}_{\frac{1}{4}}}\left|\frac{1}{n}\|W_Tx\|_2^2 - 1\right|.
	\end{equation}
	Since $c_{\min} \leq \sigma_{\min}(\Sigma) \leq \sigma_{\max}(\Sigma) \leq C_{\max}$, we have $\|\Sigma_{T, T}^{1/2}\| \leq \sqrt{C_{\max}}$ and $\|\Sigma_{T, T}^{-1/2}\| \leq 1/\sqrt{c_{\min}}$. Therefore,
	\begin{equation}\label{eq23}
		\begin{split}
		\|\frac{1}{n}X_T^\top X_T\Sigma_{T, T}^{-1} - I_s\| =& \|\Sigma_{T, T}^{1/2}(\frac{1}{n}W_T^\top W_T - I_s)\Sigma_{T, T}^{-1/2}\|\\ \leq& \|\Sigma_{T, T}^{1/2}\|\|\frac{1}{n}W_T^\top W_T - I_s\|\|\Sigma_{T, T}^{-1/2}\|\\ \leq& \sqrt{\frac{C_{\max}}{c_{\min}}}\|\frac{1}{n}W_T^\top W_T - I_s\|.
		\end{split}
	\end{equation}
    Combine \eqref{eq1}, \eqref{eq2} and \eqref{eq23} together, we have arrived at the conclusion.

\item Now we consider the proof for \eqref{ineq:H_lambda-tail-bound}. Note that $\|H_{\delta}(\alpha^\top X_T^\top X_\Omega)\|_2\geq t$ implies that there exists $\Lambda\subset \Omega$ such that all entry of $|\alpha^\top X_T^\top X_{\Lambda}|$ are greater than $\delta$, and $\left\||\alpha^\top X_T^\top X_{\Lambda}| - \delta\right\|_2 \geq t$. Thus,
\begin{equation}\label{eq32}
\begin{split}
& \bbP\left(\left\|H_\delta(\alpha^\top X_T^\top X_\Omega/n)\right\|_2 \geq t\right)\\
\leq & \bbP\left(\exists \Lambda \subseteq \Omega, \text{ such that all entries of } \left|\alpha^\top X_T^\top X_{\Lambda}/n\right| \geq \delta, \text{and } \|\alpha^\top X_T^\top X_{\Lambda}/n\|_2 \geq t\right)\\
\leq & \bbP\left(\exists \Lambda \subseteq \Omega, \sqrt{|\Lambda|}\delta\leq t, \|\alpha^\top X_T^\top X_\Lambda/n\|_2\geq t\right)\\
 & + \bbP\left(\exists \Lambda \subseteq \Omega, \sqrt{|\Lambda|}\delta > t, \text{ all entries of } \left|\alpha^\top X_T^\top X_{\Lambda}/n\right| \geq \delta\right)\\
\leq & \sum_{\substack{\Lambda\subseteq \Omega\\|\Lambda| = \lfloor (t/\delta)^2 \rfloor\wedge r}} \bbP\left(\|\alpha^\top X_T^\top X_{\Lambda}/n\|_2\geq t\right) + \sum_{\substack{\Lambda\subseteq \Omega\\|\Lambda| = \lceil (t/\delta)^2 \rceil}} \bbP\left( \text{all entries of } \left|\alpha^\top X_T^\top X_{\Lambda}/n\right| \geq \delta \right)\\
\leq & \sum_{\substack{\Lambda\subseteq \Omega\\|\Lambda| = \lfloor (t/\delta)^2 \rfloor\wedge r}} \bbP\left(\|\alpha^\top X_T^\top X_{\Lambda}/n\|_2\geq t\right) + \sum_{\substack{\Lambda\subseteq \Omega\\|\Lambda| = \lceil (t/\delta)^2 \rceil}} \bbP\left(\left\|\alpha^\top X_T^\top X_{\Lambda}/n\right\|_2 \geq t \right).
\end{split}
\end{equation}
Since $t \geq \delta \geq C\max_{i \in T^c}\|\Sigma_{i, T}\Sigma_{T, T}^{-1}\|_2\|\alpha\|_2$, we know that no matter $|\Lambda| = \lfloor (t/\delta)^2 \rfloor\wedge r$ or $\lceil (t/\delta)^2 \rceil$,
\begin{equation*}
	2C_{\max}\sqrt{\lceil (t/\delta)^2\rceil}(\max_{i \in T^c}\|\Sigma_{i, T}\Sigma_{T, T}^{-1}\|_2)\|\alpha\|_2 \leq 2C_{\max}\sqrt{2}(t/\delta)(\max_{i \in T^c}\|\Sigma_{i, T}\Sigma_{T, T}^{-1}\|_2)\|\alpha\|_2 \leq t.
\end{equation*}
By Part 3 of Lemma \ref{lm:sub-Gaussian-concentration}, for any $\Lambda \subseteq \Omega$, $t \geq 2C_{\max}\sqrt{|\Lambda|}\max_{i \in T^c}\|\Sigma_{i, T}\Sigma_{T, T}^{-1}\|_2\|\alpha\|_2$, we have
\begin{equation*}
    \begin{split}
    &\bbP\left(\left\|\alpha^\top X_T^\top X_{\Lambda}/n\right\|_2 \geq t \right)\\ \leq& \bbP\left(\left\|\alpha^\top X_T^\top X_{\Lambda}/n - \mathbb{E}\alpha^\top X_T^\top X_{\Lambda}/n\right\|_2 \geq t - \left\| \mathbb{E}\alpha^\top X_T^\top X_{\Lambda}/n\right\|_2\right)\\
    \leq& \bbP\left(\left\|\alpha^\top X_T^\top X_{\Lambda}/n - \mathbb{E}\alpha^\top X_T^\top X_{\Lambda}/n\right\|_2 \geq t - \left\| \Sigma_{\Lambda, T}\alpha\right\|_2\right)\\
    =& \bbP\left(\left\|\alpha^\top X_T^\top X_{\Lambda}/n - \mathbb{E}\alpha^\top X_T^\top X_{\Lambda}/n\right\|_2 \geq t - \left(\sum_{i \in \Lambda}(\Sigma_{i, T}\alpha)^2\right)^{1/2}\right)\\
    \leq& \bbP\left(\left\|\alpha^\top X_T^\top X_{\Lambda}/n - \mathbb{E}\alpha^\top X_T^\top X_{\Lambda}/n\right\|_2 \geq t - \sqrt{|\Lambda|}\max_{i \in T^c}|\Sigma_{i, T}\alpha|\right)\\
    \leq& \bbP\left(\left\|\alpha^\top X_T^\top X_{\Lambda}/n - \mathbb{E}\alpha^\top X_T^\top X_{\Lambda}/n\right\|_2 \geq t - \sqrt{|\Lambda|}\max_{i \in T^c}\|\Sigma_{i, T}\Sigma_{T, T}^{-1}\|_2\|\Sigma_{T, T}\|\|\alpha\|_2\right)\\
    \leq& \bbP\left(\left\|\alpha^\top X_T^\top X_{\Lambda}/n - \mathbb{E}\alpha^\top X_T^\top X_{\Lambda}/n\right\|_2 \geq t/2\right)\\
    \leq& 2\exp\left(C|\Lambda| - cn\min\left\{\frac{t^2}{\kappa^4\|\alpha\|_2^2}, \frac{t}{\kappa^2\|\alpha\|_2}\right\}\right).
    \end{split}
\end{equation*}
Combine \eqref{eq32} and the previous inequality, one obtains
\begin{equation*}
\begin{split}
& \bbP\left(\left\|H_\delta(\alpha^\top X_T^\top X_\Omega/n)\right\|_2 \geq t\right)\\
\leq & \binom{r}{\lfloor(t/\delta)^2\rfloor\wedge r}\exp\left(C\lfloor (t/\delta)^2 \rfloor\wedge r - cn\min\left\{\frac{t^2}{\kappa^4\|\alpha\|_2^2}, \frac{t}{\kappa^2\|\alpha\|_2}\right\}\right)\\
& + \binom{r}{\lceil(t/\delta)^2\rceil}_+\exp\left(C\lceil (t/\delta)^2 \rceil - cn\min\left\{\frac{t^2}{\kappa^4\|\alpha\|_2^2}, \frac{t}{\kappa^2\|\alpha\|_2}\right\}\right).
\end{split}
\end{equation*}
\end{enumerate}
\quad $\square$

\begin{Lemma}[Properties of Soft-thresholding]\label{lm:soft-thresholding}
	\begin{enumerate}[leftmargin = *]
		\item Suppose $a, b>0$, $x, y\in \mathbb{R}$, $H_{\cdot}(\cdot)$ is the soft-thresholding operator satisfying $H_a(x) = \sgn(x)\cdot (|x|-a)_+$. Then the following triangular inequality holds,
		\begin{equation}
		|H_{a+b}(x+y)| \leq |H_a(x)| + |H_b(y)|.
		\end{equation}
		\item Suppose $a, b > 0$, $x, y \in \bbR^p$, if $\|H_{a}(x)\|_{\infty, 2} \leq b$, then 
		\begin{equation}
		|\langle x, y\rangle| \leq a\|y\|_1 + b\|y\|_{1, 2}.
		\end{equation}
	\end{enumerate}
\end{Lemma}

{\bf\noindent Proof of Lemma \ref{lm:soft-thresholding}.} 
\begin{enumerate}[leftmargin = *]
	\item 
	\begin{equation*}
	\begin{split}
	|H_{a+b}(x+y)| = & (|x+y|-a-b)_+ \leq (|x|-a + |y|-b)_+ \leq (|x|-a)_+ + (|y|-b)_+\\
	= & |H_a(x)| + |H_b(y)|.
	\end{split}
	\end{equation*}
	\item 
	\begin{equation*}
        \begin{split}
        |\langle x, y\rangle| \leq& |\langle H_{a}(x), y\rangle| + |\langle x - H_{a}(x), y\rangle| = |\sum_{j = 1}^{d}\langle [H_a(x)]_{(j)}, y_{(j)}\rangle| + |\langle x - H_{a}(x), y\rangle|\\ \leq& \sum_{j = 1}^{d}\|[H_a(x)]_{(j)}\|_2\|y_{(j)}\|_2  + \|x - H_{a}(x)\|_{\infty}\|y\|_1 \leq \|H_a(x)\|_{\infty, 2}\|y\|_{1, 2} + \|x - a\|y\|_1\\
        \leq& b\|y\|_{1, 2} + a\|y\|_1.         
        \end{split}
	\end{equation*}
\end{enumerate}

\quad $\square$

\begin{Lemma}\label{lm:infinity norm}
	Suppose $X = [X_1^\top,\ldots, X_n^\top]^\top \in \mathbb{R}^{n\times p}$ is a random matrix with independent random sub-Gaussian rows satisfying Assumption \ref{as:design}, $\varepsilon_i \stackrel{i.i.d.}{\sim} N(0, \sigma^2)$.	Suppose $T \subseteq \{1, \dots, p\}$ is with cardinality $s$, $P \in \bbR^{n \times n}$ is a projection matrix and independent of $X_T$. Then, for any $t \geq \log(es)$,
	\begin{equation*}
		\bbP\left(\|X_T^\top P\varepsilon\|_{\infty} \geq C\kappa\sqrt{nt}\sigma^2\right) \leq e^{-n} + e^{-Ct}. 
	\end{equation*}
\end{Lemma}

{\bf\noindent Proof of Lemma \ref{lm:infinity norm}.} 
For fixed vector $w \in \bbR^n$, since Assumption 2 is satisfied, for $i \in T$, $X_{1i}, \dots, X_{ni}$ are independent sub-Gaussian distributed such that 
\begin{equation*}
	\mathbb{E}\exp\left(tX_{ji}\right) = \mathbb{E}\exp\left(te_i^\top \Sigma^{1/2}\Sigma^{-1/2}X_{j\cdot}^\top\right) \leq \exp\left(\frac{\kappa^2\|\Sigma^{1/2}e_i\|_2^2t^2}{2}\right) \leq \exp\left(\frac{\kappa^2\Sigma_{i, i}t^2}{2}\right) \leq \exp\left(\frac{C_{\max}\kappa^2 t^2}{2}\right).
\end{equation*}
By Hoeffding-type inequality, 
\begin{equation}\label{eq5}
	\bbP\left(|X_{\cdot i}^\top w| \geq t\|w\|_2\right) \leq 2\exp\left(-c\frac{t^2}{\kappa^2}\right).
\end{equation}
Moreover, by \cite[Lemma 1]{laurent2000adaptive}, for any $x \geq 0$,
\begin{equation*}
	\bbP\left(\sum_{i = 1}^{n}\varepsilon_i^2 \geq (n + 2\sqrt{nx} + 2x)\sigma^2\right) \leq e^{-x}.
\end{equation*}
Set $x = n$ in the last inequality, we have
\begin{equation}\label{guassian tail}
    \bbP\left(\|\varepsilon\|_2 \geq \sqrt{5n\sigma^2}\right) \leq e^{-n}.
\end{equation}
Combine \eqref{eq5} and \eqref{guassian tail} together and notice that $\|P\varepsilon\|_2 \leq \|\varepsilon\|_2$, we have
\begin{equation*}
	\begin{split}
	&\bbP\left(\|X_T^\top P\varepsilon\|_\infty \geq C\kappa\sqrt{ nt\sigma^2}\right) \leq \sum_{i \in T}\bbP\left(|X_{\cdot i}^\top P\varepsilon| \geq C\kappa\sqrt{ nt\sigma^2}\right)\\ \leq& \bbP\left(\|P\varepsilon\|_2 \geq \sqrt{5n\sigma^2}\right) + \sum_{i \in T}\bbP\left(|X_{\cdot i}^\top P\varepsilon| \geq C\kappa\sqrt{ nt\sigma^2}, \|P\varepsilon\|_2 \leq \sqrt{5n\sigma^2}\right)\\
	\leq& \bbP\left(\|\varepsilon\|_2 \geq \sqrt{5n\sigma^2}\right) + \sum_{i \in T}\bbP\left(|X_{\cdot i}^\top P\varepsilon| \geq C\kappa\sqrt{ nt\sigma^2}\bigg|\|P\varepsilon\|_2 \leq \sqrt{5n\sigma^2}\right)\\
	\leq& e^{-n} + s\cdot 2\exp\left(-Ct\right)
	\leq e^{-n} + e^{-Ct}.
	\end{split}
\end{equation*}
\quad $\square$

\begin{Lemma}\label{lm:approximate dual certificate}
	With probability at least $1 - Ce^{-cn/s}$, the approximate dual certificate defined in \eqref{eq:approximate dual certificate} can be written as $u = X^\top w$, where $\|w\|_2 \leq C\sqrt{s/n}$.
\end{Lemma}
{\bf\noindent Proof of Lemma \ref{lm:approximate dual certificate}.} 
By \eqref{eq:approximate dual certificate}, we have $u = X^\top w$, where $w = (w_1^\top, \dots, w_{l_{\max}}^\top)^\top$ and $w_l = \frac{1}{n_l}X_{I_l, T}\Sigma_{T, T}^{-1}q_{l - 1}$.
Thus $\|w\|_2^2 = \sum_{l = 1}^{l_{\max}}\|w_l\|_2^2$.
Also note that
\begin{equation*}
	\begin{split}
	\frac{1}{n_l}\|X_{I_l, T}\Sigma_{T, T}^{-1}q_{l-1}\|_2^2 =& \langle \frac{1}{n_l}X_{I_l, T}^\top X_{I_l, T}\Sigma_{T, T}^{-1}q_{l-1}, \Sigma_{T, T}^{-1}q_{l-1}\rangle\\
	=& \langle (\frac{1}{n_l}X_{I_l, T}^\top X_{I_l, T}\Sigma_{T, T}^{-1} - I_{|T|})q_{l-1}, \Sigma_{T, T}^{-1}q_{l-1}\rangle + \|\Sigma_{T, T}^{-1/2}q_{l-1}\|_2^2\\
	=& \langle -q_{l}, \Sigma_{T, T}^{-1}q_{l-1}\rangle + \|\Sigma_{T, T}^{-1/2}q_{l-1}\|_2^2\\
	\leq& \|q_l\|_2\|\Sigma_{T, T}^{-1}q_{l-1}\|_2 + \|\Sigma_{T, T}^{-1/2}q_{l-1}\|_2^2\\
	\leq& \frac{1}{c_{\min}}\|q_l\|_2\|q_{l-1}\|_2 + \frac{1}{c_{\min}}\|q_{l-1}\|_2^2 \leq \frac{2}{c_{\min}}\|q_{l-1}\|_2^2.
	\end{split}
\end{equation*}
By \eqref{ineq:q-l-2-norm}, with probability at least $1 - C\exp\left(-cn/s\right)$,
\begin{equation*}
	\begin{split}
	\|w\|_2^2 \leq& \sum_{l = 1}^{l_{\max}}\frac{C}{n_l}\|q_{l-1}\|_2^2\\ \leq& \frac{C}{n}(2\sqrt{s})^2 + \frac{C}{n}\left(2\sqrt{s/\log(es)}\right)^2 + \frac{C\log(es)}{n}\sum_{l = 3}^{l_{\max}}\left(2^{4 - l}\sqrt{s}/\log(es)\right)^2\\
	\leq& C\frac{s}{n}.
	\end{split}
\end{equation*}
\quad $\square$

{\bf\noindent Proof of Lemma \ref{lm:optimization}.} 
For any $1 \leq i \leq p$, $1 \leq j \leq d,, \Lambda \subseteq (j), |\Lambda| = k$, by Lemma \ref{lm:sub-Gaussian-concentration} with 
\begin{equation*}
	v = \Sigma^{-1}e_i, \quad U \in \bbR^{p \times k}, U_{[\Lambda, :]} = I, U_{[\Lambda^c, :]} = 0, 
\end{equation*}
we have
\begin{equation}\label{eq28}
	\bbP\left(\left\|\left(e_i\right)_{\Lambda} - \frac{1}{n}X_{\Lambda}^\top X\Sigma^{-1}e_i\right\|_2 \geq t\right) \leq 2\exp\left(Ck - cn\min\left\{\frac{t^2}{\kappa^4}, \frac{t}{\kappa^2}\right\}\right).
\end{equation}
By the same method in Lemma \ref{lm:bernstein-sub-gaussian} Part 2, 
\begin{equation*}
	\begin{split}
	&\bbP\left(\left\|H_{\alpha}\left((e_i)_{(j)} - \frac{1}{n}X_{(j)}^\top X\Sigma^{-1}e_i\right)\right\|_{2} \geq \gamma\right)\\
	\leq& \bbP\bigg(\exists \Lambda \subseteq (j), \text{ all entries of } |\left(e_i\right)_{\Lambda} - \frac{1}{n}X_{\Lambda}^\top X\Sigma^{-1}e_i| \geq \alpha \text{ and } \|\left(e_i\right)_{\Lambda} - \frac{1}{n}X_{\Lambda}^\top X\Sigma^{-1}e_i\|_2 \geq \gamma\bigg)\\
	\leq& \bbP\left(\exists \Lambda \subseteq (j), \sqrt{|\Lambda|}\alpha \leq \gamma, \|\left(e_i\right)_{\Lambda} - \frac{1}{n}X_{\Lambda}^\top X\Sigma^{-1}e_i\|_2 \geq \gamma\right)\\
	&+ \bbP\left(\exists \Lambda \subseteq (j), \sqrt{|\Lambda|}\alpha > \gamma, \text{ all entries of } |\left(e_i\right)_{\Lambda} - \frac{1}{n}X_{\Lambda}^\top X\Sigma^{-1}e_i| \geq \alpha\right)\\
	\leq& \sum_{\substack{\Lambda \subseteq (j)\\ |\Lambda| = \lfloor s/s_g\rfloor}}\bbP\left(\|\left(e_i\right)_{\Lambda} - \frac{1}{n}X_{\Lambda}^\top X\Sigma^{-1}e_i\|_2 \geq \gamma\right)\\ &+ \sum_{\substack{\Lambda \subseteq (j)\\ |\Lambda| = \lceil s/s_g\rceil}}\bbP\left(\text{ all entries of } |\left(e_i\right)_{\Lambda} - \frac{1}{n}X_{\Lambda}^\top X\Sigma^{-1}e_i| \geq \alpha\right)\\
	\leq&  \sum_{\substack{\Lambda \subseteq (j)\\ |\Lambda| = \lfloor s/s_g\rfloor}}\bbP\left(\|\left(e_i\right)_{\Lambda} - \frac{1}{n}X_{\Lambda}^\top X\Sigma^{-1}e_i\|_2 \geq \gamma\right) + \sum_{\substack{\Lambda \subseteq (j)\\ |\Lambda| = \lceil s/s_g\rceil}}\bbP\left(\|\left(e_i\right)_{\Lambda} - \frac{1}{n}X_{\Lambda}^\top X\Sigma^{-1}e_i\|_2 \geq \gamma\right).
	\end{split}
\end{equation*}
Combine \eqref{eq28} and the previous inequality together, we have
\begin{equation}\label{eq29}
	\begin{split}
	&\bbP\left(\left\|H_{\alpha}\left((e_i)_{(j)} - \frac{1}{n}X_{(j)}^\top X\Sigma^{-1}e_i\right)\right\|_{2} \geq \gamma\right)\\ \leq& \binom{b_j}{\lfloor s/s_g\rfloor}\cdot 2\exp\left(C\lfloor s/s_g\rfloor - cn\cdot C\frac{s\log(es_gb) + s_g\log(d/s_g)}{s_gn}\right)\\
	&+ \binom{b_j}{\lceil s/s_g\rceil}\cdot 2\exp\left(C\lceil s/s_g\rceil - cn\cdot C\frac{s\log(es_gb) + s_g\log(d/s_g)}{s_gn}\right)\\
	\leq& 4\left(\frac{2es_gb}{s}\right)^{2s/s_g}\exp\left(Cs/s_g - C\frac{s\log(es_gb) + s_g\log(d/s_g)}{s_g}\right)\\
	\leq& 4\exp\left(\frac{2s}{s_g}\log\left(\frac{2es_gb}{s}\right) + Cs/s_g - C\frac{s\log(es_gb) + s_g\log(d/s_g)}{s_g}\right).
	\end{split}
\end{equation}
By \eqref{eq29} and the union bound, we have
\begin{equation*}
	\begin{split}
	&\bbP\left(\max_{1 \leq i \leq p}\|H_{\alpha}(e_i - \frac{1}{n}X^\top X\Sigma^{-1}e_i)\|_{\infty, 2} \leq \gamma\right)\\
	\leq& \sum_{i = 1}^{p}\sum_{j = 1}^{d}\bbP\left(\left\|H_{\alpha}\left((e_i)_{(j)} - \frac{1}{n}X_{(j)}^\top X\Sigma^{-1}e_i\right)\right\|_{2} \geq \gamma\right)\\
	\leq& d^2b\cdot 4\exp\left(\frac{2s}{s_g}\log\left(\frac{2es_gb}{s}\right) + Cs/s_g - C\frac{s\log(es_gb) + s_g\log(d/s_g)}{s_g}\right)\\
	\leq& 4\exp\left(2\log(s_g) + 2\log(d/s_g) + \frac{3s}{s_g}\log\left(2eb\right) + Cs/s_g - C\frac{s\log(es_gb) + s_g\log(d/s_g)}{s_g}\right)\\
	\leq& 4\exp\left(- C\frac{s\log(es_gb) + s_g\log(d/s_g)}{s_g}\right).
	\end{split}
\end{equation*}
\quad $\square$
\end{document}